\documentclass[11pt]{amsart}

\usepackage[margin=1.15in]{geometry}
\usepackage{amsmath,amsfonts,amssymb,amscd,amsthm,amsbsy,epsf,graphics} 
\usepackage{stmaryrd}
\usepackage{setspace}
\usepackage{comment}
\usepackage{color}
\usepackage[dvipsnames]{xcolor}
\usepackage[colorlinks,linkcolor=blue,citecolor=blue]{hyperref} 
\usepackage{epsfig} 
\usepackage{tikz-cd}
\usepackage{bm}
\usepackage{graphicx}  
\usepackage{caption, subcaption}
\usepackage{lipsum}
\usepackage{enumitem} 

\usepackage{array} 
\usepackage[square,numbers]{natbib}

\makeatletter
\def\l@subsection{\@tocline{2}{0pt}{4pc}{5pc}{}}
\makeatother
 
\let\oldtocsection=\tocsection
\let\oldtocsubsection=\tocsubsection
\let\oldtocsubsubsection=\tocsubsubsection
 
\renewcommand{\tocsection}[2]{\hspace{0em}\oldtocsection{#1}{#2}}
\renewcommand{\tocsubsection}[2]{\hspace{0em}\oldtocsubsection{#1}{#2}}
\renewcommand{\tocsubsubsection}[2]{\hspace{2em}\oldtocsubsubsection{#1}{#2}}

\newtheorem{thm}{Theorem}[section]
\newtheorem*{thm*}{Theorem}
\newtheorem{cor}[thm]{Corollary}
\newtheorem{lemma}[thm]{Lemma}
\newtheorem{prop}[thm]{Proposition}
\newtheorem{conj}[thm]{Conjecture}

\newtheorem{innercustomthm}{Theorem}
\newenvironment{customthm}[1]
  {\renewcommand\theinnercustomthm{#1}\innercustomthm}
  {\endinnercustomthm}

\newtheorem{innercustomconj}{Conjecture}
\newenvironment{customconj}[1]
{\renewcommand\theinnercustomconj{#1}\innercustomconj}
{\endinnercustomconj}

\theoremstyle{definition}
\newtheorem{question}[thm]{Question}

\newtheorem{problem}[thm]{Problem}
\newtheorem{remark}[thm]{Remark} 
\newtheorem{remarks}[thm]{Remarks}

\newtheorem*{acknowledgement*}{Acknowledgements}

\newtheoremstyle{cases}
  {12pt plus 6 pt}
  {2pt}
  {\bfseries}   
  {}
  {\bfseries}
  {.}
  {.5em}
  {}
\theoremstyle{cases}

\newtheorem{case}{Case}

\numberwithin{subcase}{case} 
\numberwithin{subsubcase}{subcase}
\numberwithin{equation}{subsection} 

\graphicspath{{figures/}} 

\title{Cyclic branched covers of Seifert links and properties related to the $ADE$ link conjecture}

\author[Steven Boyer]{Steven Boyer} 
\thanks{Steven Boyer was partially supported by NSERC grant RGPIN 9446-2008}
\address{D\'epartement de Math\'ematiques, Universit\'e du Qu\'ebec \`a Montr\'eal, 201 President Kennedy Avenue, Montr\'eal, Qc., Canada H2X 3Y7.}
\email{boyer.steven@uqam.ca}
\urladdr{http://www.cirget.uqam.ca/boyer/boyer.html}

\author[Cameron McA. Gordon]{Cameron McA. Gordon}
\address{Department of Mathematics, University of Texas at Austin, 1 University Station, Austin, TX 78712, USA.}
\email{gordon@math.utexas.edu}
\urladdr{http://www.ma.utexas.edu/text/webpages/gordon.html}

\author[Ying Hu]{Ying Hu}
\address{Department of Mathematical and Statistical Sciences, University of Nebraska Omaha, 6001 Dodge Street, Omaha, NE 68182-0243, USA.}
\email{yinghu@unomaha.edu}
\urladdr{https://yinghu-math.github.io}



\date{\today}

\begin{document}

\begin{abstract}
In this article we show that all cyclic branched covers of a Seifert link have left-orderable fundamental groups, and therefore admit co-oriented taut foliations and are not $L$-spaces, if and only if it is not an $ADE$ link up to orientation. This leads to a proof of the $ADE$ link conjecture for Seifert links. When $L$ is an $ADE$ link up to orientation, we determine which of its canonical $n$-fold cyclic branched covers $\Sigma_n(L)$ have non-left-orderable fundamental groups. In addition, we give a topological proof of Ishikawa's classification of strongly quasipositive Seifert links and we determine the Seifert links that are definite, resp. have genus zero, resp. have genus equal to its smooth $4$-ball genus, among others. In the last section, we provide a comprehensive survey of the current knowledge and results concerning the $ADE$ link conjecture.

\end{abstract}
\maketitle

\setcounter{tocdepth}{1}
{\footnotesize 
\parskip=.2em
\tableofcontents
}

\section{Introduction} 
\label{sec: intro}
In this paper, by a {\it link} we shall mean a non-split oriented link in $S^3$ that is not the unknot. We do not distinguish between a link and its reverse (the link obtained by reversing the orientations of all its components), or its mirror image. Links fall naturally into three classes: we say that a link $L$ with exterior $X(L)$ is {\it Seifert} if $X(L)$ is a Seifert fibre space, {\it toroidal} if $X(L)$ contains an essential torus, and {\it hyperbolic} if $\mbox{int}(X(L))$ has a complete finite volume hyperbolic structure. Here we will be mainly concerned with Seifert links. The corresponding question for toroidal and hyperbolic links is studied in \cite{BGH1} and \cite{BGH2}, respectively (Also see the discussions in Section \ref{sec: $ADE$ lk conj}). We show that almost all cyclic branched covers of Seifert links have left-orderable fundamental groups. We also determine the Seifert links with various properties, for example, those that have definite symmetrized Seifert form, are braid positive, have genus zero, or have equal genus and 4-ball genus. We use the last to give a purely topological proof of Ishikawa's classification of strongly quasipositive Seifert links \cite{Ish}.

We now give a more detailed description of the results of the paper.

Burde and Murasugi have shown that a Seifert link is either a connected sum of Hopf links or a union of fibres of some Seifert structure on $S^3$  \cite{BuMu}. It is non-split, and prime unless it is a connected sum of Hopf links.

Recall that the $L$-space Conjecture (\cite[Conjecture 1]{BGW}, \cite[Conjecture 5]{Juh}) asserts that for a connected, closed, orientable, irreducible 3-manifold $M$ the following are equivalent: $M$ is not an $L$-space, $\pi_1(M)$ is left-orderable, and $M$ supports a co-orientable taut foliation. We will sometimes abbreviate these by saying that $M$ is $NLS$, $LO$, or $CTF$, respectively. There has been a good deal of work studying these properties for cyclic branched covers of knots. Similarly, one can investigate these properties for cyclic branched covers of links. However, if $n > 2$ then a link with more than one component has several $n$-fold cyclic branched covers. Here is a precise definition. 

Let $L = K_1 \cup ... \cup K_m$ be a link, $\mu_i$ an oriented meridian of $K_i$, $n$ an integer $\ge 2$, and $\psi : \pi_1(X(L)) \to \mathbb Z/n$ an epimorphism such that $\psi(\mu_i) \ne 0$ for each $i$. Then the $n$-fold cyclic cover $X_\psi(L) \to X(L)$ determined by $\psi$ extends to an $n$-fold cyclic cover $\Sigma_\psi(L) \to S^3$ branched over $L$. We say that $\Sigma_\psi(L)$ is a {\it cyclic branched cover of $L$ associated with $\psi$}. The {\it canonical} $n$-fold cyclic branched cover $\Sigma_n(L)$ is the one associated with the epimorphism $\pi_1(X(L)) \to \mathbb{Z}/n$ that sends each $\mu_i$ to $1$ (mod $n$). Note that when $L$ is a knot $K$, $\Sigma_\psi(K) \cong \Sigma_n(K)$. Also, for $n = 2$, $\Sigma_\psi(L) = \Sigma_2(L)$; in particular, $\Sigma_2(L)$ is independent of the orientation of $L$. Finally, we note that the set of all cyclic branched covers of a given link $L$ is independent of the orientation of $L$. An easy consequence of the equivariant sphere theorem is that if a (non-split) link $L$ is prime then any cyclic branched cover of $L$ is irreducible, so the $L$-space Conjecture applies to such manifolds.

The idea of considering the more general branched covers defined above arose from the authors' work on the cyclic branched covers of toroidal and hyperbolic links in \cite{BGH1} and \cite{BGH2}.

One of our main results (Theorem \ref{thm: unoriented $ADE$}) is that, with a few explicit exceptions, the fundamental group of any cyclic branched cover $\Sigma_\psi(L)$ of a prime Seifert link $L$ has a non-trivial representation into $\widetilde{PSL_2}(\mathbb R)$, and is therefore left-orderable \cite{BRW}. Since $\Sigma_\psi(L)$ is Seifert fibred and the $L$-space Conjecture is known for Seifert fiber spaces \cite{BGW, LS}, it follows that $\Sigma_\psi(L)$ is also $CTF$ and $NLS$.
The exceptions are the $ADE$ links (Figure \ref{fig: ade links}), which we now describe.

An {\it $ADE$ link} is the boundary of a plumbing of positive Hopf bands according to the tree determined by one of the Dynkin diagrams $A_m, D_m, E_6, E_7$, or $E_8$.  More precisely, denoting the $(a_1, \ldots, a_n)$ pretzel link by $P(a_1, \ldots, a_n)$ and the $(a, b)$ torus link by $T(a, b)$, the $ADE$ links are given by 
\begin{align*}
    & L(A_m) = P(1, m) = T(2, m+1), \, m \geq 1, \\
    & L(D_m) = P(-2, 2, m-2), \, m \geq 4, \\
    & L(E_6) = P(-2,3,3) = T(3,4), \\
    & L(E_7) = P(-2,3,4), \\
    & L(E_8) = P(-2,3,5) = T(3,5).
\end{align*}
A distinctive feature of the family of $ADE$ links is that the fundamental groups of their $2$-fold cyclic branched covers are precisely the finite subgroups of $SU(2)$ (Proposition \ref{prop: 2-fold}). These branched covers therefore have non-left-orderable fundamental groups, do not support co-orientable taut foliations, and as elliptic manifolds, are Heegaard Floer $L$-spaces.

\begin{thm}
\label{thm: unoriented $ADE$}
Let $L$ be a prime Seifert link that is not an $ADE$ link as an unoriented link. Then, for any cyclic branched cover $\Sigma_\psi(L)$ of $L$, $\pi_1(\Sigma_\psi(L))$ has a non-trivial representation into $\widetilde{PSL_2}(\mathbb R)$. Hence $\Sigma_\psi(L)$ is $LO$, $CTF$, and $NLS$.
\end{thm}
For links which are $ADE$ links up to orientation, we determine precisely when their canonical cyclic branched covers are $LO$, $CTF$, and $NLS$ in Theorem \ref{thm: canonical br covers}.

The $ADE$ links belong to the much more general family of fibred strongly quasipositive links. Recall that a surface in $S^3 = \partial B^4$ is {\it quasipositive} if it is isotopic rel boundary to a surface  $F$, properly embedded in $B^4\subset \mathbb{C}^2$, which is the intersection with $B^4$ of a holomorphic curve in $\mathbb C^2$. A link $L \subset S^3$ is {\it strongly quasipositive} if it (or its mirror image) has a quasipositive Seifert surface. The family of fibred strongly quasipositive links arises naturally as the set of bindings of open books which carry the tight contact structure on the $3$-sphere (\cite{Hed}). Topologically, Giroux's stabilisation theorem characterises the family as the set of fibred links whose fibre surface can be transformed into a plumbing of positive Hopf bands by a finite sequence of such plumbings (\cite{Gi, Ru3}).

  In \cite[Conjecture 1.7]{BBG1} it is conjectured that the $ADE$ links are the only prime, fibred, strongly quasipositive links $L$ for which some canonical $n$-fold cyclic branched cover $\Sigma_n(L)$ of $L$ is an $L$-space. In light of the L-space Conjecture, the present authors' results in \cite{BGH1} and \cite{BGH2}, and Theorem \ref{thm: unoriented $ADE$}, it is natural to consider the following more general version of the conjecture.

\begin{conj}[The $ADE$ Link Conjecture] 
\label{conj: $ADE$  conj}
If $L$ is a prime, fibred, strongly quasipositive link that is not an $ADE$ link then all cyclic branched covers $\Sigma_{\psi}(L)$ of $L$ are $NLS$, $LO$, and $CTF$.
\end{conj}

An almost immediate corollary of Theorem \ref{thm: unoriented $ADE$} is that the $ADE$ Link Conjecture holds for Seifert links.

\begin{thm}[The $ADE$ Link Conjecture for Seifert links]
\label{thm: $ADE$ for seifert links intro}  
If $L$ is a prime, fibred, strongly quasipositive Seifert link that is not an $ADE$ link then all cyclic branched covers $\Sigma_\psi(L)$ of $L$ are $NLS$, $LO$, and $CTF$.
\end{thm}

To deduce Theorem \ref{thm: $ADE$ for seifert links intro} from Theorem \ref{thm: unoriented $ADE$} it remains only to show that a fibred strongly quasipositive link which is an $ADE$ link up to orientation is an $ADE$ link. This is Corollary \ref{cor: $ADE$ unoriented} below. 

 We give a fuller account of the current state of knowledge on the $ADE$ Link Conjecture, and related questions about cyclic branched covers of links in Section \ref{sec: $ADE$ lk conj}. 

We now turn to the other main results of the paper. In Section \ref{sec: dsl results} we describe our general notation for Seifert links, but here we mention one example which arises in the next proposition and in Theorems \ref{thm: sqp seifert results intro} and \ref{thm: def results}. 

Let $p$ and $q$ be positive coprime integers and let $k$ be a positive even integer, and assume that $(p,q,k) \ne (1,1,2)$. Then $L(p,q;k,0)$ denotes the link consisting of $k$ parallel copies of the $(p,q)$-torus knot on the Heegaard torus in $S^3$, with half of the components oriented in one direction and the other half oriented in the opposite direction.  

The following proposition is an easy consequence of known facts. 

\begin{prop}
\label{prop: basic props}
A Seifert link is
\begin{enumerate}[leftmargin=*] 
\setlength\itemsep{0.3em}
\item[{\rm (1)}]non-split;
\item[{\rm (2)}]prime if and only if it is not a connected sum of two or more Hopf links;
\item[{\rm (3)}] fibred if and only if it is not of the form $L(p, q; k, 0)$.
\end{enumerate}
\end{prop}

It is clear that any Seifert link can be reoriented so that it is obviously the closure of a positive braid. Let $\mathcal P$ denote the set of such positive braid Seifert links (see (\ref{equ: P})). (In fact, $\mathcal P$ is precisely the set of positive braid Seifert links; see Corollary \ref{cor: bp}.) Then $L \in \mathcal P$ implies that $L$ is fibred and strongly quasipositive. The converse was proved by Ishikawa [Ish]. More precisely, he proved the following. 

\begin{thm} {\rm (Ishikawa)}
\label{thm: sqp seifert results intro}
Let $L$ be a Seifert link.

$(1)$ $L$ is fibred and strongly quasipositive if and only if $L \in \mathcal P$.

$(2)$ $L$ is strongly quasipositive if and only if either $L \in \mathcal P$ or $L$ is of the form $L(p, q; k, 0)$.
\end{thm}

We give a different proof of Theorem \ref{thm: sqp seifert results intro} in Section \ref{subsec: sqp sl} which does not involve any contact structure computations, but instead is based on classifying the Seifert links for which the genus and 4-ball genus are equal. The following theorem is proved in Section \ref{subsec: g3 = g4}.

\begin{thm}
\label{thm: g_4 = g results}
Let $L$ be a Seifert link. Then $g_{4}(L) = g(L)$ if and only if either $L \in \mathcal P$ or $g(L) = 0$.
\end{thm}

In Theorem \ref{thm: g = 0} we classify the Seifert links of genus zero, and in Proposition \ref{prop: sqp g = 0} we show that the only genus zero links that are strongly quasipositive are the links $L(p,q;k,0)$. Together with Theorem \ref{thm: g_4 = g results} this proves Theorem \ref{thm: sqp seifert results intro}. 

Recall that, denoting the number of components of a link $L$ by $|L|$, its signature $\sigma(L)$ satisfies
$$|\sigma(L)| \le 2g(L) + |L| - 1  $$
We say that $L$ is {\it definite} if this is an equality. 

Definite links arise naturally in the study of the $ADE$ Link Conjecture: a strongly quasipositive link with an $L$-space canonical cyclic branched cover is definite (\cite[Theorem 1.1]{BBG1}), and the $ADE$ links are precisely the prime, definite, positive braid links (\cite{Baa}). More generally, they coincide with the set of prime definite strongly quasipositive links whose Birman-Ko-Lee exponent is $2$ or more (\cite[Theorem 1.5]{BBG2}). They are also precisely the definite, prime, fibred, strongly quasipositive links of braid index 2 or 3; see \cite[Theorem 1.8]{BBG2}.

Since the $4$-ball genus of a definite link coincides with its genus, Theorem \ref{thm: g_4 = g results} implies that a definite Seifert link is either braid positive or has genus zero. This enables us to determine the definite Seifert links. The link $L(2,3;1,1;-,-)$ in (4) below is a 3-component fibred link of genus zero. Such links are classified in \cite{Ro} (see also \cite{EO}); $L(2,3;1,1;-,-)$ is the only one with a knotted component (a trefoil).

 \begin{thm}
\label{thm: def results}
A Seifert link is definite if and only if it is either
\begin{enumerate}[leftmargin=*] 
    \setlength\itemsep{0.3em}
    \item[{\rm (1)}] a connected sum of two or more positive Hopf links, or 
    \item[{\rm (2)}]an $ADE$ link, or
    \item[{\rm (3)}]of the form $L(p,q;2,0)$, or
    \item[{\rm (4)}] $L(2,3; 1,1; -,-)$.
\end{enumerate}
\end{thm}

\begin{cor}
\label{cor: $ADE$ unoriented}
If $L$ is a fibred strongly quasipositive link that is an $ADE$ link up to orientation then $L$ is an $ADE$ link. 
\end{cor}

\begin{proof}
A fibred strongly quasipositive link $L$ which is an $ADE$ link up to orientation is prime and also definite, since $\Sigma_2(L)$ is an $L$-space (\cite[Theorem 1.1]{BBG1}).  Theorem \ref{thm: def results} then implies that $L$ is an $ADE$ link since  
$L(p, q; 2, 0)$ is not fibered by Proposition \ref{prop: basic props} and $L(2, 3; 1, 1; -. -)$ is not strongly quasipositive by Proposition \ref{prop: sqp g = 0}.
\end{proof} 
The proof of Corollary \ref{cor: $ADE$ unoriented} shows that a prime fibred strongly quasipositive Seifert link $L$ is definite if and only if it is an $ADE$ link.  

\begin{cor}
\label{cor: definite FSQP Seifert}
For a prime fibred strongly quasipositive Seifert link $L$, the following are equivalent.
\begin{enumerate}[leftmargin=*] 
\setlength\itemsep{0.3em}
\item[{\rm (1)}] $L$ is an $ADE$ link;
\item[{\rm (2)}] $L$ is definite;
\item[{\rm (3)}] $\Sigma_{2}(L)$ is an $L$-space;
\item[{\rm (4)}] some cyclic branched cover $\Sigma_\psi(L)$ of $L$ is an $L$-space.
\end{enumerate}
\end{cor}

\begin{proof}
The equivalence of (1) and (2) is pointed out at after the proof of Corollary \ref{cor: $ADE$ unoriented}. Also, (1) implies (3), and (3) obviously implies (4). Finally, (4) implies (1) by Theorem \ref{thm: $ADE$ for seifert links intro}.
\end{proof}

\subsection*{Organisation of paper} We discuss the family of $ADE$ links in Section \ref{sec: conventions} and determine when they have a canonical cyclic branched cover with finite fundamental group. In Section \ref{sec: dsl results} we introduce our notation for Seifert links, note which correspond to the $ADE$ links, and prove some basic facts on their topology. Section \ref{sec: def seif links} contains the proof of Theorem \ref{thm: unoriented $ADE$} and the classification of which canonical cyclic branched covers of Seifert links have non-left-orderable fundamental groups. Section \ref{sec: classn sqpsls} is devoted to the proofs of Theorems \ref{thm: sqp seifert results intro}, \ref{thm: g_4 = g results} and \ref{thm: def results}, and Corollary \ref{cor: definite FSQP Seifert}. We end the paper in Section \ref{sec: $ADE$ lk conj} with a survey of the $ADE$ link conjecture.

\begin{acknowledgement*}
The authors thank Jacob White for the preparation of some of the figures. 
\end{acknowledgement*}

\section{\texorpdfstring{$ADE$}{ADE} links and their canonical cyclic branched covers}
\label{sec: conventions} 

It is known when the canonical cyclic branched covers of the $ADE$ links are $L$-spaces (\cite{GLid}, \cite[Example 1.11]{Nem}, \cite[Section 10.1]{BBG2}). In this section we show that the fundamental groups of the $L$-spaces which arise in this way are precisely the finite non-trivial subgroups of $SU(2)$ (Propositions \ref{prop: 2-fold} and \ref{prop: higher order}).

We assume throughout the paper that manifolds are compact, connected, and orientable, unless otherwise stated. The exterior of a link $L$ will be denoted by $X(L)$, though we will sometimes shorten this to $X$ when $L$ is understood. 

A {\it Seifert surface} for a link $L$ is a compact, connected, oriented surface $F$ in $S^3$ with $L = \partial F$ its oriented boundary. A Seifert matrix $S$ of $F$ is a square matrix of size $\beta_1(F)= 2g(F) + |L| - 1$ representing the Seifert form of $F$. Up to units in $\mathbb Z[t, t^{-1}]$, the {\it Alexander polynomial} of $L$ is 
$$\Delta_L(t) = \det(S - tS^T)$$ 

\begin{figure}[ht]
\centering
\begin{subfigure}{0.49\textwidth}
\centering
\includegraphics[scale=0.8]{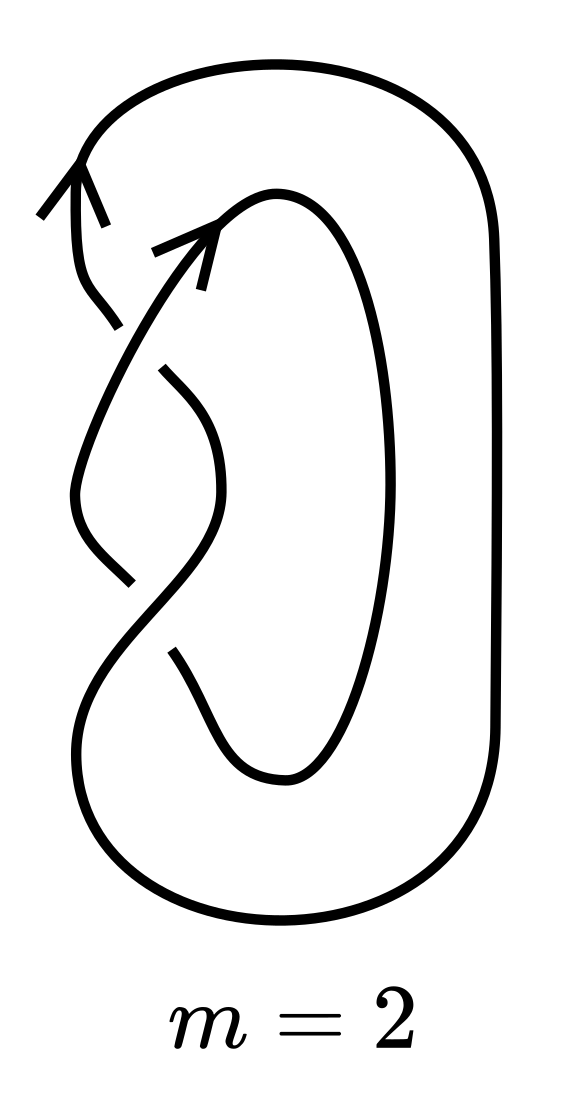}
\caption{$L(A_m) = T(2, m)$}
\end{subfigure}
\hfill
\begin{subfigure}{0.49\textwidth}
\centering
\includegraphics[scale=0.63]{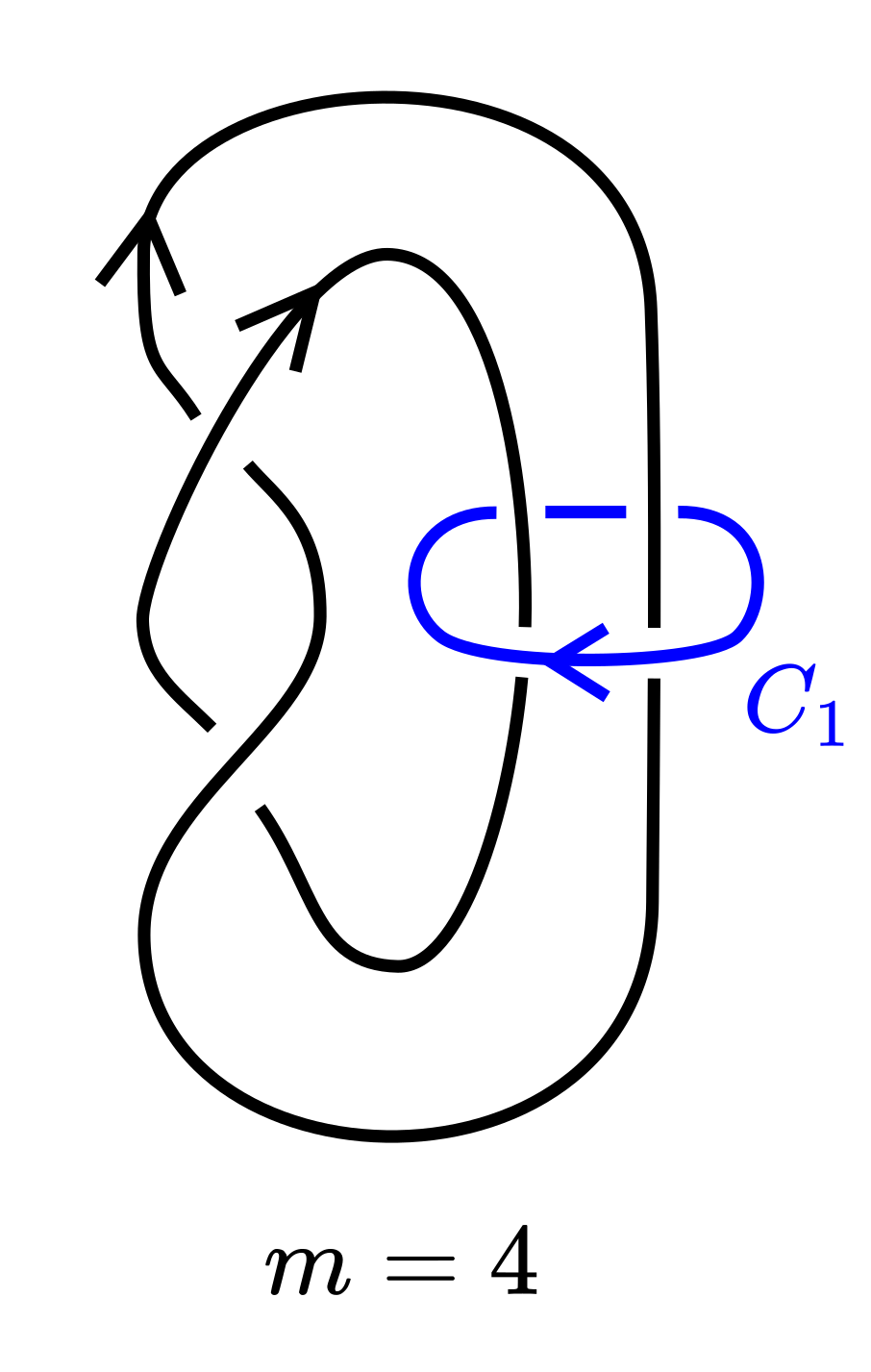}
\caption{$L(D_m) = P(-2, 2, m-2)$, $m \geq 4$} 
\end{subfigure}
\vspace{20pt}
\begin{subfigure}{0.33\textwidth}
\centering 
\vspace{10pt}
\includegraphics[scale=0.63]{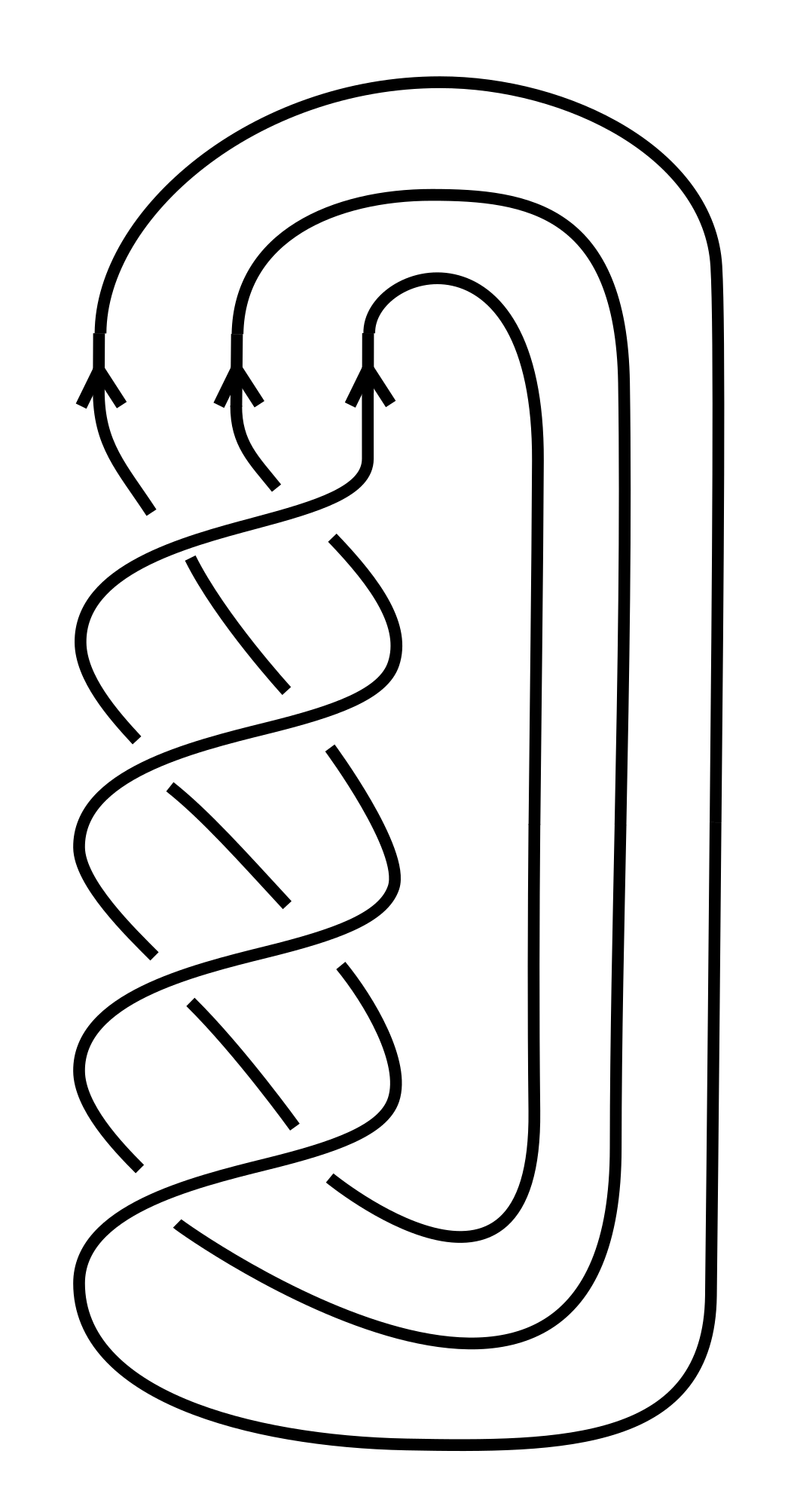}
\caption{$L(E_6) = P(-2,3,3) = T(3,4)$}
\end{subfigure}
\hfill  
\begin{subfigure}{0.31\textwidth}
\centering
\vspace{20pt}
\includegraphics[scale=0.57]{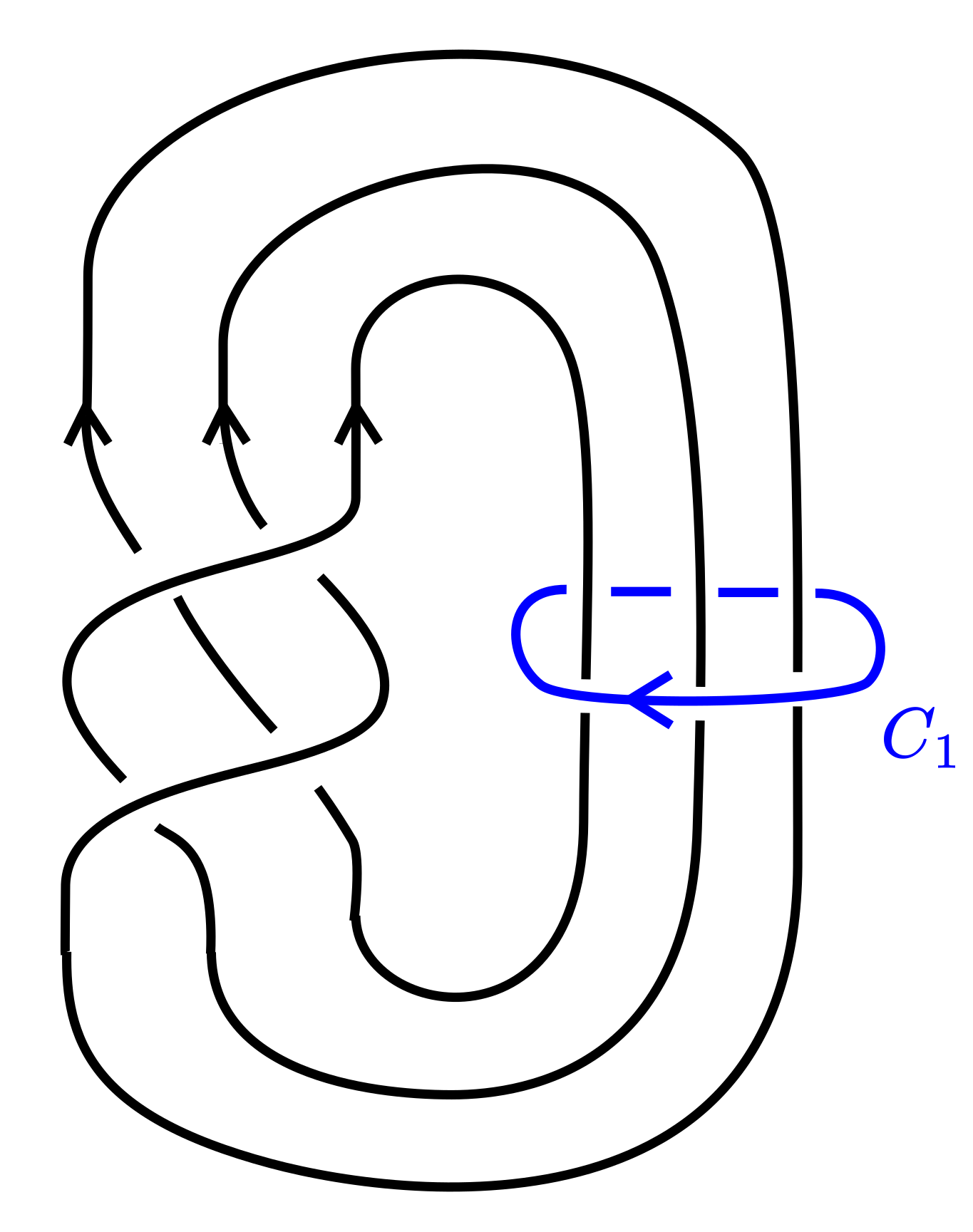}
\caption{$L(E_7) = P(-2,3,4)$}
\end{subfigure}
\hfill
\begin{subfigure}{0.33\textwidth}
\centering
\includegraphics[scale=0.63]{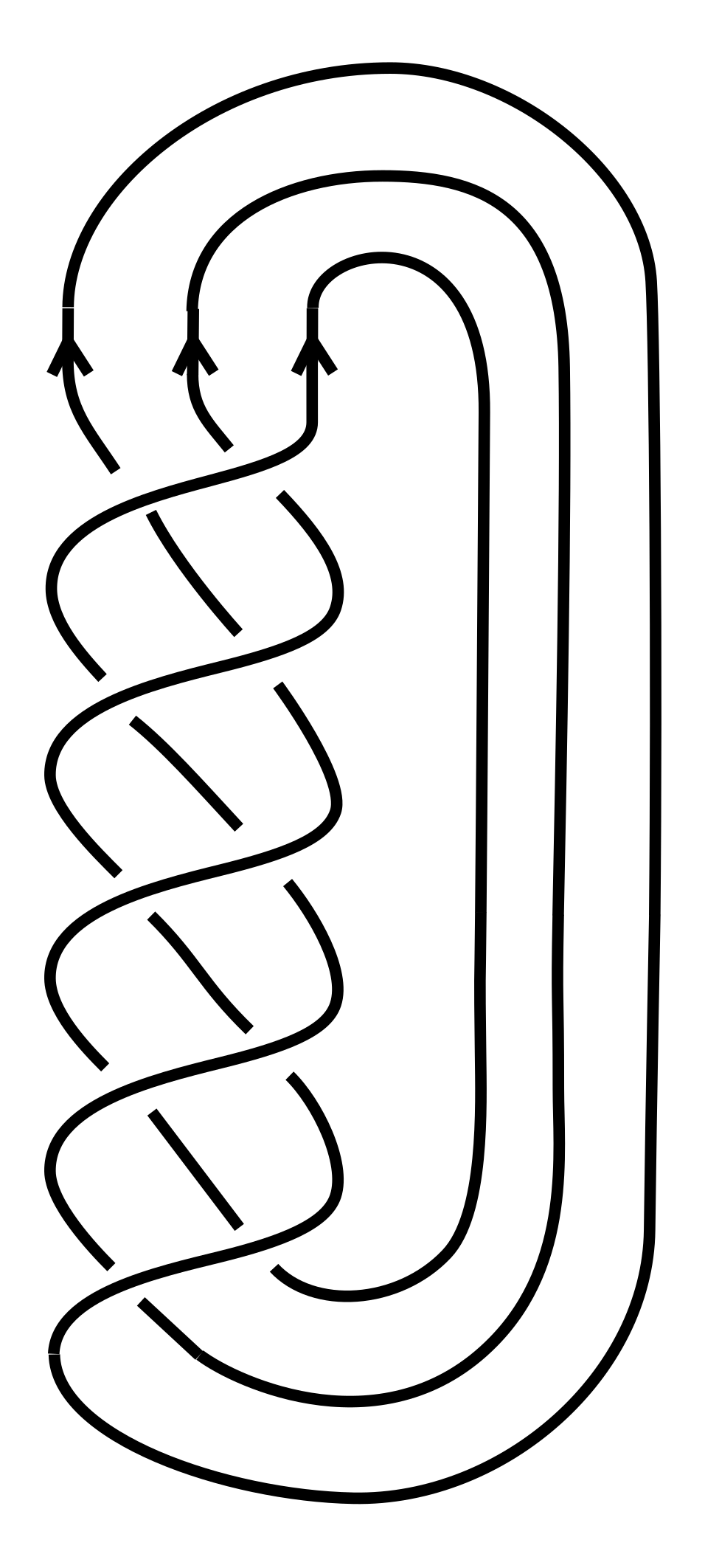}
\caption{$L(E_8) = P(-2,3,5) = T(3,5)$}       
\end{subfigure}
\caption{The $ADE$ links}
\label{fig: ade links}
\end{figure}

Recall that an {\it $ADE$  link} is the boundary of a plumbing of positive Hopf bands according to the tree determined by one of the Dynkin diagrams $A_m, D_m, E_6, E_7$, or $E_8$. These links are illustrated in Figure \ref{fig: ade links}.

It is straightforward to verify that the $ADE$ links are as identified and that the exterior of each is a Seifert fibre space (cf. Proposition \ref{prop: $ADE$ are seifert}). We claim that their $2$-fold cyclic branched covers have finite fundamental groups, and as such are Heegaard Floer $L$-spaces.  

The non-trivial finite subgroups of $SU(2)$ consist of the cyclic group $C_r$ of order $r \geq 2$, the binary dihedral group $D_r^*$ of order $4r \geq 4$, the binary tetrahedral $T^*$ of order $24$, the binary octahedral $O^*$ of order $48$, and the binary icosahedral $I^*$ of order $120$. Famously, they satisfy an $ADE$ classification via the McKay correspondence, and this classification arises naturally in the context of the $ADE$ links. 

\begin{prop}
\label{prop: 2-fold} 
The fundamental groups of the $2$-fold cyclic branched covers of the $ADE$ links are finite and given by the following table.  
\begin{center}
\begin{tabular}{|c||c|c|} \hline 
{\rm Diagram} &  $L$ & $\pi_1(\Sigma_2(L))$ \\  \hline \hline 
$A_m$ $(m \geq 1)$ & $P(1, m) = T(2,m+1)$ & $C_{m+1}$\\  \hline 
$D_m$ $(m \geq 4)$& $P(-2,2,m-2)$ & $D_{m-2}^*$\\  \hline 
$E_6$ &   $P(-2,3,3) = T(3,4)$& $T^*$ \\  \hline 
$E_7$ & $P(-2,3,4)$  & $O^*$ \\  \hline 
$E_8$ & $P(-2,3,5) = T(3,5)$ &  $I^*$ \\  \hline 
\end{tabular}
\end{center}
\end{prop}

\begin{proof}
A closed, connected $3$-manifold has a finite fundamental group $G$ if and only if it is a lens space or admits a Seifert structure with base orbifold $\mathcal{B} = S^2(a,b,c)$, where $a, b, c \geq 2$ is a Platonic triple. Further, 
\begin{itemize}

\item $G \cong D_m^*$ if and only if $\mathcal{B} = S^2(2,2,m)$ and $|H_1(G)| = \left\{ \begin{array}{ll} 2  & \mbox{ if $m$ odd} \\ 4 & \mbox{ if $m$ even} \end{array} \right.$;

\vspace{.2cm} \item $G \cong T^*$ if and only if $\mathcal{B} = S^2(2,3,3)$ and $|H_1(G)| = 3$;

\vspace{.2cm} \item $G \cong O^*$ if and only if $\mathcal{B} = S^2(2,3,4)$ and $|H_1(G)| = 2$;

\vspace{.2cm} \item $G \cong I^*$ if and only if $\mathcal{B} = S^2(2,3,5)$ and $|H_1(G)| = 1$.

\end{itemize}
See Lemma 5.1 of \cite{BZ}. 

For integers $a_i \ne 0$, the $2$-fold cyclic branched cover of the $(a_1, a_2, \ldots, a_n)$-pretzel link admits a Seifert fibre space structure with base orbifold $S^2(|a_1|, |a_2|, \ldots, |a_n|)$, so given our identification of an $ADE$ link $L$ with a pretzel link, $\Sigma_2(L)$ admits a Seifert structure whose base orbifold $\mathcal{O}_2(L)$ is determined. For instance, $L(A_m) = P(1, m)$ so $\Sigma_2(L(A_m))$ has a Seifert structure with base orbifold $\mathcal{O}_2(L(A_m)) \cong S^2(1, m)$. It follows that $\Sigma_2(L(A_m))$ has Heegaard genus one, so has a cyclic fundamental group whose order is 
$$|\pi_1(\Sigma_2(L(A_m)))| = |H_1(\Sigma_2(L(A_m)))| = |\Delta_{L(A_m)}(-1)|$$  
Since $L(A_m) = T(2, m+1)$ has Alexander polynomial is $(t^{m+1} + (-1)^m)/(t+1)$, we conclude that  
$$\pi_1(\Sigma_2(L(A_m))) \cong C_{m+1}$$
Given the first paragraph of the proof, the cases $L(D_m), L(E_6), L(E_7), L(E_8)$ can be dealt with using similar arguments and 
the table \footnote{The Alexander polynomials of these links are calculated in Section \ref{subsec: alex ply}. The notation for the $ADE$ links used in that section can be found in Table \ref{tab: examples of seifert links}.}:  

\begin{center}
\begin{tabular}{|c||c|c|c|} \hline 
$L$ & $\mathcal{O}_2(L)$ & $\Delta_L(t)$ & $|H_1(\Sigma_2(L))| = |\Delta_L(-1)|$ \\  \hline \hline 
&&&  \\
$\begin{array}{c} L(D_m) \\ (m \geq 4) \end{array}$ &$S^2(2,2, m-2)$ &   $\left\{ \begin{array}{cl} t-1  & \mbox{ if $m$ odd} \\ (t-1)(t^{m-1} - 1)& \mbox{ if $m$ even} \end{array} \right.$  & $\left\{ \begin{array}{ll} 2  & \mbox{ if $m$ odd} \\ 4 & \mbox{ if $m$ even} \end{array} \right.$  \\  
&&&  \\ \hline 
$L(E_6)$ &$S^2(2,3,3)$&   $t^6 - t^5 + t^3 - t + 1$  & $3$  \\  \hline 
$L(E_7)$ &$S^2(2,3,4)$&   $(t-1)(1 + t^3 + t^6)$  & $2$  \\  \hline 
$L(E_8)$ &$S^2(2,3,5)$&   $t^8 - t^7 + t^5 - t^4 + t^3 - t + 1$  & $1$  \\  \hline 
\end{tabular}
\end{center}
\end{proof}

A similar proof yields the following proposition.  

\begin{prop}
\label{prop: higher order}
If $L$ is an $ADE$ link, then $\pi_1(\Sigma_n(L))$ is finite for some $n \geq 3$ if and only if $L$ is either $T(2,2), T(2,3), T(2, 4)$, or $T(2,5)$. Further, 
\begin{enumerate}[leftmargin=*] 
\setlength\itemsep{0.3em}
\item[{\rm (a)}] $\pi_1(\Sigma_n(T(2,2)) \cong \mathbb Z/n$ for each $n \geq 3$;
\item[{\rm (b)}] When $m = 3,4,5$, the values of $n \geq 3$ for which $\pi_1(\Sigma_n(T(2, m))$ is finite and the associated groups are listed in the following table.
\begin{center}
\begin{tabular}{|c||c|c|c|} \hline 
$L$ & $\pi_1(\Sigma_3(L))$ & $\pi_1(\Sigma_4(L))$ & $\pi_1(\Sigma_5(L))$ \\  \hline \hline 
$T(2,3)$  & $D_2^*$  & $T^*$ & $I^*$ \\  \hline 
$T(2,4)$ &  $T^*$  & $*$ & $*$\\  \hline 
$T(2,5)$ &  $I^*$  & $*$ & $*$\\  \hline 
\end{tabular}
\end{center}
\end{enumerate}
\end{prop}

We show in Section \ref{sec: def seif links} that if $L$ is an $ADE$ link and $\pi_1(\Sigma_n(L))$ is infinite, then $\pi_1(\Sigma_n(L))$ has a non-trivial representation into $\widetilde{PSL_2}(\mathbb R)$ and is therefore left-orderable. See Proposition \ref{prop: pos bd case}. 

\section{Seifert links}
\label{sec: dsl results}
Consider the genus $1$ Heegaard splitting of $S^3$, $S^3 = V_1 \cup _T V_2$, where $V_1$ and $V_2$ are solid tori and $T = \partial V_1 = \partial V_2$. Let $p$ and $q$ be positive coprime integers and let $C_{p, q}$ be a simple closed curve on $T$ that has algebraic intersection numbers $q$ with a meridian disk of $V_1$ and $p$ with a meridian disk of $V_2$. Let $C_i$ be a core of $V_i$, $1 = 1, 2$.

By \cite{BuMu}, a Seifert link is isotopic (as an unoriented link, and up to mirror image) to one of the following: 
\begin{itemize}
    \setlength\itemsep{0.3em}
\item[(I)] the union of $C_1$ and $k \geq 2$ parallel (with respect to the $0$-framing) copies of $C_2$. This is a connected sum of $k$ Hopf links.
\item[(II)] $k \geq 1$ parallel copies of $C_{p, q}$ on $T$, which we denote by $kC_{p, q}$.
\item[(III)] $kC_{p, q} \cup C_1$, $k\geq 1$.
\item[(IV)] $kC_{p, q} \cup C_1 \cup C_2$, $k\geq 1$.
\end{itemize}

By \cite{JWW}, if $L$ is a Seifert link that is not a Hopf link then the exterior of $L$ has a unique Seifert fibration up to isotopy. 

Taking $p = q = 1$ and $k = 2$ in case II gives the Hopf link, so we exclude $k = 1$ in case I. 

Choose an orientation on the curve $C_{p, q}$ and call it the {\it positive} orientation. In cases II, III, and IV, $L$ will have $k_+$ positively oriented copies of $C_{p, q}$ and $k_-$ negatively oriented copies, where $k = k_+ + k_-$. Since we do not distinguish between a link and its reverse we may assume that $w = k_+ - k_- \geq 0$.
\begin{itemize}
\setlength\itemsep{0.3em}
\item In case II we write $L = L(p, q; k, w)$;  
\item In case III we write $L = L(p, q; k, w; \varepsilon)$, where $\varepsilon = \pm$ so that the linking number of $C_1$ with a positively oriented copy of $C_{p, q}$ is $\varepsilon p$; 
\item In case IV we write $L = L(p, q; k, w; \varepsilon_1, \varepsilon_2)$, where the linking number of a positively oriented copy of $C_{p,q}$ with $C_1$ is $\varepsilon_1 p$, and with $C_2$ is $\varepsilon_2 q$. 
\end{itemize}

Let $H_+ (H_-)$ denote the positive (negative) Hopf link.  We point out the redundancies in the above notation in Remark \ref{rem: notation conventions}.

\begin{remark}[Notational redundancies]
    \label{rem: notation conventions}
\end{remark}
\vspace{-.5cm}
\begin{enumerate}
\setlength\itemsep{0.3em}
\item  In case II we exclude $L(1,1;2,0) = H_-$.  Similarly, we  exclude $L(1,q;1,1;\varepsilon) = H_{\varepsilon}$ from case III.

\item   In case II, and in case IV when $\varepsilon_1  = \varepsilon_2$, the link $L$ is symmetric in $p$ and $q$, so in these cases, we may assume that $p \leq q$.

\item  Since the unknot is excluded by hypothesis, in case II we assume that $p > 1$ in $L(p,q;1,1) = T(p,q)$.

\item  In case III we have $L(p, 1; k, w; \varepsilon) = L(p, 1; k+1, w+ \varepsilon)$ (or the reverse of $L(p, 1; k+1, 1)$ if $w = 0$ and $\varepsilon = -)$, so we may assume $q > 1$. 

\item  Similarly in case IV we may assume $p, q > 1$. Further, we have $L(p,q;k,w;-,+) = L(q,p;k,w;+,-)$.

\item Finally, note that $L(p,q;k,0;+)$ is the reverse of $L(p,q;k,0;-)$ and $L(p,q;k,0;\varepsilon_1,\varepsilon_2)$ is the reverse of $L(p,q;k,0;-\varepsilon_1, -\varepsilon_2)$, and recall that we do not distinguish between a link and its reverse.
\end{enumerate}

In cases II, III, and IV we will say that the link is {\it $0$-core, $1$-core}, or {\it  $2$-core}, respectively.

\begin{prop}
\label{prop: $ADE$  are seifert}
  The $ADE$  links are Seifert links.
\end{prop}

\begin{proof}
The descriptions of the $ADE$ links as Seifert links are:
 
\begin{itemize}
\setlength\itemsep{0.3em}
\item $T(2,q)$, $q$ odd: $L(2,q;1,1)$;
\item $T(2,2q)$: $L(1,q;2,2)$;
\item $T(3,4)$ and $T(3,5)$: $L(3,4;1,1)$ and $L(3,5;1,1)$;
\item $P(-2,2,q)$, $q \ge 3$ odd: $L(2,q;1,1;+)$;
\item $P(-2,2,2q)$, $q > 1$: $L(1,q;2,2;+)$;
\item $P(-2,2,2) = T(3,3)$: $L(1,1;3,3)$; 
\item $P(-2,3,4)$: $L(3,2;1,1;+)$. 
\end{itemize}
\end{proof}

Table \ref{tab: examples of seifert links} below lists the Seifert links which will be frequently referenced in this article. Each is an $ADE$ link up to orientation.

\begin{small} 
{\renewcommand{\arraystretch}{1.3}  
\begin{table}[ht]
\centering 
\begin{tabular}{| m{2.2cm} || m{2.2cm} | m{5.7cm} | l | m{2.2cm} |}
\hline
Seifert link & Also denoted & Description & ADE? & $ADE$ up to orientation? \\ 
\hline \hline
$L(1,q;2,2)$ & $T(2,2q)$ & The $2$-component $(2,2q)$ torus link & Yes &  \\     
\hline 
$L(1,q;2,0)$ & $T(2,2q)'$ & The $2$-component $(2, 2q)$ torus link with the orientation of one of its components reversed, $q \ge 2$   & No & Yes \\ 
\hline  
$L(1,q;2,2;+)$ & $P(-2,2,2q)$ & The $3$-component $(-2,2,2q)$ pretzel link. When $q = 1$, we have $P(-2,2,2)=T(3,3)$, which is also equivalent to $L(1,1;3,3)$ & Yes & \\ 
\hline  
$L(1,q;2,0; +)$ & $P(-2,2,2q)'$ & A reoriented version of the pretzel link $P(-2,2,2q)$, which is equivalent to $L(1,q;2,0; -)$. It is also equivalent to $L(1,1;3,1)$, when $q = 1$. & No & Yes\\
\hline  
$L(1,q;2,2; -)$ & $P(-2,2,2q)''$ & A reoriented version of the pretzel link $P(-2,2,2q)$. When $q = 1$, it is equivalent to $L(1,1;3,1)$ & No & Yes\\
\hline 
$L(2,q;1,1) \;\;$ $q$ odd& $T(2,q)$& The  $(2,q)$ torus knot &  Yes &  \\ 
\hline  
$L(2,q;1,1;+)$ $q$ odd &$P(-2,2,q)$ & The $2$-component $(-2,2,q)$ pretzel link & Yes & \\ 
\hline 
$L(2,q;1,1;-)$ $q$ odd & $P(-2,2,q)'$& The $2$-component $(-2,2,q)$ pretzel link with the orientation of one of its components reversed & No &  Yes \\ 
\hline 
$L(3,2;1,1;+)$ & $P(-2,3,4)$&The $2$-component $(-2,3,4)$ pretzel link & Yes & \\
\hline  
$L(3,2;1,1;-)$ & $P(-2,3,4)'$& The $2$-component $(-2,3,4)$ pretzel link with the orientation of one of its components reversed & No & Yes \\
\hline  
$L(3,4;1,1)$ & $T(3, 4)$ & The $(3, 4)$ torus knot & Yes &  \\
\hline  
$L(3,5;1,1)$ & $T(3, 5)$& The $(3, 5)$ torus knot  & Yes &  \\
\hline
\end{tabular}
\caption{The $ADE$ links up to orientation}
\label{tab: examples of seifert links}
\end{table}
}
\end{small}

Next, we prove Proposition \ref{prop: basic props}.

\begin{proof}[Proof of Proposition \ref{prop: basic props}]
Assertion (1) follows from the fact that a Seifert fibre space with non-empty boundary is irreducible. If $L$ is a composite Seifert link than its exterior $X(L)$ contains an essential separating meridional annulus $A$. Since $A$ is separating, it is vertical with respect to the Seifert structure on $X(L)$. Otherwise $X(L)$ would be a twisted $I$-bundle over the Klein bottle, which is impossible. But if $L$ is of type II, III or IV, and is not a Hopf link, then no meridian of $L$ is a Seifert fibre of $X(L)$. This proves (2). Finally, (3) follows from \cite{EN}; see the discussion in \cite[Proof of Corollary 1.3]{Ish}.
\end{proof}

The proof of Proposition 10.2 of \cite{BGH1} now shows that:

\begin{cor}
If $L$ is a prime Seifert link then $\Sigma_\psi(L)$ is irreducible for all cyclic branched covers $\Sigma_\psi(L)$ of $L$. 
\end{cor}

Finally, $\mathcal{P}$ will denote the following family of Seifert links:
\begin{equation}
    \label{equ: P}
\mathcal{P} = \{\#_{k} H_+, \; L(p, q; k, k), \;L(p, q; k, k; +), \;L(p, q; k, k; +, +) \}
\end{equation}

Clearly these are braid positive, and every Seifert link can be reoriented to be a member of $\mathcal{P}$.

\section{Left-orderable cyclic branched covers of Seifert links}
\label{sec: def seif links}

Our goal in this section is to prove Theorem \ref{thm: unoriented $ADE$} and classify the Seifert links for which some canonical cyclic branched cover has a non-left-orderable fundamental group. One of our main tools is the following theorem which, after noting that $PSL_2(\mathbb{R})$ is a subgroup of $\mbox{Homeo}_+(S^1)$ and $\widetilde{PSL_2}(\mathbb{R})$ is a subgroup of $\widetilde{\mbox{Homeo}}_+(S^1)\leq \mbox{Homeo}_+(\mathbb{R})$,     follows from the proof of \cite[Theorem 5.1]{BGH2}.

\begin{thm}
\label{thm: LO cbc} 
Let $L = K_1\cup \cdots \cup K_m$ be a link and let $\psi: \pi_1(X(L)) \to \mathbb Z/n$ be an epimorphism such that $\psi(\mu_i) = a_i \ne 0$, $1 \le i \le m$. Suppose that there exists a representation $\rho: \pi_1(X(L)) \to PSL_2(\mathbb{R})$ with non-cyclic image such that $\rho(\mu_i)$ is conjugate to rotation by $2\pi a_i/n$, $1 \le i \le m$. Then there exists a non-trivial representation $\pi_1(\Sigma_\psi(L)) \to \widetilde{PSL_2}(\mathbb{R})$ and therefore $\pi_1(\Sigma_\psi(L))$ is left-orderable.
\end{thm}

\begin{proof}
Let  $n_i = n/(a_i,n)$ be the order of $\psi(\mu_i)$ in $\mathbb Z/n$ and set $n_* = (n_1, n_2, \ldots, n_m)$ and $\mu_* = (\mu_1,\mu_2, \ldots, \mu_m)$. Denote the fundamental group of the orbifold $X(L)(\mu_*; n_*)$ (see \cite[Section 3.1]{BGH2}) by $G$. 

Our hypotheses imply that $\psi$ and $\rho$ factor through homomorphisms $\bar \psi: G \to \mathbb Z/n$ and $\bar \rho: G \to PSL_2(\mathbb{R})$ where, by construction, the cover of $X(L)(\mu_*; n_*)$ corresponding to $\bar \psi$ is $\Sigma_\psi(L)$. Hence $\pi_1(\Sigma_\psi(L)) \cong \mbox{ker}(\bar \psi) \leq G$. 

If $\Sigma_\psi(L)$ has positive first Betti number the conclusion is clear, so assume that $\Sigma_\psi(L)$ is a rational homology $3$-sphere. It is shown in \cite[Proof of Theorem 5.1]{BGH2} that the Euler class $e(\bar\rho|_{\pi_1(\Sigma_\psi(L))})$ is zero, and hence $\bar \rho|_{\pi_1(\Sigma_\psi(L))}$ lifts to a representation $\pi_1(\Sigma_\psi(L)) \to \widetilde{PSL_2}(\mathbb{R})$. Then $\pi_1(\Sigma_\psi(L))$ is left-orderable by \cite[Theorem 1.1]{BRW}. 
\end{proof}

\subsection{Strategy and results}
Suppose that $L$ is a prime Seifert link.
We orient $L$ as a positive braid link, i.e. as a member of $\mathcal{P}$, and take $\mu_i$ to be a positively oriented meridional class of the $i^{th}$ component $K_i$ of $L$ for each $i$. 

In Section \ref{subsec: seq of orbs} we introduce a sequence of orbifolds $\bar{\mathcal{B}}_2(L)$, $\bar{\mathcal{B}}_3(L)$, $\bar{\mathcal{B}}_4(L), \ldots$ depending only on the underlying unoriented link,  with the property that 
\begin{equation} 
\label{eqn: decreasing} 
\chi(\bar{\mathcal{B}}_2(L)) > \chi(\bar{\mathcal{B}}_3(L)) > \cdots > \chi(\bar{\mathcal{B}}_n(L)) > \cdots 
\end{equation} 
Theorem \ref{thm: unoriented $ADE$} follows immediately from the next proposition. 

\begin{prop} 
\label{prop: sfs 2 implies all}
Suppose that $L$ is a prime Seifert link.
\begin{enumerate}[leftmargin=*] 
\setlength\itemsep{0.3em}
\item[{\rm (1)}] $\chi(\bar{\mathcal{B}}_2(L)) > 0$ if and only if $\pi_1(\Sigma_2(L))$ is finite and this occurs if and only if $L$ is an $ADE$ link up to orientation.
\item[{\rm (2)}] If $\chi(\bar{\mathcal{B}}_2(L)) \leq 0$ then the fundamental group of any cyclic branched cover $\Sigma_\psi(L)$ of $L$ has a non-trivial representation into $\widetilde{PSL_2}(\mathbb{R})$, so $\pi_1(\Sigma_\psi(L))$ is left-orderable.  
\end{enumerate}
\end{prop}

The proof of this proposition will occupy most of this section; Proposition \ref{prop: sfs 2 implies all}(1) will follow from Proposition \ref{prop: chi > 0} and Proposition \ref{prop: sfs 2 implies all}(2) will be proved in Section \ref{subsec: links with chi(B) <= 0}. 

If $\chi(\bar{\mathcal{B}}_n(L)) > 0$ for some $n \geq 2$, then Sequence (\ref{eqn: decreasing}) and Proposition \ref{prop: sfs 2 implies all}(1) imply that $L$ is an $ADE$ link up to orientation, and we will see in Proposition \ref{prop: chi > 0} that $\pi_1(\Sigma_j(L))$ is finite for $1 \leq j \leq n$ and therefore not left-orderable. 

Suppose, on the other hand, that $\chi(\bar{\mathcal{B}}_n(L)) \leq 0$ for some $n \geq 2$. When $n = 2$, Proposition \ref{prop: sfs 2 implies all}(2) shows that $\pi_1(\Sigma_n(L))$ is left-orderable and so it is natural to ask whether $\pi_1(\Sigma_n(L))$ is left-orderable when $\chi(\bar{\mathcal{B}}_n(L)) \leq 0$ for some $n \geq 3$. But this is not true. Dabkowski, Przytycki and Togha have shown that all cyclic branched covers of $L = L(1,q; 2,0) = T(2, 2q)'$ have non-left-orderable fundamental groups (\cite[Theorem 1(a)]{DPT}), although  $\chi(\bar{\mathcal{B}}_n(L)) \leq 0$ for $n, q \geq 3$. In Section \ref{subsec: canonical branched covers not lo} we classify the prime Seifert links $L$ and values of $n \geq 3$ such that $\chi(\bar{\mathcal{B}}_n(L))) \leq 0$ and $\pi_1(\Sigma_n(L))$ not left-orderable. Here is the result. 
\begin{prop}
\label{prop: exceptions}
Suppose that $L$ is a prime oriented Seifert link. Then $\chi(\bar{\mathcal{B}}_n(L))) \leq 0$ and $\pi_1(\Sigma_n(L))$ is not left-orderable for some $n \ge 2$ if and only if $L$ and $n$ are either  
\begin{itemize}
\setlength\itemsep{0.3em} 
\item $L(1, q; 2, 0) = T(2, 2q)'$ and $n \geq 4$ when $q = 2$ or $n \geq 3$ when $q \geq 3$, or
\item $L(1,q; 2, 0; +) = P(-2, 2, 2q)'$, $q \ge 1$ and $n \geq 3$, or
\item $L(2, 3; 1, 1; -) = P(-2,2,3)'$ and $n = 3$.
\end{itemize}
Further, $L$ is an $ADE$ link up to orientation, though not an $ADE$ link.   
\end{prop} 

The results of this section lead to a complete description of the status of the canonical $n$-fold cyclic branched covers $\Sigma_n(L)$ of Seifert links $L$ with respect to the properties $NLS$, $LO$, and $CTF$. We state the result in terms of the torus and pretzel link notation of Table \ref{tab: examples of seifert links} for $ADE$ links up to orientation. Let $\ast$ denote any of the properties $NLS$, $LO$, or $CTF$. The following theorem will be proved at the end of Section \ref{subsec: ccbcs with nlo}.

\begin{thm}
\label{thm: canonical br covers}
Let $L$ be a prime Seifert link.
\begin{enumerate}[leftmargin=*] 
\setlength\itemsep{0.3em} 
\item[{\rm (1)}] If $L$ is not an $ADE$ link up to orientation then $\Sigma_n(L)$ is $\ast$ for all $n \ge 2$.
\item[{\rm (2)}] If $L$ is $T(2,2)$ or $q \geq 1$ and $L$ is either $T(2,2q)'$ or $P(-2,2,2q)'$, then $\Sigma_n(L)$ is not $\ast$ for all $n \ge 2$.
\item[{\rm (3)}] If $L$ is $P(-2,2,3)'$, $P(-2,3,4)'$, or $P(-2,2,2q)''$ with $q \ge 2$, or is an $ADE$ link other than $T(2,2)$, then $\Sigma_n(L)$ is not $*$ if and only if $2 \le n \le N$ where 
$$N = \left\{ \begin{array}{ll} 
5 & \mbox{ for } L = T(2, 3),    \\ 
3& \mbox{ for } L = T(2, 4), T(2, 5), \mbox{ or } P(-2,2,3)',  \\ 
2& \mbox{ otherwise.} 
\end{array} \right.$$
\item[{\rm (4)}] $\Sigma_n(L)$ is $*$ if and only if there exists a non-trivial representation $\pi_1(\Sigma_n(L)) \to \widetilde{PSL_2}(\mathbb{R})$.
\end{enumerate}
\end{thm}

\subsection{A sequence of orbifolds associated to a Seifert link}
\label{subsec: seq of orbs} 

A prime Seifert link $L$ is of the form $L(p, q; k, w), L(p, q; k, w; \varepsilon_1)$, or $L(p, q; k, w; \varepsilon_1, \varepsilon_2)$, $k \ge 1$. We denote its components by  $K_i$, where $kC_{p,q} = K_1 \cup \cdots \cup K_k$, $K_{k+1}$ corresponds to $C_1$ if it occurs in $L$, and $K_{k+2}$ corresponds to $C_2$ if it occurs. 

Orient the meridional class $\mu_i \in \pi_1(X(L))$ of $K_i$ positively with respect to the orientation on $K_i$. We use $L^*$ to denote $L$ re-oriented as a member of $\mathcal{P}$ and $\mu_i^* \in \{\mu_i^{\pm 1}\}$ the meridional class of $K_i$ positively oriented with respect to $L^*$. 
Without loss of generality we can suppose that 
\vspace{-.2cm}
\begin{itemize}

\item $\mu_i^* = \mu_i$ for $1 \leq i \leq k_+$;

\vspace{.2cm} \item $\mu_i^* = \mu_i^{-1}$ for $k_+ + 1 \leq i \leq k$.

\end{itemize}

If $P_r$ denotes a compact, connected subsurface of $S^2$ with $r \geq 0$ boundary components, the base orbifold $\mathcal{B}_L$ of the exterior of $L$ is given by
\vspace{-.2cm}
\begin{itemize}

\item $P_k(p, q)$ when $L = L(p, q; k, w)$;

\vspace{.2cm} \item $P_{k+1}(p)$ when $L = L(p, q; k, w; \varepsilon)$;

\vspace{.2cm} \item $P_{k+2}$ when $L =  L(p, q; k, w; \varepsilon_1, \varepsilon_2)$.

\end{itemize}
Let $b_i \cong S^1$ denote the boundary component of $\mathcal{B}_L$ corresponding to $K_i$, oriented coherently with some fixed orientation on $\mathcal{B}_L$, and let $\gamma_i \in \pi_1(\mathcal{B}_L)$ be the class carried by $b_i$. The reader will verify that up to changing the orientation on $\mathcal{B}_L$, the image of  $\mu_i^*$ under the epimorphism $\pi_1(X(L)) \to \pi_1(\mathcal{B}_L)$ is
\vspace{-.2cm}
\begin{itemize}

\item $\gamma_i$ for $1 \leq i \leq k$;

\vspace{.2cm} \item $\gamma_{k+1}^{q}$ for $i = k+1$, if $C_1$ occurs in $L$;

\vspace{.2cm} \item $\gamma_{k+2}^{p}$ for $i = k+2$, if $C_2$ occurs in $L$.

\end{itemize}
Fix $n \geq 2$ and let $\bar{\mathcal{B}}_n(L)$ be the orbifold obtained by attaching orbidisks $\mathcal{D}_i$ to $\mathcal{B}_L$ along the $b_i$,  where
\vspace{-.2cm}
\begin{itemize}

\item $\mathcal{D}_i = D^2(n)$ for $1 \leq i \leq k$; 

\vspace{.2cm} \item $\mathcal{D}_i = D^2(nq)$ if $i = k+1$ and $C_1$ occurs in $L$;  

\vspace{.2cm} \item $\mathcal{D}_i = D^2(np)$ if $i = k+ 2$ and $C_2$ occurs  in $L$.

\end{itemize}
In other words, $\mathcal{D}_i = D^2(nm_i)$ where $m_i$ is the distance $\Delta(\mu_i, h_i)$ of $\mu_i$ to a regular fibre $h_i$ on $\partial N(K_i)$ of the Seifert fibring of $X(L)$.

Thus $\bar{\mathcal{B}}_n(L)$ is a $2$-sphere with either $k, k+1$, or $k+2$ cone points depending on the values of $p$ and $q$. Moreover,  
$$\chi(\bar{\mathcal{B}}_n(L)) = \chi(\mathcal{B}_L) + \sum_{i} \frac{1}{nm_i},$$
where the index $i$ in the sum corresponds to the $i^{th}$ component of $L$. 
>From above we see that
\begin{equation}
\label{eqn: stable} 
\chi(\bar{\mathcal{B}}_2(L)) > \chi(\bar{\mathcal{B}}_3(L)) > \cdots > \chi(\bar{\mathcal{B}}_n(L)) > \cdots
\end{equation}

Let $X_n(L) \to X(L)$ be the $n$-fold cyclic cover associated to the homomorphism 
\begin{equation}
    \label{equ: epi}
    H_1(X(L)) \to \mathbb Z/n, \alpha \mapsto \mbox{lk}(\alpha, L)
\end{equation}
If $\tilde \mu_i$ is the inverse image of $\mu_i$ in $X_n(L)$, then $\Sigma_n(L)$ is the Dehn filling of $X_n(L)$ along the $\tilde \mu_i$. Let $\mathcal{O}_n(L)$ be the base orbifold of the Seifert manifold $\Sigma_n(L)$.

\begin{lemma}
\label{lemma: cover}
There is a cover $\mathcal{O}_n(L) \to \bar{\mathcal{B}}_n(L)$ of degree $n/r$, where $r$ is the order in $\mathbb Z/n$ of the image of the class of a regular fibre of the Seifert structure on $X(L)$ under the map in (\ref{equ: epi}). 
\end{lemma}

\begin{proof}
The cover $X_n(L) \to X(L)$ induces a Seifert structure on $X_n(L)$ and a cover $\widetilde{\mathcal{B}}_n(L) \to \mathcal{B}_L$, where $\widetilde{\mathcal{B}}_n(L)$ is the base orbifold of $X_n(L)$. Since the inverse image in $X_n(L)$ of a regular fibre of $X(L)$ has $n/r$ components, the cover $\widetilde{\mathcal{B}}_n(L) \to \mathcal{B}_L$ is of degree $n/r$. Hence if $C_i$ is the boundary component of $\mathcal{B}_n(L)$ corresponding to $K_i$ and $\widetilde C_i$ its inverse image in $\widetilde{\mathcal{B}}_n(L)$, the cover of base orbifolds restricts to a degree $n/r$ cover $\widetilde C_i \to C_i$. 

The base orbifold of $\Sigma_n(L)$ is obtained by attaching orbi-disks $D^2(\tilde m_i)$ to $\widetilde{\mathcal{B}}_n(L)$ along the $\widetilde C_i$, where $\tilde m_i$ is the distance $\Delta(\tilde \mu_i, \tilde h_i)$ of $\tilde \mu_i$ to a regular fibre $\tilde h_i$ of $X_n(L)$ on its boundary component containing $\tilde \mu_i$. Recall $m_i = \Delta(\mu_i, h_i)$, so as the inverse image of $h_i$ in $X_n(L)$ has $n/r$ components, the reader will verify that 
$$\tilde m_i = rm_i$$
Thus the degree $n/r$ cover $\widetilde C_i \to C_i$ extends to a cyclic cover of $D^2(\tilde m_i) = D^2(r m_i) \to D^2(nm_i)$ and therefore the degree $n/r$ cover $\widetilde{\mathcal{B}}_n(L) \to \mathcal{B}_L$ extends to a cover $\mathcal{O}_n(L) \to \bar{\mathcal{B}}_n(L)$.  
\end{proof}

\begin{remarks}
\label{rem: r}
$\;$ 

(1) The image of the class of a regular fibre of $X(L)$ under the homomorphism $\psi: H_1(X(L)) \to \mathbb Z/n, \alpha \mapsto \mbox{lk}(\alpha, L)$ is congruent (mod $n$) to 
$$s = \left\{
\begin{array}{ll}
wpq & \mbox{ if } L = L(p, q; k, w) \\
wpq + \varepsilon_1 p & \mbox{ if } L = L(p, q; k, w;\varepsilon_1) \\
wpq + \varepsilon_1 p + \varepsilon_2 q& \mbox{ if } L = L(p, q; k, w;\varepsilon_1, \varepsilon_2)  
\end{array} \right.$$
Thus
$$r = n/\gcd(n, s)$$

(2) $\bar{\mathcal{B}}_n(L)$ is the base orbifold of the Seifert $3$-orbifold $L(\mu_*; n)$ with underlying space $S^3$ and singular set $L$ with $\mathbb Z/n$ isotropy. The cover $\mathcal{O}_n(L) \to \bar{\mathcal{B}}_n(L)$ is induced by the Seifert-fibre-preserving cover $\Sigma_n(L) \to L(\mu_*; n)$. 

\end{remarks}

\subsection{Seifert links with \texorpdfstring{$\chi(\bar{\mathcal{B}}_n(L)) > 0$}{chi(bar Bn(L))> 0}} 
\label{subsec: spherical case} 

The $2$-orbifolds with underlying space $S^2$ and positive Euler characteristic are $S^2(a,b)$, where $a, b \geq 1$, $S^2(2,2, c)$ where $c \geq 2$, $S^2(2,3,3)$, $S^2(2,3,4)$, and $S^2(2,3, 5)$. As such it is easy to determine the prime Seifert links with $\chi(\mathcal{B}_n(L)) > 0$. After removing redundancies using Remark \ref{rem: notation conventions}, the results are compiled in the following table where $L^*$ denotes $L$ re-oriented as a member of $\mathcal{P}$, which in each case is an $ADE$ link (see Table \ref{tab: examples of seifert links}), and the final column lists the associated $ADE$ graph $\Gamma$. :

{\small
\begin{center}
\begin{tabular}{|c||c|c|c|c|} \hline 
& $L$ & $\bar{\mathcal{B}}_n$ &$L^*$ & $\Gamma$  \\  \hline \hline 

$n \geq 2$   & $L(1,1; 2,w)$  &  $S^2(n, n)$ &  $L(1,1; 2,2)$ & $A_1$     \\  \hline 

$n = 2$  &$L(1, q; 2, w)$& $S^2(2,2,q)$ &$L(1,q; 2,2)$ & $A_{2q-1}$ \\ 
& $L(2,q; 1,1)$  &  $S^2(2,2,q)$   &$L(2,q; 1,1)$& $A_{q-1}$      \\  
&$L(1,q; 2, w; \varepsilon_1)$ & $S^2(2,2,2q)$ &$L(1,q; 2, 2; +)$ &$D_{2q+2}$ \\ 
&$L(2,q; 1,1; \varepsilon_1)$&$S^2(2,2,2q)$&$L(2,q; 1,1; +)$& $D_{q+2}$ \\
& $L(3,4; 1,1)$    &  $S^2(2,3,4)$   & $L(3,4; 1,1)$&  $E_6$  \\ 
&$L(3,2; 1,1; \varepsilon_1)$& $S^2(2,3, 4)$ &$L(3,2; 1,1; +)$& $E_7$ \\
&$L(3,5; 1,1)$ &$S^2(2,3,5)$&$L(3,5; 1,1)$& $E_8$ \\ \hline 
$n = 3$ &      $L(2,3; 1,1)$      & $S^2(2,3,3)$ & $L(2,3; 1,1)$& $A_2$     \\   
&$L(1, 2; 2, w)$ & $S^2(2,3,3)$&$L(1, 2; 2, 2)$& $A_{3}$ \\ 
&$L(2,5; 1,1)$&$S^2(2,3,5)$& $L(2,5; 1,1)$& $A_4$ \\  \hline 
$n = 4$ &$L(2,3; 1,1)$ & $S^2(2,3, 4)$     & $L(2,3; 1,1)$& $A_2$     \\  \hline 
$n = 5$ &$L(2,3; 1,1)$ &   $S^2(2,3,5)$   &$L(2,3; 1,1)$ & $A_2$  \\  \hline 
\end{tabular}
\end{center}
}
\bigskip

Proposition \ref{prop: sfs 2 implies all}(1) follows from the next result. 

\begin{prop}
\label{prop: chi > 0}
Let $L$ be a prime Seifert link. 
\begin{enumerate}[leftmargin=*] 
\setlength\itemsep{0.3em}
\item[{\rm (1)}]  $\chi(\bar{\mathcal{B}}_n(L))  > 0$ for some $n \ge 2$ if and only if $L$ is an $ADE$ link up to orientation.
\item[{\rm (2)}] $\pi_1(\Sigma_n(L))$ is finite if and only if $\chi(\bar{\mathcal{B}}_n(L))  > 0$. 
\item[{\rm (3)}]If $\pi_1(\Sigma_n)$ is finite then $\pi_1(\Sigma_r(L))$ is finite for $2 \leq r \leq n$.
\end{enumerate}
\end{prop}

\begin{proof}
The first assertion follows from the table immediately above. For the second, note that the base orbifold of a closed Seifert fibred manifold has a positive Euler characteristic if and only if it is either reducible or has a finite fundamental group \cite[Theorem 5.3(ii)]{Scott}, so as the cyclic branched covers $\Sigma_n(L)$ are irreducible manifolds, $\pi_1(\Sigma_n(L))$ is finite if and only if $\chi(\mathcal{O}_n(L)) > 0$. Lemma \ref{lemma: cover} shows that the latter holds if and only if $\chi(\bar{\mathcal{B}}_n(L))  > 0$, which proves (2). Assertion (3) then follows from (\ref{eqn: stable}). 
\end{proof}

\subsection{Seifert links with \texorpdfstring{$\chi(\bar{\mathcal{B}}_n(L)) = 0$}{ chi(\bar{\mathcal{B}}n(L)) = 0}} 
\label{sec: chi = 0}
Since $S^2(a_1, a_2, \ldots, a_r)$ is Euclidean if and only if $\chi(\mathcal{\bar{B}}_n(L))= 0$ and hyperbolic if and only if $\chi(\mathcal{\bar{B}}_n(L)) < 0$, (\ref{eqn: stable}) combines with Proposition \ref{prop: chi > 0} to show that:

\begin{prop}
\label{prop: chi leq 0}
Let $L$ be a prime Seifert link. 
\begin{enumerate}[leftmargin=*] 
\setlength\itemsep{0.3em}
\item[{\rm (1)}] If $\bar{\mathcal{B}}_n(L)$ is Euclidean then $\pi_1(\Sigma_r(L))$ is finite for $r < n$ and $\bar{\mathcal{B}}_r(L)$ is hyperbolic for $r > n$. 
\item[{\rm (2)}]If $\chi(\bar{\mathcal{B}}_n(L)) <  0$ then $\bar{\mathcal{B}}_r(L)$ is hyperbolic for $r \geq n$.  
\end{enumerate}
\end{prop}

The Euclidean $2$-orbifolds with underlying space $S^2$ are $S^2(2,3,6), S^2(2,4,4), S^2(3,3,3)$, and $S^2(2,2,2,2)$. As such it is easy to determine when $\mathcal{B}_n(L)$ is Euclidean. The results are compiled in following table, where we have removed redundancies using Remark \ref{rem: notation conventions} (see Table \ref{tab: examples of seifert links}). 

{\small
\begin{center}
\begin{tabular}{|c||c|c|c|c|} \hline 
 & $L$ & $\bar{\mathcal{B}}_n$ & $L^*$ & $\Delta_{L^*}(e^{2\pi i / n})$ \\  \hline \hline 
$n = 2$ & $L(1,2; 3,w)$  &  $S^2(2,2,2,2)$   & $L(1,2; 3,3)$ & $0$  \\  
             & $L(1,1; 4,w)$    &  $S^2(2,2,2,2)$   & $L(1,1; 4,4)$ & $0$ \\ \hline 
$n = 3$ &      $L(1,3; 2,w)$                           & $S^2(3,3,3)$ &   $L(1,3; 2,2)$ & $0$  \\   
              &$L(1,1; 3,w)$&$S^2(3,3,3)$&$L(1,1; 3,3)$ & $0$  \\  \hline 
$n = 4$ &$L(1,2; 2,w)$ & $S^2(2,4,4)$     & $L(1,2; 2,2)$& $0$  \\  \hline 
$n = 6$ & $L(2,3; 1,1)$ &   $S^2(2,3,6)$   & $L(2,3; 1,1)$ & $0$  \\  \hline 
\end{tabular}
\end{center}
}
\bigskip
Here the final two columns list $L$ re-oriented as a positive braid link $L^*$ and the value of the Alexander polynomial of $L^*$ evaluated at $e^{2\pi i / n}$. (See Proposition \ref{prop: II Delta L} for a calculation of $\Delta_{L^*}(t)$.) 

\begin{prop}
\label{prop: euclidean}
Given a prime Seifert link $L$, there is at most one $n$ for which $\bar{\mathcal{B}}_n(L)$ is Euclidean. Further, if $\bar{\mathcal{B}}_n(L)$ is Euclidean and either $n = 2$ or $L = L^*$ then $\Sigma_n(L)$ has a positive first Betti number.

\end{prop}
 
\begin{proof}
The first assertion is Proposition \ref{prop: chi leq 0}(1). The last column of the table shows that $\Sigma_n(L)$ has a positive first Betti number when $L = L^*$. Since $\Sigma_2(L) \cong \Sigma_2(L^*)$, this shows that $\Sigma_2(L))$ has a positive first Betti number when $\mathcal{B}_2(L)$ is Euclidean.
\end{proof}

\subsection{The fundamental group of \texorpdfstring{$\bar{\mathcal{B}}_n(L)$}{\bar{\mathcal{B}}n(L)}}
\label{subsec: pi_1 B_n}

Recall that $\bar{\mathcal{B}}_n(L)$ is a $2$-sphere with cone points whose orders are given in Section \ref{subsec: seq of orbs}.

With a slight abuse of notation we think of $\gamma_i$ as an element of $\pi_1(\bar{\mathcal{B}}_n(L))$. Then,  
$$\pi_1(\bar{\mathcal{B}}_n(L)) = \left\{
\begin{array}{ll} 
\langle \gamma_1, \gamma_2, \ldots, \gamma_{k+2} \; | \; \gamma_i^n \; (1 \leq i \leq k), \gamma_{k+1}^q, \gamma_{k+2}^p, \gamma_1 \gamma_2 \cdots \gamma_{k+2} \rangle & \mbox{ if } L =  L(p, q; k, w) \\ 

\langle \gamma_1, \gamma_2, \ldots, \gamma_{k+2} \; | \; \gamma_i^n \; (1 \leq i \leq k), \gamma_{k+1}^{nq}, \gamma_{k+2}^p,  \gamma_1 \gamma_2 \cdots\gamma_{k+2}  \rangle & \mbox{ if } L =  L(p, q; k, w; \varepsilon_1) \\ 

\langle \gamma_1, \gamma_2, \ldots, \gamma_{k+2} \; | \; \gamma_i^n \; (1 \leq i \leq k), \gamma_{k+1}^{nq}, \gamma_{k+2}^{np}, \gamma_1 \gamma_2 \cdots\gamma_{k+2} \rangle & \mbox{ if } L =  L(p, q; k, w; \varepsilon_1, \varepsilon_2)

\end{array} \right.$$
In other words, $\pi_1(\bar{\mathcal{B}}_n(L))$ is one of the groups 
$$\Delta(a_1, a_2, \ldots, a_r) = \langle \gamma_1, \gamma_2, \ldots, \gamma_{r} \; | \; \gamma_1^{a_1}, \gamma_2^{a_2}, \ldots, \gamma_r^{a_r}, \gamma_1 \gamma_2 \cdots \gamma_{r} \rangle$$
where the $a_i$ are positive integers. 
We say that $\Delta(a_1, a_2, \ldots, a_r)$ is spherical, Euclidean or hyperbolic if $S^2(a_1, a_2, \ldots, a_r)$ has that property.  
When $r \geq 3$ and each $a_i$ is $2$ or more, $\Delta(a_1, a_2, \ldots, a_r)$ is 
\vspace{-.2cm}
\begin{itemize}

\item spherical if and only if $r = 3$ and $(a_1, a_2, a_3)$ is a Platonic triple;

\vspace{.2cm} \item Euclidean if and only if $r = 3$ and $(a_1, a_2, a_3)$ is a Euclidean triple or $r = 4$ and each $a_i$ is $2$;

\vspace{.2cm} \item hyperbolic, and therefore a discrete subgroup of $PSL(2, \mathbb R)$,  otherwise.

\end{itemize}

In our case  $r = k+2$ and 
$$\pi_1(\bar{\mathcal{B}}_n(L)) = \left\{
\begin{array}{ll} 
\Delta(n, n, \ldots, n, q, p) & \mbox{ if } L =  L(p, q; k, w) \\ \\ 

\Delta(n, n, \ldots, n, nq, p) & \mbox{ if } L =  L(p, q; k, w; \varepsilon_1) \\ \\ 

\Delta(n, n, \ldots, n, nq, np) & \mbox{ if } L =  L(p, q; k, w; \varepsilon_1, \varepsilon_2)

\end{array} \right.$$
\medskip

\begin{prop} 
\label{prop: pos bd case}
If $L = L^*$ then the following are equivalent.
\item[{\rm (1)}] $\pi_1(\Sigma_n(L))$ has a non-trivial representation into $\widetilde{PSL_2}(\mathbb R)$;
\item[{\rm (2)}] $\pi_1(\Sigma_n(L))$ is left-orderable;
\item[{\rm (3)}] $\chi(\bar{\mathcal{B}}_n(L)) \leq 0$.
\end{prop}

\begin{proof}
As observed earlier, (1) implies (2).
  
Next, note that $\chi(\bar{\mathcal{B}}_n(L)) > 0$ if and only if $\pi_1(\Sigma_n(L))$ is finite by Proposition \ref{prop: chi > 0}(2), so (2) implies (3).

Suppose (3) holds. If $\chi(\bar{\mathcal{B}}_n(L)) = 0$, then the conclusion holds by Proposition \ref{prop: euclidean} so assume that $\chi(\bar{\mathcal{B}}_n(L)) < 0$. We saw in Section \ref{subsec: seq of orbs} that the image under the epimorphism $\pi_1(X(L)) \to \pi_1(\mathcal{B}_L)$ of the positively oriented meridional class $\mu_i^*$ is

\begin{itemize}
\setlength\itemsep{0.3em}
\item $\gamma_i$ for $1 \leq i \leq k$;
\item $\gamma_{k+1}^{q}$ for $i = k+1$, if $C_1$ occurs in $L$;
\item $\gamma_{k+2}^{p}$ for $i = k+2$, if $C_2$ occurs in $L$.
\end{itemize}
Since we have assumed that $L = L^*$, we have $\mu_i^* = \mu_i$. 
On the other hand, the hyperbolicity of $\bar{\mathcal{B}}_n(L)$ implies that there is a faithful representation $\varphi: \pi_1(\bar{\mathcal{B}}_n(L)) \to PSL_2(\mathbb R)$ such that $\varphi(\gamma_i)$ is conjugate to a rotation of angle

\begin{itemize}

\item $2\pi/n$ for $1 \leq i \leq k$;

\vspace{.2cm} \item $\left\{ \begin{array}{ll} 
2\pi/nq & \mbox{ for $i = k+1$, if $C_1$ occurs in $L$} \\
2\pi/q & \mbox{ for $i = k+1$, if $C_1$ does not occur in $L$}
\end{array} \right.$

\vspace{.2cm} \item $\left\{ \begin{array}{ll} 
2\pi/np & \mbox{ for $i = k+2$, if $C_2$ occurs in $L$} \\
2\pi/p & \mbox{ for $i = k+2$, if $C_2$ does not occur in $L$}
\end{array} \right.$

\end{itemize}
Hence the composition $\pi_1(X(L)) \to \pi_1(\bar{\mathcal{B}}_n(L)) \xrightarrow{\; \varphi \;}  PSL_2(\mathbb R)$ satisfies the hypotheses of 
Theorem \ref{thm: LO cbc} to show that $\pi_1(\Sigma_n(L))$ has a non-trivial representation into $\widetilde{PSL_2}(\mathbb R)$. This shows that (3) implies (1).
\end{proof}

\begin{cor}
\label{cor: B2 hyp implies $LO$}
If $L$ is a Seifert link then $\chi(\bar{\mathcal{B}}_2(L)) \leq 0$ if and only if $\pi_1(\Sigma_2(L))$ has a non-trivial representation into $\widetilde{PSL_2}(\mathbb R)$ and this holds if and only if  $\pi_1(\Sigma_2(L))$ left-orderable.
\end{cor}

\begin{proof}
Since $\Sigma_2(L)$ and $\bar{\mathcal{B}}_2(L)$ are independent of the orientation on $L$, this is an immediate consequence of Proposition \ref{prop: pos bd case}.   
\end{proof}

\subsection{Constructing representations of \texorpdfstring{$\pi_1(X(L))$}{pi1(X(L))}} 
\label{subsec: jn}

In this section we orient $L$ to be a member of $\mathcal P$.

Let $\Sigma_\psi(L)$ be a cyclic branched cover of $L$, where $\psi$ sends $\mu_i$ to $a_i$ (mod $n$), $0 \le a_i < n$. Our goal is to determine sufficient conditions for $\pi_1(\Sigma_\psi(L))$ to have a non-trivial representation into $ \widetilde{PSL_2}(\mathbb{R})$ based on Theorem \ref{thm: LO cbc}. We follow the same plan as that used in the proof of Proposition \ref{prop: pos bd case}, the key being to construct an appropriate representation $\pi_1(\bar{\mathcal{B}}_n(L)) \xrightarrow{\; \varphi \;}  PSL_2(\mathbb R)$.

To describe how we proceed first note that $\bar{\mathcal{B}}_n(L)$ has $r$ cone points $x_1, x_2, \ldots x_r$, where 
$$r =  \left\{ 
\begin{array}{ll}
k & \mbox{ if } L = L(1,1; k, w)  \\
k+1 & \mbox{ if } L = L(1,q; k, w) \mbox{ with } q > 1 \\ 
k+1 & \mbox{ if } L =  L(p,1; k, w) \mbox{ with } p > 1 \\ 
k+1 & \mbox{ if } L(1,q; k, w; \varepsilon_1) \\ 
k+2 & \mbox{ otherwise,}
\end{array} \right. $$
whose orders are listed in Section \ref{subsec: seq of orbs}. 

Recall from Section \ref{subsec: pi_1 B_n} that $\pi_1(\bar{\mathcal{B}}_n(L))$ is generated by classes $\gamma_1, \gamma_2, \ldots, \gamma_r$ corresponding to the boundary components of $\mathcal{B}$ and satisfying the relation $\gamma_1 \gamma_2 \cdots \gamma_r = 1$. Each $\gamma_i$ has finite order dividing the order of the cone point $x_i$. 

A homomorphism $\pi_1(\bar{\mathcal{B}}_n(L)) \xrightarrow{\; \varphi \;}  PSL_2(\mathbb R)$ can be described by elements $f_1, f_2, \ldots, f_{r} \in \pi_1(\Sigma_\psi(L)) \leq \mbox{Homeo}_+(\mathbb R)$, where $f_i$ is a lift to $\pi_1(\Sigma_\psi(L))$ of $\varphi(\gamma_i) \in PSL_2(\mathbb R)$, and the product $f_1 \circ f_2 \circ \cdots \circ f_{r}$ is translation by an integer. Hence, to construct a homomorphism $\rho$ obtained as a composition $\pi_1(X(L)) \to \pi_1(\bar{\mathcal{B}}_n(L)) \to PSL_2(\mathbb R)$ such that $\rho(\mu_i)$ is conjugate to rotation by $2\pi a_i/n$ we have to choose the $f_i$ to be conjugates of translations by appropriate $\theta_i \in \mathbb R$. 

We noted in Section \ref{subsec: seq of orbs} that the image of  $\mu_i$ under the epimorphism $\pi_1(X(L)) \to \pi_1(\mathcal{B})$ is

\begin{itemize}
\setlength\itemsep{0.3em}
\item $\gamma_i$ for $1 \leq i \leq k$,

\item $\gamma_{k+1}^{q}$ for $i = k+1$, if $C_1$ occurs in $L$,

 \item $\gamma_{k+2}^{p}$ for $i = k+2$, if $C_2$ occurs in $L$,

\end{itemize}
By taking the $\theta_i \in [0, 1)$ to be  

\begin{itemize}
    \setlength\itemsep{0.3em}
\item $\theta_i = a_i/n$ for $1 \leq i \leq k$ 

\item $\theta_{k+1} = \left\{ 
\begin{array}{ll}
0 & \mbox{ if $q= 1$ and } L = L(p,q; k, w)  \\
1/q & \mbox{ if $q > 1$ and } L = L(p,q; k, w)  \\
a_{k+1}/nq &  \mbox{ if $C_1$ occurs in $L$}
\end{array} \right.$

\vspace{.2cm} \item $\theta_{k+2} = \left\{ 
\begin{array}{ll} 
0 & \mbox{ if $p = 1$ and } L = L(p,q; k, w) \mbox{ or } L(p,q; k, w; \varepsilon_1) \\
1/p & \mbox{ if $p > 1$ and } L = L(p,q; k, w) \mbox{ or } L(p,q; k, w; \varepsilon_1) \\
a_{k+2}/np &  \mbox{ if $C_2$ occurs in $L$}
\end{array} \right.$
\end{itemize}
we obtain elements  $f_1, f_2, \ldots, f_{r} \in \widetilde{PSL}_2(\mathbb R)$ which determine a homomorphism $\pi_1(\bar{\mathcal{B}}_n(L)) \xrightarrow{\; \varphi \;}  PSL_2(\mathbb R)$ if their composition $f_1 \circ f_2 \circ \cdots \circ f_{r}$ is translation by an integer. Moreover, by the choice of the $\theta_i$ the image of $\mu_i$ under the composition $\rho$ of $\varphi$ with $\pi_1(X(L)) \to \pi_1(\bar{\mathcal{B}}_n(L))$ is conjugate to a rotation by $2\pi a_i/n$. Work of Jankins and Neumann determines when this can be done. 

By construction, $\theta_i \in [0, 1)$, so if we define
$$\sigma = \sum_i \theta_i,$$
then $0 \leq \sigma < r$. The following proposition is a consequence of \cite[Corollary 2.3]{JN}. 

\begin{prop} {\rm (Jankins-Neumann)}
\label{prop: jn}
Suppose that the $\theta_i$ are as above and $r \geq 3$. Then there are elements $f_1, f_2, \ldots,$ $f_{r} \in \widetilde{PSL}_2(\mathbb R) \leq \mbox{Homeo}_+(\mathbb R)$ where $f_i$ is a conjugate in $\widetilde{PSL}_2(\mathbb R)$ to translation by $\theta_i$ and the product $f_1 \circ f_2 \circ \cdots \circ f_{r}$ is translation by an integer. Moreover, we can choose the $f_i$ so that the image of the associated homomorphism $\pi_1(\bar{\mathcal{B}}_n(L)) \xrightarrow{\; \varphi \;}  PSL_2(\mathbb R)$ is not isomorphic to $\mathbb Z/n$ if  
$$\left\{
\begin{array}{l}
r = 3 \mbox{ and either } \sigma < 1 \mbox{ or } \sigma > 2, \mbox{or }\\ 
r = 4 \mbox{ and } \sigma \ne 2, \mbox{or } \\
r \geq 5\\
\end{array} \right.$$
\end{prop}

\begin{proof}
The existence of $f_1, f_2, \ldots, f_{r}$ is guaranteed by \cite[Corollary 2.3]{JN}. Moreover we can obtain such $f_i$ whose product is translation by any of the integers in 
$$Z = \left\{
\begin{array}{ll}
\{1, 2, \ldots , r-2\} & \mbox{if } \sigma \leq 1 \\ 
\{2, 3, \ldots , r-2\} & \mbox{if } 1 < \sigma < r-1 \\
\{2, 3, \ldots , r-1\} & \mbox{if } \sigma \geq r-1 \\
\end{array} \right.$$
Hence we can produce $r-2$ non-conjugate representations unless $1 <  \sigma<  r-1$, in which case we can produce $r-3$.  

If the image of the associated homomorphism $\pi_1(\bar{\mathcal{B}}_n(L)) \xrightarrow{\; \varphi \;}  PSL_2(\mathbb R)$ is cyclic it would conjugate into $PSO(2)$ and therefore the $f_i$ would simultaneously conjugate into the group of translations of $\mathbb R$. In this case $f_1 \circ f_2 \circ \cdots \circ f_{r}$ would be translation by $\sum_i \theta_i = \sigma$ which, by construction, is an integer. Thus if $Z \ne \{\sigma\}$ we can choose $f_1, f_2, \ldots, f_r$ so that the image of $\varphi$ is not cyclic, and this occurs if either $r = 3$ and $\sigma < 1$ or $\sigma > 2$, or $r = 4$ and $\sigma \ne 2$, or $r \geq 5$. 
\end{proof}

Recall the $n$-fold cyclic branched cover $\Sigma_\psi(L) \to S^3$ with branch set $L$ defined by $\psi$ at the beginning of this subsection. 

\begin{cor}
\label{cor: sn is lo}
Suppose that $\bar{\mathcal{B}}_n(L)$ has $r$ cone points where one of the following conditions holds: 

\begin{itemize}
\setlength\itemsep{0.2em}
\item $r = 3$ and $\sigma < 1$ or $\sigma > 2;$
 \item $r = 4$ and $\sigma \ne 2;$
\item $r \geq 5$. 

\end{itemize}
Then $\pi_1(\Sigma_\psi(L))$ has a non-trivial representation into $\widetilde{PSL_2}(\mathbb R)$ for any cyclic branched cover $\Sigma_\psi(L)$ of $L$, so is left-orderable.  
\end{cor}
In the next subsection we show that when $\chi(\bar{\mathcal{B}}_2(L)) \leq 0$, $\pi_1(\Sigma_\psi(L))$ is left-orderable even when $r= 4$ and $\sigma = 2$ or $r = 3$ and $1 \leq \sigma \leq 2$. 

\subsection{Seifert links with \texorpdfstring{$\chi(\bar{\mathcal{B}}_2(L)) \leq 0$}{chi(\bar{\mathcal{B}}2(L))< 0}} 
\label{subsec: links with chi(B) <= 0}

The goal of this section is to prove Proposition \ref{prop: sfs 2 implies all}(2). We continue to use the notation of the previous subsection. In particular, $L$ is a prime Seifert link oriented as a positive braid link and $\psi: \pi_1(X(L)) \to \mathbb{Z}/n$ is an epimorphism which sends $\mu_i$ to $a_i$ (mod $n$), $0 < a_i < n$. We saw in Corollary \ref{cor: sn is lo} that $\pi_1(\Sigma_\psi(L))$ has a non-trivial representation into $\widetilde{PSL_2}(\mathbb R)$ unless $r = 4$ and $\sigma = 2$ or $r = 3$ and $1 \leq \sigma \leq 2$. Our goal here is to show that the same conclusion holds in these remaining cases if we assume that $\chi(\bar{\mathcal{B}}_2(L)) \leq 0$. Note that this is a necessary hypothesis: if $\chi(\bar{\mathcal{B}}_2(L)) > 0$ then $\pi_1(\Sigma_2(L))$ is finite (Proposition \ref{prop: chi > 0}(2)). The two cases are dealt with in Lemmas \ref{lemma: r = 4, sigma = 2} and \ref{lemma: r=3}, respectively. Proposition \ref{prop: sfs 2 implies all}(2) then follows immediately from these two lemmas together with Corollary \ref{cor: sn is lo}.

We assume below that 
$$\chi(\bar{\mathcal{B}}_2(L)) \leq 0$$
If $n = 2$, $\Sigma_\psi(L) = \Sigma_2(L)$, whose fundamental group has a non-trivial representation into $\widetilde{PSL_2}(\mathbb R)$ by Corollary \ref{cor: B2 hyp implies $LO$}. Assume then that 
$$n \geq 3$$
As above set 

\begin{itemize}
    \setlength\itemsep{0.5em}
\item $\theta_i = a_i/n$ for $1 \leq i \leq k$ 

\item $\theta_{k+1} = \left\{ 
\begin{array}{ll}
0 & \mbox{ if $q= 1$ and } L = L(p,q; k, w)  \\
1/q & \mbox{ if $q > 1$ and } L = L(p,q; k, w)  \\
a_{k+1}/nq &  \mbox{ if $C_1$ occurs in $L$}
\end{array} \right.$

 \item $\theta_{k+2} = \left\{ 
\begin{array}{ll} 
0 & \mbox{ if $p = 1$ and } L = L(p,q; k, w) \mbox{ or } L(p,q; k, w; \varepsilon_1) \\
1/p & \mbox{ if $p > 1$ and } L = L(p,q; k, w) \mbox{ or } L(p,q; k, w; \varepsilon_1) \\
a_{k+2}/np &  \mbox{ if $C_2$ occurs in $L$}
\end{array} \right.$
\end{itemize}
and
$$\sigma = \sum_i \theta_i \in [0, r)$$

\begin{lemma}
\label{lemma: r = 4, sigma = 2}
If $\chi(\bar{\mathcal{B}}_2(L)) \leq 0$, $r = 4$ and $\sigma = 2$, then $\pi_1(\Sigma_\psi(L))$ has a non-trivial representation into $\widetilde{PSL_2}(\mathbb R)$.
\end{lemma}

\begin{proof}
Since $r = 4$, $L$ is one of 

$$\left\{ 
\begin{array}{ll}
L(1,1; 4, 4) \mbox{ with } \bar{\mathcal{B}}_n(L) = S^2(n,n,n,n)  \\
L(1,q; 3, 3) \mbox{ with } q > 1 \mbox{ and } \bar{\mathcal{B}}_n(L) = S^2(q,n,n,n)\\ 
L(p,1; 3, 3) \mbox{ with } p > 1  \mbox{ and } \bar{\mathcal{B}}_n(L) = S^2(p,n,n,n)\\ 
L(p,q; 2, 2) \mbox{ with } p, q > 1 \mbox{ and } \bar{\mathcal{B}}_n(L) = S^2(p, q,n,n) \\ 
L(1,q; 3, 3; +)\mbox{ with } \bar{\mathcal{B}}_n(L) = S^2(nq,n,n,n) \\ 
L(p,q; 2, 2; +)\mbox{ with } p > 1  \mbox{ and }  \bar{\mathcal{B}}_n(L) = S^2(p,nq,n,n) \\ 
L(p,q; 2, 2; +, +) \mbox{ with } \bar{\mathcal{B}}_n(L) = S^2(np,nq,n,n)
\end{array} \right. $$

In the case that  $q > 1$ and $L = L(p,q; k, k; +)$ or $L(p,q; k, k; +,  +)$, replace $\theta_3$ by $\theta_3' = \frac{a_3 + n}{nq}$ 
to produce a new representation with $\sigma \ne 2$. The result then follows from Corollary \ref{cor: sn is lo}. 

Similarly if $L = L(p,q; k, k; +,  +)$ with $p > 1$, replace $\theta_4$ by $\theta_4' = \frac{a_4 + n}{np}$ 
to produce a new representation with $\sigma \ne 2$ and apply Corollary \ref{cor: sn is lo}.

Using Remark \ref{rem: notation conventions}, we are reduced to considering when $L$ is one of 
$$ \left\{ 
\begin{array}{ll}
L(p,q; 2, 2) \mbox{ with } p, q > 1 \mbox{ a link with $2$ components and }  2 = \frac{a_1}{n} + \frac{a_2}{n} + \frac{1}{p} + \frac{1}{q} \\ 
L(1,q; 3, 3) \mbox{ with } q > 1 \mbox{ a link with $3$ components and } 2 = \frac{a_1}{n} + \frac{a_2}{n} + \frac{a_3}{n} + \frac{1}{q}  
\\
L(1,1; 4, 4)  \mbox{ a link with $4$ components and } 2 = \frac{a_1}{n} + \frac{a_2}{n} + \frac{a_3}{n} + \frac{a_4}{n} 
\end{array} \right. $$
We'll consider these three possibilities one by one. 

\setcounter{case}{0}

\begin{case} 
{\rm $L = L(p,q; 2, 2)$ with $p, q > 1$ and $2 = \frac{a_1}{n} + \frac{a_2}{n} + \frac{1}{p} + \frac{1}{q}$.}
\end{case} 

If we replace $\psi$ by $-\psi$, which does not alter $\Sigma_{\psi}(L)$, each $a_i$ is replaced by $a_i' = n - a_i$ and $2 = \sigma = \frac{a_1}{n} + \frac{a_2}{n} + \frac{1}{p} + \frac{1}{q}$ is replaced by $\sigma' = 2 - \frac{a_1}{n} - \frac{a_2}{n} + \frac{1}{p} + \frac{1}{q} = 2 - \sigma + \frac{2}{p} + \frac{2}{q} =  \frac{2}{p} + \frac{2}{q} < 2 $. Corollary \ref{cor: sn is lo} then shows that $\pi_1(\Sigma_\psi(L))$ has a non-trivial representation into $\widetilde{PSL_2}(\mathbb R)$. 

\begin{case} 
{\rm $L = L(1,q; 3, 3)$ with $q > 1$ and $2 = \frac{a_1 + a_2 +a_3}{n} + \frac{1}{q}$.}
\end{case}

Replacing $\psi$ by $-\psi$, means that $\sigma$ is replaced by 
$\sigma' = 3 - \frac{a_1 + a_2 +a_3}{n} + \frac{1}{q} = 3 - \sigma + \frac{2}{q} = 1 + \frac{2}{q}$. We are done if $\sigma' \ne 2$ by Corollary \ref{cor: sn is lo}. Otherwise 
$2 = 1 + \frac{2}{q}$ and therefore $q = 2$. It follows that $\frac{a_1 + a_2 +a_3}{n} = \frac32$. Hence $n$ is even and so $n \geq 4$. 

Consider the restriction $\widetilde P_4 \to P_4$ of the $n$-fold cyclic cover $\widetilde{B}_n \to \mathcal{B}$ associated to $\psi$. 
Since $q = 2$ the reader will verify that the class $h \in H_1(X(L))$ of an appropriately oriented regular Seifert fibre equals $2(\mu_1 + \mu_2 + \mu_3)$. Thus 
$\psi(h) \equiv 2(\psi(\mu_1) + \psi(\mu_2) + \psi(\mu_3)) \equiv 2(a_1 + a_2 + a_3) = 3n \equiv 0$ (mod $n$). It follows that $h$ lifts to the cover $X_n(L)$ of $X(L)$. 
Thus $\widetilde P_4 \to P_4$ is an $n$-fold cyclic cover and therefore $\chi(\widetilde P_4) = -2n$. 

By assumption, the order $r_i$ of $a_i$ (mod $n$) satisfies $2 \leq r_i \leq n$. Then 
$|\partial \widetilde P_4| = n(\frac{1}{r_1} + \frac{1}{r_2} + \frac{1}{r_3}) \leq 3n/2$ and therefore the genus $\tilde g$ of $\widetilde P_4$ is
$$\tilde g =  (1/2)(2 - |\partial \widetilde P_4| - \chi(\partial \widetilde P_4)) >  (1/2)(2 - 3n/2 + 2n) = 1 + n/4 > 0$$
Hence $\tilde g > 0$ and therefore the underlying surface of the base orbifold of $\Sigma_\psi(L)$ has positive genus. It follows that $\Sigma_\psi(L)$ has positive first Betti number, and therefore $\pi_1(\Sigma_\psi(L))$ has a non-trivial representation into $\widetilde{PSL_2}(\mathbb R)$.

\begin{case} 
{\rm  $L = L(1,1; 4, 4)$ and $2 = (a_1 + a_2 + a_3 + a_4)/n$.}
\end{case}

As above we have $h = \mu_1 + \mu_2 + \mu_3 + \mu_4$ and so $\psi(h) = (a_1 + a_2 + a_3 + a_4) = 2n \equiv 0$ (mod $n$). Thus $h$ lifts to $X_n(L)$ so that the associated cover $\widetilde P_4 \to P_4$ is $n$-fold cyclic. Hence $\chi(\widetilde P_4) = -2n$. 

By assumption, the order $r_i$ of $a_i$ (mod $n$) satisfies $2 \leq r_i \leq n$. Then 
$|\partial \widetilde P_4| = n(\frac{1}{r_1} + \frac{1}{r_2} + \frac{1}{r_3} + \frac{1}{r_4}) \leq 2n$. It follows that the genus $\tilde g$ of $\partial \widetilde P_4$ is given by
$$\tilde g = (1/2)(2 + 2n - |\partial \widetilde P_4|) \geq (1/2)(2 + 2n - 2n) = 1$$
Then $\Sigma_\psi(L)$ has positive first Betti number, and therefore $\pi_1(\Sigma_\psi(L))$ has a non-trivial representation into $\widetilde{PSL_2}(\mathbb R)$. 
\end{proof}

\begin{lemma}
\label{lemma: r=3}
If $r = 3, 1 \leq \sigma \leq 2$ and $\chi(\bar{\mathcal{B}}_2(L)) \leq 0$, then $\pi_1(\Sigma_\psi(L))$ has a non-trivial representation into $\widetilde{PSL_2}(\mathbb R)$.
\end{lemma}

\begin{proof}
Since $r = 3$, $L$ is one of 

$$\left\{ 
\begin{array}{ll}
L(1,1; 3,3) \mbox{ with } \bar{\mathcal{B}}_n(L) = S^2(n,n,n)  \\
L(1,q; 2, 2) \mbox{ with } q > 1 \mbox{ and } \bar{\mathcal{B}}_n(L) = S^2(q,n,n)\\ 
L(p,1; 2, 2) \mbox{ with } p > 1  \mbox{ and } \bar{\mathcal{B}}_n(L) = S^2(p,n,n)\\ 
L(p,q; 1, 1) \mbox{ with } p, q > 1 \mbox{ and } \bar{\mathcal{B}}_n(L) = S^2(p, q,n) \\ 
L(1,q; 2, 2; +)\mbox{ with } \bar{\mathcal{B}}_n(L) = S^2(nq,n,n) \\ 
L(p,q; 1, 1; +) \mbox{ with } p > 1  \mbox{ and }  \bar{\mathcal{B}}_n(L) = S^2(p,nq,n) \\ 
L(p,q; 1,1; +, +) \mbox{ with } \bar{\mathcal{B}}_n(L) = S^2(np,nq,n)
\end{array} \right. $$

Then as $\chi(\bar{\mathcal{B}}_2(L)) > 0$ when $L$ is $L(1,1; 3,3), L(1,q; 2, 2), L(p,1; 2, 2), L(1,q; 2, 2; +)$, $L(2,q; 1, 1; +)$, $L(p,1; 1, 1; +)$, $L(1,q; 1,1; +, +)$, or $L(p,1; 1,1; +, +)$, we are reduced to considering  
$$L = \left\{ 
\begin{array}{ll}
L(p,q; 1, 1) \mbox{ with } p, q > 1 \\ 
L(p,q; 1, 1; +) \mbox{ with } p > 2  \mbox{ and }  q > 1 \\ 
L(p,q; 1,1; +, +) \mbox{ with } p, q > 1
\end{array} \right. $$
We consider these three cases separately.

\setcounter{case}{0}

\begin{case} 
\label{case: torus knot}
{\rm $L = L(p, q; 1, 1)$ with $p, q > 1$.}
\end{case} 

Then $L$ is the torus knot $T(p, q)$ and the surjectivity of $\psi$ implies that $\Sigma_\psi(L) = \Sigma_n(T(p,q))$. Since $n \ge 3$, $\chi(\bar{\mathcal{B}}_n(L)) < \chi(\bar{\mathcal{B}}_2(L)) \le 0$. The result then follows from Proposition \ref{prop: pos bd case}.

We remark here that Propositions \ref{prop: pos bd case} and \ref{prop: chi > 0}(2) recover the result of \cite{GLid} that the fundamental group $\pi_1(\Sigma_n(T(p,q)))$ is left-orderable if and only if it is infinite. 

\begin{case} 
\label{case: L(p,q; 1, 1; +)}
{\rm $L = L(p,q; 1, 1; +)$ with $p > 2, q > 1$.}
\end{case} 

First note that $(p,q) \ne (3,2)$ as otherwise $\bar{\mathcal{B}}_2(L) = S^2(3,4,2)$ is spherical, contrary to our hypotheses.  

If we replace $\psi$ by $-\psi$, $\sigma$ is replaced by $\sigma' = (1 + \frac{2}{p} + \frac{1}{q}) - \sigma$. Then 
$\sigma + \sigma' = 1 + \frac{2}{p} + \frac{1}{q}$, and given that $(p, q) \ne (3, 2)$ we have $\sigma + \sigma' < 2$. Therefore as $\sigma \geq 1$ we have $\sigma' < 1$ and the result follows from Corollary \ref{cor: sn is lo}.

\begin{case} 
{\rm $L = L(p,q; 1, 1; +, +)$ with $p, q > 1$. }
\end{case}

We can assume that $2 \leq q < p$ since $L(p,q; 1, 1; +, +) = L(q,p; 1, 1; +, +)$. 

If we replace $\psi$ by $-\psi$, $\sigma$ is replaced by $\sigma' = \frac{a_1'}{n} + \frac{a_2'}{np} + \frac{a_3'}{nq} = (1 + \frac{1}{p} + \frac{1}{q}) - \sigma$. Hence
$\sigma + \sigma' = 1 + \frac{1}{p} + \frac{1}{q} < 2$, so as $\sigma \geq 1$ we have $\sigma' < 1$. Corollary \ref{cor: sn is lo} then gives the result. 
\end{proof}

\subsection{Canonical cyclic branched covers with non-left-orderable fundamental groups}
\label{subsec: canonical branched covers not lo}
\label{subsec: ccbcs with nlo}

Here we prove Proposition \ref{prop: exceptions} by determining exactly when a canonical cyclic branched cover $\Sigma_n(L)$ of a prime oriented Seifert link $L$ has a non-left-orderable fundamental group. This is certainly the case if $\chi(\bar{\mathcal{B}}_n(L)) > 0$ (Proposition \ref{prop: chi > 0}) and never the case if $\chi(\bar{\mathcal{B}}_2(L)) \leq 0$ (Proposition \ref{prop: sfs 2 implies all}). Thus we restrict our analysis to the case that $\chi(\bar{\mathcal{B}}_2(L)) > 0$ and $\chi(\bar{\mathcal{B}}_n(L)) \leq 0$. In other words, $L$ is an $ADE$ link up to orientation (Proposition \ref{prop: sfs 2 implies all}(1)) and $n \geq 3$. Under these conditions we can also suppose that $L$ is not an $ADE$ link as otherwise $\pi_1(\Sigma_n(L))$ would be left-orderable by Proposition \ref{prop: pos bd case}. With these constraints, the possibilities for $L$ and $n$ are  
\begin{itemize}
\setlength\itemsep{0.3em}
\item $L(1, q; 2, 0) = T(2,2q)'$ with $q \ge 2$ and $n \geq 3$;
\item $L(1,1;3,1) = P(-2,2,2)'$ and $n \ge 3$;
\item $L(1,q; 2, 0; +) = P(-2,2,2q)'$ with $q \ge 2$ and $n \ge 3$;
\item $L(1,q; 2, 2; -) = P(-2,2,2q)''$ with $q \ge 2$ and $n \geq 3$; 
\item $L(2, q; 1,1; -) = P(-2, 2, q)'$ with $q \geq 1$ odd and $n \geq 4$ if $q = 1$, $n \geq 3$ if $q \geq 3$; 
\item $L(3, 2; 1, 1; -) = P(-2, 3, 4)'$ and $n \geq 3$. 
\end{itemize}
Since $L(1,1; 3,1), L(1,1; 2, 0; +)$, and $L(1,1; 2, 2; -)$ are the same oriented link up to reversal, we suppose that $q \geq 2$ when dealing with  $L(1,q; 2, w; \varepsilon)$. Similarly, since $L(2, 1; 1,1; -) = L(1,2; 2, 0)$, we assume that $q \geq 3$ when dealing with  $L(2, q; 1,1; -)$.

\begin{prop}
\label{prop: chi not pos and lo final}
$\;$

\begin{enumerate}[leftmargin=*] 
\setlength\itemsep{0.3em} 
\item[{\rm (1)}] If $n \geq 3$ and $L$ is either $L(3, 2; 1, 1; -) = P(-2, 3, 4)'$ or $L(1,q; 2, 2; -) = P(-2,2,2q)''$ with $q \geq 2$, then $\pi_1(\Sigma_n(L))$ has a non-trivial representation into $\widetilde{PSL_2}(\mathbb R)$, and is therefore left-orderable. 
\item[{\rm (2)}] If $n, q \geq 3$ with $q$ odd and $L = L(2, q; 1, 1; -) = P(-2, 2, q)'$, then $\pi_1(\Sigma_n(L))$ has a non-trivial representation into $\widetilde{PSL_2}(\mathbb R)$, and is therefore left-orderable, as long as $(n, q) \ne (3,3)$.
\end{enumerate}
\end{prop}

\begin{proof}
We continue to use the notation of the previous subsection. In each of the cases under consideration, $\bar{\mathcal{B}}_n(L)$ has three cone points. That is, $r = 3$. Our goal is to show that we can choose $\theta_1, \theta_2, \theta_3$ so that their sum $\sigma$ is strictly less than $1$, in which case Corollary \ref{cor: sn is lo} implies that $\pi_1(\Sigma_n(L))$ is left-orderable. Since we are only considering canonical cyclic branched covers the value of $a_i$ is $\pm 1$ (mod $n$), where the sign depends on the orientation of $K_i$, and these numbers determine the $\theta_i$ as in Section \ref{subsec: links with chi(B) <= 0}. 

If $L = L(3, 2; 1, 1; -)$, then $\bar{\mathcal{B}}_n(L) = S^2(n, 2n, 3)$ and we take $\theta_1 = 1/n, \theta_2 = (n-1)/2n$ and $\theta_3 = 1/3$. Then $\sigma = \frac{1}{n} + \frac{n-1}{2n} + \frac{1}{3} = \frac{5n + 3}{6n}$, which is less than $1$ if $n \geq 4$. On the other hand, Proposition \ref{prop: III Delta L} shows that $\Delta_L(t) = 1-t^3$, so $b_1(\Sigma_3(L)) > 0$ and therefore $\Sigma_3(L)$ has a left-orderable fundamental group. 

If $L = L(1, q; 2,2; -)$, then  $\bar{\mathcal{B}}_n(L) = S^2(n,n,qn)$ and we take $\theta_1 = \theta_2 = 1/n, \theta_3 = (n-1)/nq$. Then $\sigma =   \frac{2}{n} + \frac{n-1}{nq} = \frac{2q + n - 1}{nq}$, which is strictly less than $1$ if either $n, q \geq 3$ or $q = 2$ and $n \geq 4$.  When $q = 2$, the Alexander polynomial of $L$ is $(1-t)(1 - t^3)$ (Proposition \ref{prop: III Delta L}), so $\Sigma_3(L)$ has a positive first Betti number and consequently a left-orderable fundamental group. Thus (1) holds. 

If $L = L(2, q; 1, 1; -)$ with $q \geq 3$, then $\bar{\mathcal{B}}_n(L) = S^2(n, qn, 2)$ and we take $\theta_1 = 1/n, \theta_2 = (n-1)/qn$ and $\theta_3 = 1/2$. Then $\sigma = \frac{1}{n} + \frac{n-1}{qn} + \frac{1}{2}$ and the reader will verify that $\sigma < 1$ if and only if $(n-2)(q-2) > 2$. Since $n, q \geq 3$ with $q$ odd though not $n = q = 3$, we are done unless $q = 3$ and $n = 4$. But in this case, $L = L(2, 3; 1,1; -)$ has Alexander polynomial of $(1-t)(1+t^2)$ by Proposition \ref{prop: III Delta L}, so $\Sigma_4(L)$ has a positive first Betti number and consequently a left-orderable fundamental group. This completes the proof of (2). 
\end{proof}

Our next proposition shows that $\pi_1(\Sigma_n(L))$ is non-left-orderable in the remaining cases.

\begin{prop}
\label{prop: chi not pos and nlo}
$\pi_1(\Sigma_n(L))$ is not left-orderable if
\begin{itemize}
\setlength\itemsep{0.3em}
\item $L(1, q; 2, 0) = T(2,2q)'$ where $q \geq 1$ and $n \geq 3$; 
\item $L(1,1; 3,1) = P(-2, 2, 2)'$ and $n \geq 3$; 
\item $L(1,q; 2, 0; +) = P(-2,2,2q)'$ where $q \geq 2$ and $n \geq 3$; 
\item $L(2, 3; 1,1; -) = P(-2, 2, 3)'$ and $n = 3$.
\end{itemize}
\end{prop}

\begin{proof}
When $L = L(1, q; 2, 0)$, it is shown in \cite[Theorem 1(a)]{DPT} that $\Sigma_n(L)$ has a non-left-orderable fundamental group for each $n \geq 3$. To deal with the remaining cases, recall that $r$ is the order of the image $s$ of the class of a regular fibre of $X(L)$ under the homomorphism 
$$\psi: H_1(X(L)) \to \mathbb Z/n, \alpha \mapsto \mbox{lk}(\alpha, L)$$ 
and that there is a degree $n/r$ cover between the base orbifold $\mathcal{O}_n(L)$ of $\Sigma_n(L)$ and $\bar{\mathcal{B}}_n(L)$. Hence if $r = n$ then $\mathcal{O}_n \cong \bar{\mathcal{B}}_n(L)$. The following table lists $\bar{\mathcal{B}}_n(L)$, $s, r$ and $\mathcal{O}_n(L)$ for the links considered in the proposition. 

\medskip
{\footnotesize 
\begin{center}
\begin{tabular}{|c||c|c|c|c|c|} \hline 
$L$ & $n$ & $\bar{\mathcal{B}}_n(L)$ & $s$ (mod $n$)& $r$ & $\mathcal{O}_n(L)$ \\  \hline \hline 

$L(1,1; 3, 1)$    & $\geq 3$  &  $S^2(n,n,n)$&  $1$ & $n$   & $S^2(n,n,n)$ \\  \hline 
$L(1,q; 2, 0; +)$ & $\geq 3$& $S^2(n,n,nq)$     & $1$ &$n$ & $S^2(n,n,nq)$ \\  \hline 
$L(2, 3; 1, 1; -)$ &$3$&$S^2(2,3,9)$&$-1 \;\;\;$ &$3$& $S^2(2,3,9)$ \\ \hline 
\end{tabular}
\end{center}
}

Let $\mu_i, \lambda_i \in H_1(\partial N(K_i))$ be the standardly oriented meridian, longitude pair for the component $K_i$ of $L$. We denote the Seifert fibre class in $H_1(\partial N(K_i))$ by $h_i$.  

\setcounter{case}{0}
\begin{case}
$L = L(1,1; 3, 1)$ and $n \geq 3$.
\end{case}

Here $\mu_i \mapsto 1$ (mod $n$) under $\psi$. Moreover, in $H_1(X(L))$ we have 
$$\left\{ \begin{array}{l}
\lambda_1 = \mu_2 -  \mu_3 \mapsto 0 \; (\mbox{mod } n) \\ 
\lambda_2 = \mu_1 - \mu_3  \mapsto 0 \; (\mbox{mod } n)  \\
\lambda_3 = -\mu_1 - \mu_2  \mapsto -2 \; (\mbox{mod } n) 
\end{array} \right. \;\;\; \mbox{ and } \;\;\; \left\{ \begin{array}{l}
h_1 = \mu_1 + \lambda_1  \mapsto 1 \; (\mbox{mod } n)\\ 
h_2 = \mu_2 + \lambda_2  \mapsto 1 \; (\mbox{mod } n)\\
h_3 = -\mu_3 - \lambda_3  \mapsto 1 \; (\mbox{mod } n)
\end{array} \right.$$
We take 
$$\left\{ \begin{array}{l}
h_1^* =  \lambda_1  \mapsto 0 \; (\mbox{mod } n)\\ 
h_2^* =  \lambda_2  \mapsto 0 \; (\mbox{mod } n)\\
h_3^* =  2 \mu_3 + \lambda_3 \mapsto 0 \; (\mbox{mod } n)
\end{array} \right.$$
and note 
$$\mu_i = \varepsilon_i(h_i  - h_i^*)$$ 
where $\varepsilon_i = 1$ if $i = 1,2$ and $-1$ if $i = 3$. The reader will verify that 
\vspace{-.2cm} 
\begin{itemize}

\item $\Delta(h_i, h_i^*) = 1$ for each $i$; 

\vspace{.2cm} \item $h_1^* + h_2^* + h_3^* = 0 \in H_1(X(L))$.  

\end{itemize}
\vspace{-.2cm}  
It can be shown that there is an oriented horizontal pair of pants $P$ properly embedded in $X(L)$ whose oriented boundary components represent $h_1^*, h_2^*$ and $h_3^*$. 
Now $P$ lifts to a planar surface $\tilde P$ in $X_n(L)$ which can be used to determine the Seifert invariants of $\Sigma_n(L)$. Indeed, if $\tilde h_1^*, \tilde h_2^*, \tilde h_3^*$ denote the oriented boundary components of $\tilde P$ then $\tilde h_i^* \mapsto h_i^*, \tilde h_i \mapsto n h_i$, so $\Delta(\tilde h_i, \tilde h_i^*) = 1$ for each $i$. Further since $\tilde \mu_i \mapsto n \mu_i$ we have
$$\tilde \mu_i = \varepsilon_i(h_i  - nh_i^*)$$ 
Hence the Seifert invariants of $\Sigma_n(L)$ (expressed in Jankins-Neumann notation) are 
$$(0; 1/n, 1/n, -1/n) \equiv (0; 1/n, 1/n, 1 - 1/n)$$ 
It is shown in \cite{Na} that $\Sigma_n(L)$ does not support a co-oriented taut foliation, so it has a non-left-orderable fundamental group (\cite{BGW}).

\begin{case}
$L = L(1,q; 2, 0; +)$ where $q \geq 2$ and $n \geq 3$.
\end{case}

Here $\mu_i \mapsto 1$ (mod $n$) under $\psi$. Moreover, 
$$\left\{ \begin{array}{l}
\lambda_1 = -q \mu_2 +  \mu_3 \mapsto 1-q \; (\mbox{mod } n) \\ 
\lambda_2 = -q \mu_1 - \mu_3  \mapsto -(1 + q) \; (\mbox{mod } n)  \\
\lambda_3 = \mu_1 - \mu_2  \mapsto 0 \; (\mbox{mod } n) 
\end{array} \right. \mbox{ and } \;\;\; \left\{ \begin{array}{l}
h_1 = q \mu_1 + \lambda_1  \mapsto 1 \; (\mbox{mod } n)\\ 
h_2 = -q \mu_2 - \lambda_2  \mapsto 1 \; (\mbox{mod } n)\\
h_3 = \mu_3 + q\lambda_3  \mapsto 1 \; (\mbox{mod } n)
\end{array} \right.$$
Set 
$$\left\{ \begin{array}{l}
h_1^* =   (1-q) \mu_1 - \lambda_1 \; \mapsto 0 \; (\mbox{mod } n)\\ 
h_2^* =  -(1 + q)\mu_2- \lambda_2 \; \mapsto 0 \; (\mbox{mod } n)\\
h_3^* =  - \lambda_3 \; \mapsto 0 \; (\mbox{mod } n)
\end{array} \right.$$ 
and note that 
$$\left\{ \begin{array}{l}
\mu_1 = h_1 + h_1^* \\ 
\mu_2 = h_2  - h_2^* \\
\mu_3 = h_3 + qh_3^*
\end{array} \right.$$

The reader will verify that 
\vspace{-.2cm} 
\begin{itemize}
\setlength\itemsep{0.3em}
\item $\Delta(h_i, h_i^*) = 1$ for each $i$; 
\item $h_1^* + h_2^* + h_3^* = 0 \in H_1(X(L))$.  
\end{itemize}
There is a pair of pants $P$ properly embedded in $X(L)$ whose oriented boundary components represent $h_1^*, h_2^*$ and $h_3^*$. 
Now $P$ lifts to a planar surface $\tilde P$ in $X_n(L)$. If $\tilde h_1^*, \tilde h_2^*, \tilde h_3^*$ denote the oriented boundary components of $\tilde P$ then $\tilde h_i^* \mapsto h_i^*$ for each $i$. Hence as $\tilde h_i \mapsto n h_i$ and $\tilde \mu_i \mapsto n \mu_i$, we have 
$$\tilde \mu_i = \left\{ \begin{array}{c}
\tilde h_1 + n\tilde h_1^* \\ 
\tilde h_2 -n  \tilde h_2^* \\
\tilde h_3 + nq\tilde h_3^* 
\end{array} \right.$$
and therefore the Seifert invariants of $\Sigma_n(L)$ are 
$$(0; -1/n, 1/n, -1/nq) \equiv (2; 1 - 1/n, 1/n, 1 - 1/nq) \equiv (1; 1/n, 1 - 1/n, 1/nq)$$ 
It follows from \cite{Na} that $\Sigma_n(L)$ does not support a co-oriented taut foliation, so it has a non-left-orderable fundamental group (\cite{BGW}).

\begin{case}
$L = L(2, 3; 1, 1; -)$ and $n = 3$
\end{case}

Here $\mu_1, \mu_2 \mapsto 1$ (mod $3$) under $\psi$. Orient $C_2$ to link $K_1$ positively and let $\mu_3, \lambda_3 \in H_1(\partial N(C_2))$ be its oriented meridional and longitudinal classes. Then $\mu_3 \mapsto 0$ (mod $3$) under $\psi$. 
Further, in $H_1(X(L \cup C_2))$ we have 
$$\left\{ \begin{array}{l}
\lambda_1 = -2 \mu_2 + 3 \mu_3  \mapsto 1 \; (\mbox{mod } 3)\\ 
\lambda_2 = -2 \mu_1 - \mu_3   \mapsto 1 \; (\mbox{mod } 3)\\
\lambda_3 = 3 \mu_1 - \mu_2   \mapsto -1 \; (\mbox{mod } 3)
\end{array} \right. \; \mbox{ and } \;\;\; \left\{ \begin{array}{l}
h_1 = 6 \mu_1 + \lambda_1   \mapsto 1 \; (\mbox{mod } 3)\\ 
h_2 = -2 \mu_2 -3 \lambda_2   \mapsto 1 \; (\mbox{mod } 3)\\
h_3 = 3 \mu_3 + 2\lambda_3   \mapsto 1 \; (\mbox{mod } 3)
\end{array} \right.$$
Set 
$$\left\{ \begin{array}{l}
h_1^* = 5 \mu_1 + \lambda_1 \mapsto 0 \; (\mbox{mod } 3)\\ 
h_2^* = - \mu_2  -2 \lambda_2 \mapsto 0 \; (\mbox{mod } 3) \\
h_3^* = -5 \mu_3 +-3 \lambda_3 \mapsto 0 \; (\mbox{mod } 3) 
\end{array} \right.$$ 
and note that 
$$\left\{ \begin{array}{l}
\mu_1 = h_1 - h_1^* \\ 
\mu_2 = -2 h_2 + 3h_2^* \\
\mu_3 = -3 h_3 - 2h_3^*
\end{array} \right.$$

The reader will verify that 
\vspace{-.2cm} 
\begin{itemize}

\item $\Delta(h_i, h_i^*) = 1$ for each $i$; 

\vspace{.2cm} \item $h_1^* + h_2^* + h_3^* = 0 \in H_1(X(L \cup C_2))$.

\end{itemize}
There is an oriented horizontal pair of pants $P$ properly embedded in $X(L \cup C_2)$ whose oriented boundary components represent $h_1^*, h_2^*$ and $h_3^*$. 
Since $\tilde h_i \mapsto 3h_i, \tilde h_i^* \mapsto h_i^*$ and  
$$\tilde \mu_i \mapsto \left\{ \begin{array}{l}
3\mu_i \mbox{ for } i = 1, 2 \\ 
\mu_3 \mbox{ for } i = 3
\end{array} \right.$$
we have 
$$\left\{ \begin{array}{l}
\tilde \mu_1 = \tilde h_1 - 3\tilde h_1^* \\ 
\tilde \mu_2 = -2 \tilde h_2 + 9 \tilde h_2^* \\
\tilde \mu_3 = -3 \tilde h_3 - 2\tilde h_3^* 
\end{array} \right.$$

Hence the Seifert invariants of $\Sigma_3(L)$ are 
$$(0; 1/3, 2/9, -3/2) \equiv (2; 1/3, 2/9, 1/2) \equiv (1; 2/3, 7/9, 1/2)$$ 
Then by \cite{Na}, $\Sigma_n(L)$ does not support a co-oriented taut foliation, so it has a non-left-orderable fundamental group (\cite{BGW}). 
\end{proof}

\begin{proof}[Proof of Theorem \ref{thm: canonical br covers}] 

Since the $L$-space conjecture holds for Seifert fibre spaces (\cite{BRW}, \cite{LS}, \cite{BGW}), it suffices to verify the theorem in the case that $* = LO$. Part (1) of Theorem \ref{thm: canonical br covers} then follows from Proposition \ref{prop: sfs 2 implies all}.

Parts (2) and (3) concern the case that $L$ is an $ADE$ link up to orientation, and in this case $\pi_1(\Sigma_2(L))$ is a finite group, hence non-left-orderable. 
Proposition \ref{prop: higher order}(a) shows that $\pi_1(\Sigma_n(T(2,2)))$ is non-left-orderable for all $n \geq 3$, and the same conclusion holds for $q \geq 1$ when $L$ is $T(2, 2q)' = L(1, q; 2, 0)$ or $P(-2, 2, 2q)' = L(1, q; 2, 0; +)$ by Proposition \ref{prop: chi not pos and nlo}. (Recall that $L(1,1; 3,1)$  coincides with $L(1, 1; 2, 0; +)$.) Thus part (2) of Theorem \ref{thm: canonical br covers} holds. 

For part (3), let $N$ be the integer defined in Theorem \ref{thm: canonical br covers}. If $L$ is an $ADE$ link other than $T(2,2)$, Propositions \ref{prop: 2-fold}  and \ref{prop: higher order} show that $\pi_1(\Sigma_n(L))$ is finite, or equivalently $\chi(\bar{\mathcal{B}}_n(L)) > 0$ (Proposition \ref{prop: chi > 0}), if and only if $2 \leq n \leq N$. Then $\chi(\bar{\mathcal{B}}_n(L)) \leq 0$ for $n > N$, in which case Proposition \ref{prop: pos bd case} implies that $\pi_1(\Sigma_n(L))$ is left-orderable. 

Proposition \ref{prop: chi not pos and lo final}(1) shows that $\Sigma_n(L)$ has a left-orderable fundamental group for $n \geq 3$ when $L$ is either $P(-2,3,4)' = L(3,2; 1,1; -)$ or $P(-2,2,2q)'' = L(1, q: 2, 2; -)$ with $q \geq 2$. Finally, when $L = P(-2,2,3)' = L(2,3; 1,1;-)$, Proposition \ref{prop: chi not pos and nlo} shows that the fundamental group $\Sigma_n(L)$ is non-left-orderable for $n = 2, 3$ while Proposition \ref{prop: chi not pos and lo final}(2) shows that it is left-orderable for $n \geq 4$. 

Finally, consider (4). If $\pi_1(\Sigma_n(L))$ admits a non-trivial representation with values in $\widetilde{PSL_2}(\mathbb R)$, a subgroup of  $\mbox{Homeo}_+(\mathbb R)$,  then it is left-orderable by \cite[Theorem 1.1]{BRW}. Conversely, the reader will note that in each case we verified the left-orderability of $\pi_1(\Sigma_n(L))$ by either constructing a non-trivial representation of it with values in $\widetilde{PSL_2}(\mathbb R)$ or showing that it has infinite abelianisation and in the latter case there is clearly a non-trivial representation $\pi_1(\Sigma_n(L)) \to \widetilde{PSL_2}(\mathbb R)$. Consequently (4) holds.
\end{proof}

\section{The classification of strongly quasipositive Seifert links} 
\label{sec: classn sqpsls}
Ishikawa classified strongly quasipositive Seifert links (\cite{Ish}). In the fibred case he did this by determining when their open books induce the tight contact structure on $S^3$. In this section, we give a different, more topological proof, which is based on comparing $3$-sphere and $4$-ball genera. In particular it does not involve any contact structure computations. We also classify the Seifert links that have equal 3-sphere and 4-ball genera (Theorem \ref{thm: g_4 = g results}), have 3-sphere genus zero (Theorem \ref{thm: g = 0}), and are definite (Theorem \ref{thm: def results}). 
 
\subsection{Strategy and results}
To describe our approach in more detail, recall that the {\it genus} of a link $L$ is $g(L) = \min \{g(F) \mid  F$ a Seifert surface for $L\}$. A Seifert matrix $S$ of $F$ is a square matrix of size $\beta_1(F)= 2g(F) + |L| - 1$ representing the Seifert form of $F$. The {\it signature} $\sigma(L)$ and {\it nullity} $\eta(L)$ of $L$ are the signature and nullity of the matrix $S + S^T$. Since $\Delta_{L}(-1) = \det(S + S^T)$, we have that $\eta(L) = 0$ if and only if $\Delta_{L}(-1) \ne 0$. The signature $\sigma(L)$ satisfies
$$ \sigma(L) \le 2g(L) + |L| - 1  $$
We say that $L$ is {\it definite} if this is an equality. 

The (smooth) {\it $4$-ball genus} of $L$ is 
\begin{align}
g_{4}(L) = \min  \{g(F) \mid F &\text{ is a compact, connected, oriented surface, } \nonumber \\
&\text{smoothly and properly embedded in } B^4, \nonumber \\
&\text{such that } \partial F = L \}\nonumber
\end{align}

Clearly $g_{4}(L) \le g(L)$. Since $\beta_{1}(F) = 2g(F) + |L| - 1$, it is sometimes more convenient to work with $h(L)$, defined as for $g(L)$ but with $g(F)$ replaced by $\beta_{1}(F)$:
$$h(L) = \min \{\beta_1(F) \; | \; F \; \text{a Seifert surface for} \; L\}$$ 
Then $|\sigma(L)| \leq h(L)$  and $L$ is definite if and only if this is an equality. Similarly $h_4(L)$ is defined as for $g_{4}(L)$, but with $g(F)$ replaced by $\beta_{1}(F)$. Clearly $h_4(L) \leq h(L)$.

\begin{lemma}
\label{lemma: conn qp}
If a link $L$ bounds a connected quasipositive surface then $g_4(L) = g(L)$.
\end{lemma}

\begin{proof}
Let $F$ be a quasipositive surface in $S^3$ with $\partial F = L$. Then by \cite{Ru1},
\begin{align}
\chi(F) = \max\{\chi(F') \mid F' & \text{ a surface with no closed components} \nonumber \\
& \text{smoothly and properly embedded in } B^4 \text{ with } \partial {F'} = L\} \nonumber 
\end{align}
If $F$ is connected it follows that $g_4(L) = g(F) = g(L)$.
\end{proof}

Note that Lemma \ref{lemma: conn qp} applies if $L$ is fibred and strongly quasipositive.
Recall the family of Seifert links defined in (\ref{equ: P}):
   $$ \mathcal{P} = \{\#_{k} H_+, \; L(p, q; k, k), \;L(p, q; k, k; +), \;L(p, q; k, k; +, +) \}$$

These are clearly braid positive. (In fact they are the only braid positive Seifert links; see Corollary \ref{cor: bp}.) As such they are fibred and strongly quasipositive, and therefore satisfy $g_{4}(L) = g(L)$ by Lemma \ref{lemma: conn qp}. Note that any link $L$ with $g(L) = 0$ also satisfies $g_{4}(L) = g(L)$.
Here are the main results of this section. 

\begin{customthm}{\ref{thm: g_4 = g results}}
Let $L$ be a Seifert link. Then $g_{4}(L) = g(L)$ if and only if either $L \in \mathcal{P}$ or $g(L) = 0$.
\end{customthm}

We explicitly determine the Seifert links of genus $0$. 

\begin{thm}
\label{thm: g = 0}
Let $L$ be a Seifert link. Then $g(L) = 0$ if and only if $L$ is either
\begin{enumerate}[leftmargin=*] 
\setlength\itemsep{0.3em}
\item[{\rm (1)}] a connected sum of Hopf links; 
\item[{\rm (2)}] $L(p, q; k, 0)$ or $L(1, 1; k, 1)$; 
\item[{\rm (3)}] $L(1, q; k, 0; \varepsilon)$, $q>1$, or $L(1, 2; k, 1; -)$;
\item[{\rm (4)}] $L(p, p+1; k, 0; -, +)$ or $L(2, 3; k, 1; -, -)$.
\end{enumerate}
\end{thm}

We show in Proposition \ref{prop: sqp g = 0} that of the links listed in parts (2), (3), and (4) of Theorem \ref{thm: g = 0}, only these in the $L(p,q;k,0)$ family are strongly quasipositive. Together with Theorem \ref{thm: g_4 = g results} this gives the following.

\begin{customthm} {\ref{thm: sqp seifert results intro}} {\rm (Ishikawa)}
\begin{enumerate}[leftmargin=*] 
\setlength\itemsep{0.3em} 
\item[{\rm (1)}] $L$ is fibred and strongly quasipositive if and only if $L \in \mathcal P$;
\item[{\rm (2)}] $L$ is strongly quasipositive if and only if either $L \in \mathcal P$ or $L$ is of the form $L(p, q; k, 0)$.
\end{enumerate}
\end{customthm}

\begin{cor}
\label{cor: bp}
  A Seifert link $L$ is braid positive if and only if $L \in \mathcal{P}$.
\end{cor}

\begin{proof}
If $L \in \mathcal P$, then $L$ is braid positive. Conversely, if $L$ is braid positive, then it is fibred and strongly quasipositive. Therefore $L \in \mathcal P$ by Theorem \ref{thm: sqp seifert results intro}.
\end{proof}

Baader has shown that a prime, braid positive link is definite if and only if it is an $ADE$  link \cite{Baa}. Also, we determine the definite links among those listed in Theorem \ref{thm: g = 0} (see Proposition \ref{prop: def g = 0}). Combining these with Theorem \ref{thm: g_4 = g results} gives the following classification of definite Seifert links.

\begin{customthm}{\ref{thm: def results}}
A Seifert link is definite if and only if it is either
\begin{enumerate}[leftmargin=*] 
\setlength\itemsep{0.3em}  
\item[{\rm (1)}] $\#_{k} H_+$, $k \ge 2$,
\item[{\rm (2)}] an $ADE$ link,
\item[{\rm (3)}] $L(p,q;2,0)$,
\item[{\rm (4)}] $L(2,3; 1,1; -,-)$.
\end{enumerate}
\end{customthm}

We prove Theorem \ref{thm: g_4 = g results} in Section \ref{subsec: g3 = g4}, Theorem \ref{thm: g = 0} is proved in Section \ref{subsec: alex ply}. Theorems \ref{thm: sqp seifert results intro} and \ref{thm: def results} are proved in Sections \ref{subsec: sqp sl} and \ref{subsec: def sl}, respectively. 

\subsection{Alexander polynomials of Seifert links and Seifert links of genus zero}
\label{subsec: alex ply}
In this section we compute the Alexander polynomials of Seifert links and prove Theorem \ref{thm: g = 0}.  We use $\Delta_L(t)$ or simply $\Delta_L$ to denote the Alexander polynomial of a link $L$. If $|L|=m > 1$, then $\Delta_L(t_1, \cdots, t_m)$ denotes the multivariable Alexander polynomial of $L$.

A Seifert link $L$ of type II $L(p,q; k, w)$, type III $L(p,q; k, w; \epsilon)$ or type IV $L(p,q; k, w; \epsilon_1, \epsilon_2)$ is obtained from the corresponding link $L'$ with $k = 1$ by replacing $C_{p,q}$ with $k$ parallel copies of $C_{p,q}$ on $T$. Alternatively, let $V$ be an unknotted solid torus in $S^3$. Let $J = J_1 \cup \cdots \cup J_k \subset \mbox{int}(V)$ be the unoriented closure of the braid $\delta^{kpq}$ in $B_k$ where $\delta = \sigma_1 \cdots \sigma_{k-1}$ is the Garside element in $B_k$, so $\delta^k$ corresponds to the positive full twist. Note that the linking number between two such (positively oriented) copies of $C_{p,q}$ is $pq$. So orienting the components of $J$ appropriately, $L$ is  a satellite of $L'$ with pattern $J$. Expressing this in terms of the {\it splicing} operation \cite{EN} will enable us to compute the Alexander polynomial of $L$.

We begin by orienting all the components of $J$ in the same direction. Let $C$ be the boundary of a meridian disk of $V$, oriented so that the linking number of $C$ with each $J_i$ is $1$. 

We then have oriented links $J$ and $L'' = J \cup C$ in $S^3$. Let $L_0''$ be the corresponding link in $S^3$ with no meridional twists in $V$. Thus $L_0''$ is a connected sum of $k$ copies of $H_+$. Then $H_1(S^3 \setminus L_0'') \cong \mathbb Z^{k+1}$, with basis $\{s, t_1, \ldots, t_k\}$, where $t_i$ is the class of a positively oriented meridian of $J_i$, $1 \leq i \leq k$, and $s$ is the class of a positively oriented meridian of $C$. One computes (e.g. using Fox's free differential calculus) that the multivariable Alexander polynomial of $L_0''$ is given by
$$\Delta_{L_0''}(s, t_1, \ldots , t_k) = (1 - s)^{k-1}$$

There is a homeomorphism from $S^3 \setminus L_0''$ to $S^3 \setminus L''$ whose induced map on first homology sends $t_i \mapsto t_i$, $1 \leq i \leq k$, and $s \mapsto s + pq(\sum_{i=1}^k t_i)$. Hence by \cite{To}, 

\begin{lemma}
\label{lemma: Delta L'' multi}
$\Delta_{L''}(s, t_1, \ldots , t_k) = (1 - sT^{pq})^{k-1}$, where $T = \prod_{i=1}^k t_i$. 
\end{lemma}

Note that $[C] = \sum_{i=1}^k t_i \in H_1(S^3 \setminus J)$. Hence by \cite{To} the Alexander polynomial of $J = L'' \setminus C$ is given by the following.

\begin{lemma}
\label{lemma: Delta J multi}
$\Delta_{J}(t_1, \ldots , t_k) = (1 - T^{pq})^{k-1}/(1 - T)$, where $T = \prod_{i=1}^k t_i$. 
\end{lemma}

Now let $J(w)$ and $L''(w)$ be obtained from $J$ and $L''$ by reversing the orientations of $J_1, \ldots , J_{k_-}$ for some $0 \leq k_- \leq k/2$, where $w$ denotes $k - 2k_- = \mbox{lk}(C, J(w))$, the homological winding number of $J(w)$ in $V$.

Since for any link $L$ with $|L| > 1$ we have the Alexander polynomial $\Delta_L(t) = (1-t) \Delta_L(t, \ldots , t)$ by \cite[Lemma 10.1]{Mil}, where $\Delta_L(t, \ldots , t)$ is the multivariable Alexander polynomial with all variables substituted with $t$, Lemma \ref{lemma: Delta J multi} implies

\begin{cor}
\label{cor: Delta J(w)}
$\Delta_{J(w)}(t) = (1-t)(1 - t^{wpq})^{k-2}(1 + t^w + \cdots + t^{w(pq - 1)})$. 
\end{cor}

Consider the $0$-core case. Here the corresponding link $L'$ with $k = 1$ is the torus knot $T(p, q)$. Then $L = L(p, q; k, w)$ is the splice of $L''(w)$ and $T(p, q)$, spliced along $C$. Using Corollary \ref{cor: Delta J(w)}, the formula for the Alexander polynomial of the splice of a link and a knot \cite[Theorem 5.2]{EN}, \cite[Corollary 3.3]{Ci} gives

\begin{prop}
\label{prop: II Delta L}
Let $L = L(p, q; k, w)$. Then 
$$\Delta_{L}(t) = (1-t)(1 - t^{wpq})^{k-2}(1 + t^w + \cdots + t^{w(pq - 1)})\Delta_{T(p, q)}(t^w)$$
\end{prop}

When $w= 0$ this is to be interpreted as $pq(1-t)$ if $k = 2$ and $0$ if $k \ge 4$.

In the $1$-core case, $L = L(p, q; k, w; \varepsilon)$ is the splice of $L''(w)$ and the $2$-component link $L' = L(p, q; 1, 1; \varepsilon)$, along $C$ and $C_{p,q}$. Let the $2$-variable Alexander polynomial of $L'$ be $\Delta_{L'}(s, t)$, where $s$ corresponds to the component $C_1$, and $t$ to the component $C_{p,q}$. This Alexander polynomial is computed by Sumners and Woods in \cite{SW}. 

\begin{lemma} {\rm (\cite[Theorem 6.2]{SW})}
\label{lemma: III Delta L''}
Let $L' = L(p, q; 1, 1; \varepsilon)$. Then $$\Delta_{L'}(s, t) = (1 - (s^\varepsilon t^q)^p)/(1 - s^\varepsilon t^q)$$
\end{lemma}

Combining this with Lemma \ref{lemma: Delta L'' multi}, the formula for the Alexander polynomial of a splice \cite[Theorem 5.2]{EN} gives

\begin{prop}
\label{prop: III Delta L}
Let $L = L(p, q; k, w; \varepsilon)$. Then 
$$\Delta_{L}(t) = (1-t)(1 - t^{(wq + \varepsilon)p})^{k-1}(1 + t^{wq + \varepsilon} + \cdots + t^{(wq  + \varepsilon)(p - 1)})$$
\end{prop}

For the  $2$-core case, $L = L(p, q; k, w; \varepsilon_1, \varepsilon_2)$ is the splice of $L''(w)$ and $L' = L(p, q; 1, 1; \varepsilon_1, \varepsilon_2)$ along $C$ and $C_{p,q}$. The $3$-variable Alexander polynomial of $L'$ is also computed by Sumners and Woods; here $s_i$ corresponds to $C_i$, $i = 1, 2$, and $t$ corresponds to $C_{p,q}$.

\begin{lemma}
{\rm (\cite[Theorem 6.3]{SW})} 
\label{lemma: IV multi Delta L'}
Let $L' = L(p, q; 1, 1; \varepsilon_1, \varepsilon_2)$. Then 
$$\Delta_{L'}(s_1, s_2, t) = 1 - s_1^{\varepsilon_1 p} s_2^{\varepsilon_2 q} t^{pq}$$
\end{lemma}

The splicing formula \cite[Theorem 5.2]{EN} then gives, using Lemma \ref{lemma: Delta L'' multi}

\begin{prop}
\label{prop: IV Delta L}
Let $L = L(p, q; k, w; \varepsilon_1, \varepsilon_2)$. Then
$$\Delta_L(t) = (1-t)(1 - t^{(wpq + \varepsilon_1 p + \varepsilon_2 q)})^k$$
\end{prop}

\begin{lemma}
\label{lemma: deg delta L}
The degree of $\Delta_L$ is given by
$$\deg \Delta_L = \left\{
\begin{array}{ll} 
1 + w(kpq - p - q) & \mbox{ if $L$ is either $0$-core with $w \ne 0$ or $w = 0$ and $k = 2$} \\
1 + (kp - 1)|wq + \varepsilon| & \mbox{ if $L$ is $1$-core} \\
1 + k|wpq + \varepsilon_{1} p + \varepsilon_{2} q| & \mbox{ if $L$ is $2$-core}
\end{array} \right.$$
\end{lemma}

\begin{proof}

  This follows from Propositions \ref{prop: II Delta L}, \ref{prop: III Delta L}, and \ref{prop: IV Delta L}, respectively.
\end{proof}

\begin{proof}[Proof of Theorem \ref{thm: g = 0}] 

If $L$ is a connected sum of Hopf links or of the form $L(p,q;k,0)$ then $g(L) = 0$ (for the latter see Lemma \ref{lemma: 0-core w = 0}).

Let $L$ be a Seifert link other than those just mentioned. Since $L$ is fibred by Proposition \ref{prop: basic props}, $\deg \Delta_L = 2g(L) + |L| - 1$, so $g(L) = 0$ if and only if $\deg \Delta_L = |L| - 1$. By Lemma \ref{lemma: deg delta L}, in the $0$-,  $1$-, or  $2$-core case this happens if and only if
$$
\begin{array}{ll}
\mbox{(i) } 1 + w(kpq -p -q) = k - 1 & \mbox{ if $L$ is $0$-core with $w \ne 0$ }\\
\mbox{(ii) }  1 + (kp - 1)|wq + \varepsilon| = k & \mbox{ if $L$ is $1$-core} \\
\mbox{(iii) } 1 + k|wpq + \varepsilon_{1} p + \varepsilon_{2} q| = k + 1 & \mbox{ if $L$ is $2$-core}
\end{array}$$

First consider the $0$-core case: $L = L(p, q; k, w)$ with $w \ne 0$. If $k = 1$ then $L = L(p, q; 1,1)$ is the torus knot $T(p,q)$, so $g(L) > 0$ since $L$ is not the unknot by hypothesis. If $k = 2$, then (i) holds if and only if $p = q = 1$.  Since $w\neq 0$, $w=2$ and $L = L(1, 1; 2, 2)$ is the positive Hopf link, which is included in case (1).

Now assume $k > 2$. If $w = 1$ then the left-hand side minus the right-hand side of (i) is 
$$(p - 1)(q - 1) + (k - 1)(pq - 1)$$ 
which is positive unless $p = q = 1$. We can exclude the case that $w > 1$ since if $d$ is the difference of the two sides of (i), then 
\begin{eqnarray}
d  \geq  1+ 2(kpq - p - q) + 1 - k  
& =  & 2(p-1)(q-1) + 2(k-1)pq - k \nonumber \\ 
& \ge  & 2(k-1) - k \nonumber \\ 
& = & k - 2 \nonumber \\  
&> & 0 \nonumber 
\end{eqnarray}

In the $1$-core case, by Remark \ref{rem: notation conventions}(4) we may assume that $q > 1$. We see that (ii) holds if and only if $p = 1$ and $|wq + \varepsilon| = 1$, and that the latter holds if and only if either $w= 0$ or $w = 1$, $q = 2$, and $\varepsilon = -$.

In the  $2$-core case, (iii) holds if and only if $|wpq + \varepsilon_{1} p + \varepsilon_{2} q| = 1$. If $w = 0$, after removing the redundancies using Remark \ref{rem: notation conventions}(5) and (6), we get $\varepsilon_1 = -$, $\varepsilon_2 = +$, and $q = p + 1$. If $w \ge 1$ then it is easy to check that $wpq + \varepsilon_{1} p + \varepsilon_{2} q > 1$ unless $\varepsilon_1 = \varepsilon_2 = -$, $w = 1$, $p = 2$, and $q = 3$.
\end{proof}

\subsection{Seifert links with \texorpdfstring{$g_3 = g_4$}{g3 = g4}}
\label{subsec: g3 = g4}
In this section, we prove Theorem \ref{thm: g_4 = g results}. Since a connected sum of Hopf links has genus 0, to prove Theorem \ref{thm: g_4 = g results} we must show that if $L$ is 0-, 1-, or 2-core and $L \notin \mathcal P$ then $g_4(L) = g(L)$ if and only if $g(L) = 0$. We distinguish three cases: (1) $w = 0$, (2) $0 < w < k$, and (3) $w = k$ and $L \notin \mathcal P$. We dispose of case (1) in Proposition \ref{prop: w = 0} below, by showing that in that case $g_4(L) = 0$. First we have

\begin{lemma}
\label{lemma: 0-core w = 0}
If $L = L(p,q;k,0)$ then $g(L) = 0$.
\end{lemma}

\begin{proof}
  By pairing oppositely oriented components it is easy to see that $L$ bounds a disjoint union of $k/2$ annuli in $S^3$. Tubing these together shows that $g(L) = 0$.
\end{proof} 

\begin{prop}
\label{prop: w = 0}
  Let $L$ be a $0$-, $1$-, or $2$-core Seifert link with $w = 0$. Then $g_{4}(L) = 0$.
\end{prop}

\begin{proof}

  If $L = L(p,q;k,0)$ then $g(L) = 0$ by Lemma \ref{lemma: 0-core w = 0}.

 Let $L = L(p,q;k,0;\varepsilon)$. Let the copies of $C_{p,q}$ in $L$ be $C_{p,q}^{(i)}$, $1 \leq i \leq k$, numbered in the order they lie on the Heegaard torus $T$. We may suppose that $C_{p,q}^{(i)}$ is positively oriented if $i$ is odd and negatively oriented if $i$ is even. Band $C_{p,q}^{(2i-1)}$ and $C_{p,q}^{(2i)}$ together by a band that lies in the subannulus of $T$ that they cobound, for $1 \leq i \leq k/2$. The resulting $k/2$ components, together with $C_1$, is the $(k/2 + 1)$-component unlink. Capping off the components of this unlink with disks produces, in $B^4$, the disjoint union of $k/2$ annuli and a disk. Tubing these together shows that $g_4(L) = 0$.

  In the  $2$-core case, the analogous banding construction produces the split sum of the $k/2$-component unlink and $C_1 \cup C_2 = H_{\pm}$. Since $H_{\pm}$ bounds an annulus in $S^3$, we get that $L$ bounds a disjoint union of annuli in $B^4$, showing that $g_4(L) = 0$ as in the previous case.
\end{proof}  

Cases (2) and (3) are proved in Section \ref{subsubsec: 0 < w < k} and Section \ref{subsubsec: w = k not P} respectively.  In case (2) we use a banding argument similar to that used in the proof of Proposition \ref{prop: w = 0} to get an upper bound on $g_4(L)$. A lower bound on $g(L)$ comes from the degree of the Alexander polynomial of $L$, calculated above in Section \ref{subsec: alex ply}. We use a similar approach in case (3), but here, since $w = k$, the banding construction is more complicated. In the 1-core case, and in the 2-core case when $\varepsilon_2 = +$, we band $C_1$ to a copy of $C_{p,q}$, and in the 2-core case when $\varepsilon_2 = -$, we band both $C_1$ and $C_2$ to a copy of $C_{p,q}$. The strategy is explained in more detail in subsection \ref{subsubsec: w = k not P}.

\subsubsection{The case \texorpdfstring{$0 < w < k$}{0 < w < k}}
\label{subsubsec: 0 < w < k}
In this section we prove Theorem \ref{thm: g_4 = g results} for $0$-,  $1$- and  $2$-core links with $0 < w < k$. We will denote such a link simply by $L(k,w)$.

\begin{lemma}
  \label{lemma: diff deg}
  
  Let $L = L(k,w)$ be a $0$-,  $1$-, or  $2$-core Seifert link with $0 < w < k$ and let $L' = L(w,w)$. Then $\deg \Delta_L - \deg \Delta_{L'} > (k-w)$ unless $g(L) = 0$.

\end{lemma}

\begin{proof}
By Lemma \ref{lemma: deg delta L}, 
 $$\deg \Delta_L - \deg \Delta_{L'} = \left\{ 
 \begin{array}{ll} 
 (k-w)wpq & \mbox{ if $L$ is $0$-core} \\ 
 (k-w)p(wq + \varepsilon)  & \mbox{ if $L$ is $1$-core} \\
 (k-w)(wpq + \varepsilon_{1}p + \varepsilon_{2}q) & \mbox{ if $L$ is $2$-core} 
 \end{array} \right.$$
By the computation in the $w \ne 0$ case in the proof of Theorem  \ref{thm: g = 0}, we see that this is larger than $(k-w)$ unless $w = p = q = 1$, or $w = 1, q = 2, \varepsilon = -$, or $w = 1, p = 2, q = 3, \varepsilon_1 = \varepsilon_2 = -$, respectively, i.e. unless $g(L) = 0$. 
\end{proof}

\begin{prop}
\label{prop: 0 < w}
  Let $L$ be a $0$-,  $1$-, or  $2$-core Seifert link with $0 < w < k$. Then $g_{4}(L) = g(L)$ if and only if $g(L) = 0$.

\end{prop}

\begin{proof}
Banding each of the negatively oriented $C_{p,q}$ components of $L = L(k,w)$ with a positively oriented component  and then capping off the resulting $(k-w)/2$-component unlink with disk gives a cobordism in $S^3 \times I$ consisting of $(k-w)/2$ annuli together with $L' \times I$, where $L' = L(w,w)$. Let $F'$ be a compact, connected, oriented surface properly embedded in $B^4$ with $\partial F' = L'$ and $h_4(L') = \beta_1(F')$. Then in the 4-ball $(S^3 \times I) \cup B^4$, where $S^3 \times \{1\}$ is identified with $\partial B^4$, $L$ bounds the disjoint union of $(L' \times I) \cup F'$ and $(k - w)/2$ annuli. Tubing these together gives a connected surface $F$ with $\partial F = L$. Hence 
$$h_4(L) \le \beta_1(F) = \beta_1(F') + (k - w) = h_4(L') + (k - w)$$ 
Since $L$ and $L'$ are fibred, $h(L) = \deg \Delta_L$ and $h(L') = \deg \Delta_{L'}$. Then, since $h(L') \ge h_{4}(L')$, we have
$$h(L) - h_{4}(L) \ge (h(L) - h(L')) + (h_{4}(L') - h_{4}(L)) \ge (\deg \Delta_L - \deg \Delta_{L'}) - (k - w)$$
which by Lemma \ref{lemma: diff deg} is strictly positive unless $g(L) = 0$.
  \end{proof}

\subsubsection{The case \texorpdfstring{$w = k, L \notin \mathcal P$}{w = k, L \not\in \mathcal P}}
\label{subsubsec: w = k not P}

Since we assume that $L\notin \mathcal{P}$ (see (\ref{equ: P})), we only need to consider links of the form $L(p,q;k,k;-)$, $L(p,q;k,k;-,+)$, or $L(p,q;k,k;-,-)$. In this section we determine which links of these forms satisfy $g_4(L) = g(L)$. The results are stated in the following three propositions.

\begin{prop}
  \label{prop: III -}
  Let $L = L(p,q;k,k;-)$. Then $g_4(L) < g(L)$.
\end{prop}

\begin{prop}
\label{prop: IV +}
Let $L = L(p,q;k,k;-,+)$. Then $g_4(L) < g(L)$.
\end{prop}
\begin{prop}
\label{prop: IV -}
Let $L = L(p,q;k,k;-,-)$. Then $g_4(L) = g(L)$ if and only if $L = L(2, 3; 1,1; -, -)$.
\end{prop}

  Our approach is the following. Consider the case in Proposition \ref{prop: III -}. We band $C_1$ to a $C_{p,q}$ component $C$ to get a link $L'$ with $|L'| = |L| - 1 = k$. After some manipulation we show that $L'$ is the closure of a braid with braid index $b$ and $c$ crossings, where we explicitly compute $b$ and $c$ as functions of the parameters $p,q$ and $k$. Let $L''$ be the link $L \setminus (C_1 \cup C)$. The band move defines a cobordism $F''$ in $S^3 \times I$ from $L \times \{0\}$ to $L' \times \{1\}$, where $F''$ is the disjoint union of $L'' \times I$ and a pair of pants $P$ with $P \cap (S^3 \times \{0\}) = (C_1 \cup C) \times \{0\}$.

  Let $F'$ be the Seifert surface for $L'$ defined by the braid representative; so $\chi(F') = b - c$. Then $F = F'' \cup (F' \times \{1\}) \subset S^3 \times I \subset B^4$ is an oriented surface with $\partial F = L$. Therefore
  $$h_4(L) \le \beta_1(F) = \beta_1 (F') + 1 = 2 - \chi(F') = 2 - b + c$$
  Hence
\begin{equation} 
\label{eqn: 1}
h(L) - h_4(L) \ge \deg \Delta_{L} -2 + b - c
\end{equation}
We show that this is positive .

For Proposition \ref{prop: IV +} we do the same construction, while for Proposition \ref{prop: IV -}, in order to get a braid diagram, we first band $C_2$ to $C$ and then band $C_1$ to the resulting component $C'$.

\begin{proof}[Proof of Proposition \ref{prop: III -}]
\begin{figure}[ht]
    \centering
   \includegraphics[scale=1.2]{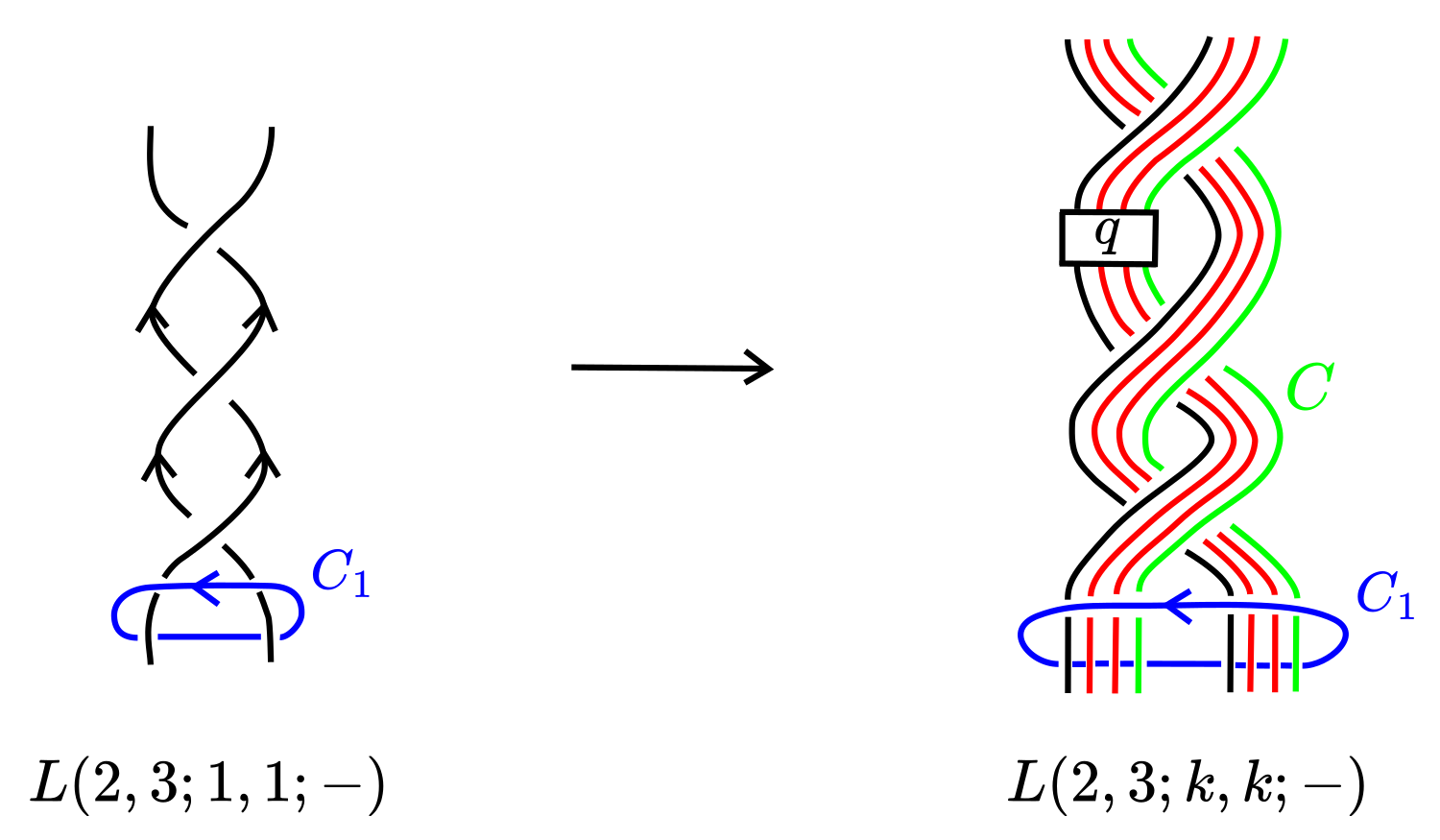}
    \caption{The link $L(p, q; 1, 1; -)$ with $p=2, q= 3$ and the diagram $D$ of $L(p, q; k, k; -)$ when $p=2, q= 3$ and $k=4$.  In the figure, the box with the letter $q$ inside represents $q$ full twists. }
    \label{fig: diagram D}
\end{figure}

$L(p,q;1,1;-)$ has a diagram as shown in Figure \ref{fig: diagram D}, which illustrates the case $p = 2$, $q = 3$, where the braided subdiagram is understood to be closed off in the usual way. 


 A diagram $D$ of $L= L(p,q;k, k; -)$ is obtained from this by taking $k$ parallel (with respect the ``blackboard'' framing) copies of each strand of the braid, and then inserting $q$ full twists in these strands, as indicated in Figure \ref{fig: diagram D} (where $k = 4$). (Before putting in these twists, the number of crossings in the diagram between two copies of $C_{p,q}$ is $2q(p-1)$, so their linking number is $q(p-1)$; putting in the twists gives the correct linking number $pq$.)

The diagram $D$ has $k^{2}q(p-1)$ crossings between the $C_{p,q}$ components, $2kp$ crossings of $C_1$ with the $C_{p,q}$ components, and $qk(k-1)$ crossings coming from the twists, making a total of $c(D) = k^{2}pq + 2kp - kq$ crossings.

Let $C$ be the rightmost copy of $C_{p,q}$, and band $C_1$ to $C$ as shown in  Figure \ref{fig: 4}. This gives a link $L'$ with $|L| - 1 = k$ components. Let $C'$ be the component of $L'$ that is created by banding $C_1$ to $C$.
\begin{figure}[ht]
    \centering
   \includegraphics[scale=1.2]{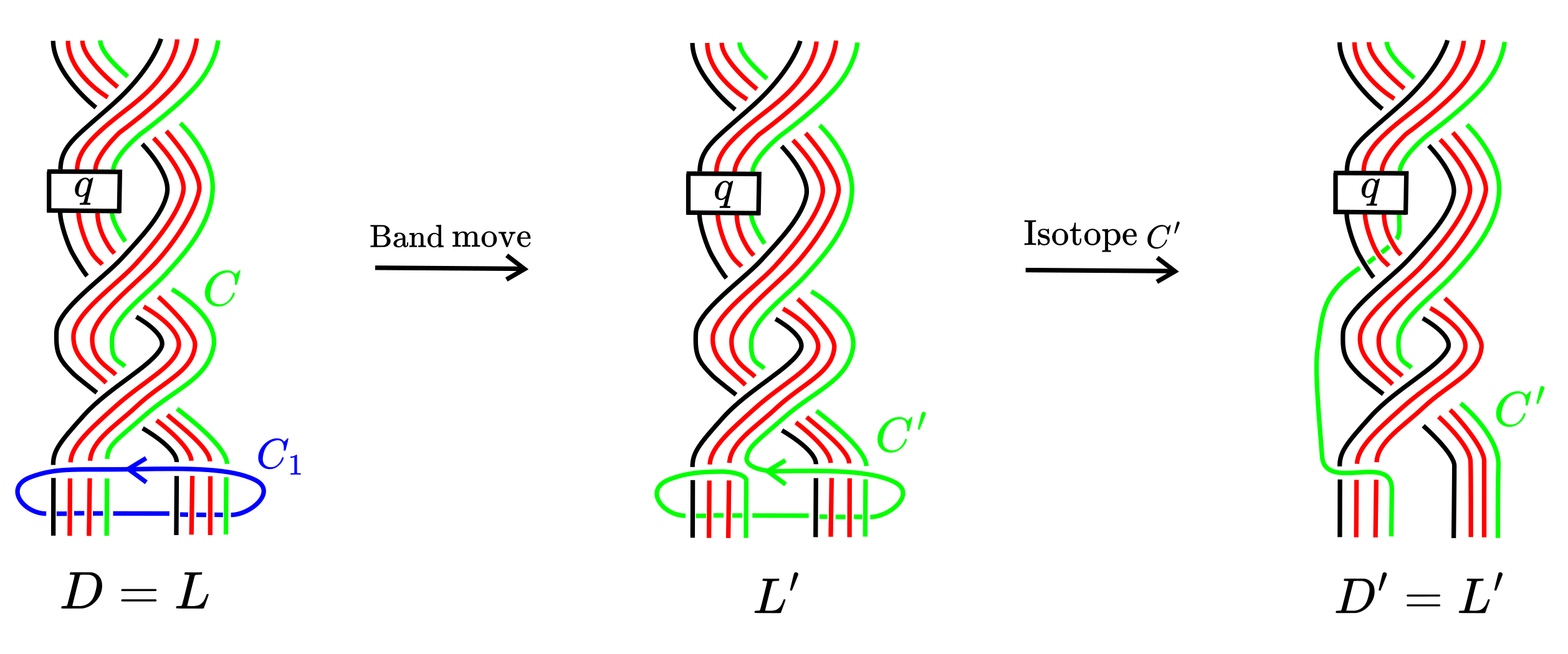}
    \caption{Band $C_1$ to $C$ which results in the component $C'$ in $L'$. The diagram also shows how to isotope $C'$ in the case that $p<q$.    }
    \label{fig: 4}
\end{figure}
Clearly $L'$ has a braid diagram with braid index $kp$. Further manipulations will give a braid index $kp$ diagram which has $l$, say, fewer crossings than $D$. Then (\ref{eqn: 1}) gives
$$h(L) - h_4(L) \ge \deg \Delta_L -2 + kp -c(D) + l$$
By Lemma \ref{lemma: deg delta L} $\deg \Delta_L = 1 + (kp-1)(kq-1)$. Therefore, recalling the formula for $c(D)$, we have
$$\deg \Delta_L -2 + kp - c(D) = -2kp$$
giving
\begin{equation}
\label{eqn: 2}
h(L) - h_4(L) \ge l - 2kp
\end{equation}

To compute $l$ we distinguish two cases.

\noindent
{\bf Case 1.} $p < q$.

In this case, we can further isotope the component $C'$ as shown in Figure \ref{fig: 4}. We denote the new diagram of $L'$ by $D'$. We claim that going from $D$ to $D'$, there is a total loss of $4k(p-1) + 2k$ crossings. We use Figure \ref{fig: 4 old} to demonstrate the computation. 

In Figure \ref{fig: 4 old}, the black arcs together with the red arc $A$ is the component $C'$ of the diagram $D'$ for  $p = 5, q = 8$ (cf. Figure \ref{fig: 4}).  As in Figure \ref{fig: 4}, we can isotope $C'$ by isotoping the subarc $A$ rel endpoints to the arc $B$ colored blue. Note that for each crossing marked by red or green dots in Figure \ref{fig: 4 old}, and for the crossing of $D$ where $C_1$ passed over $C$ at which the band move is made, we lose $k$ crossings in going from $D$ to $D'$. This results in a total loss of $4k(p-1) + 2k$ crossings.

\begin{figure}[ht]
    \centering
   \includegraphics[scale=0.5]{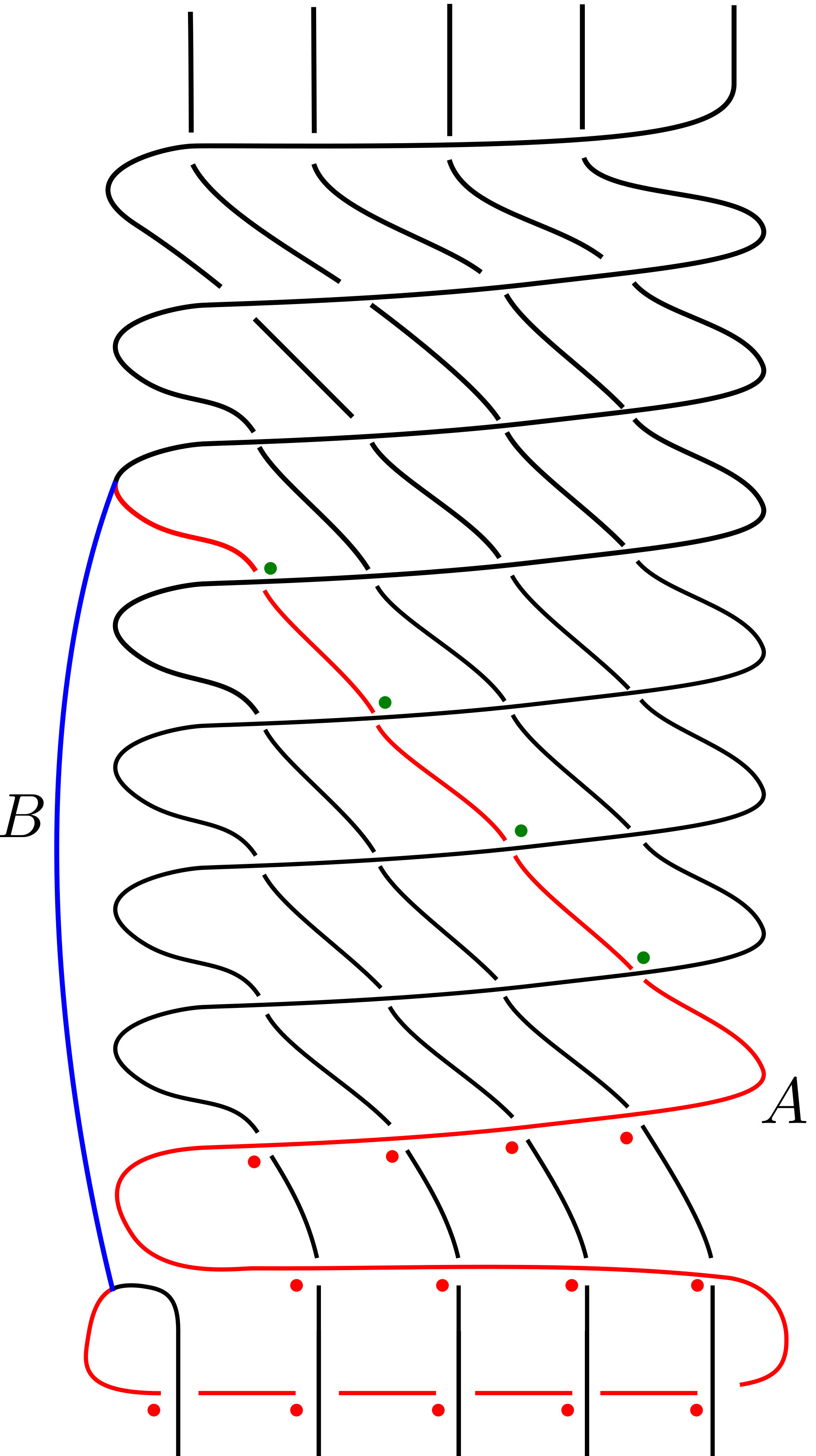}
    \caption{}
    \label{fig: 4 old}
\end{figure}

Next, the bottom left $(k-1)$ crossings in $D'$ as shown in Figure \ref{fig: 4}) can be pushed downwards and along the diagram until they appear just above the $q$ full twists, so that there the diagram now looks like Figure \ref{fig: 7}. Finally, the diagram in Figure \ref{fig: 7} can be isotoped to look like Figure \ref{fig: 8}. Since a full twist in $r$ strands gives rise to $r(r-1)$ crossings, this final isotopy results in a further loss of $k(k-1) - (k-1)(k-2) = 2(k-1)$ crossings.

\medskip
\begin{figure}[ht]
\centering
  \begin{minipage}{0.4\textwidth}
    \centering
   \includegraphics[scale=1]{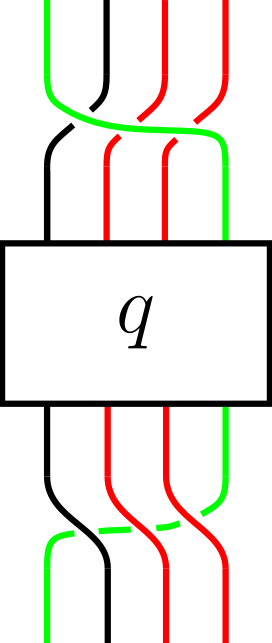}
    \caption{}
    \label{fig: 7}
  \end{minipage}
  \vspace{30pt}
  \begin{minipage}{0.4\textwidth}
    \centering
    \vspace{28pt}
   \includegraphics[scale=1.1]{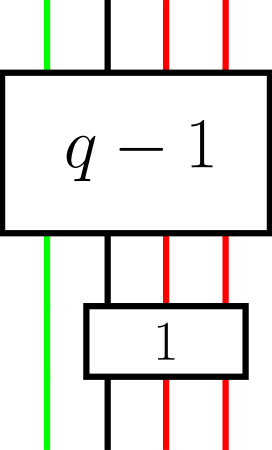}
    \caption{}
    \label{fig: 8}
  \end{minipage}
 \end{figure}

We end up with a braid diagram of $L'$, with braid index $kp$ and $4k(p-1) + 2k + 2(k-1) = 4kp -2 = 2(2kp-1)$ fewer crossings than $D$. Then (\ref{eqn: 2}) gives
$$h(L) - h_4(L) \ge 2(2kp-1) - 2kp = 2(kp-1)$$

Hence $h_4(L) < h(L)$ unless $k = p = 1$,  which is the negative Hopf link and is excluded by notational convention (see Remark \ref{rem: notation conventions}(1)). 

\noindent
{\bf Case 2.} $p > q$.

In this case, we do the same construction as in case (1).  Figure \ref{fig: 9} illustrates the isotopy we perform on the component $C'$ in the case $p = 7$ and $q = 4$.

\begin{figure}[ht]
    \centering
   \includegraphics[scale=0.4]{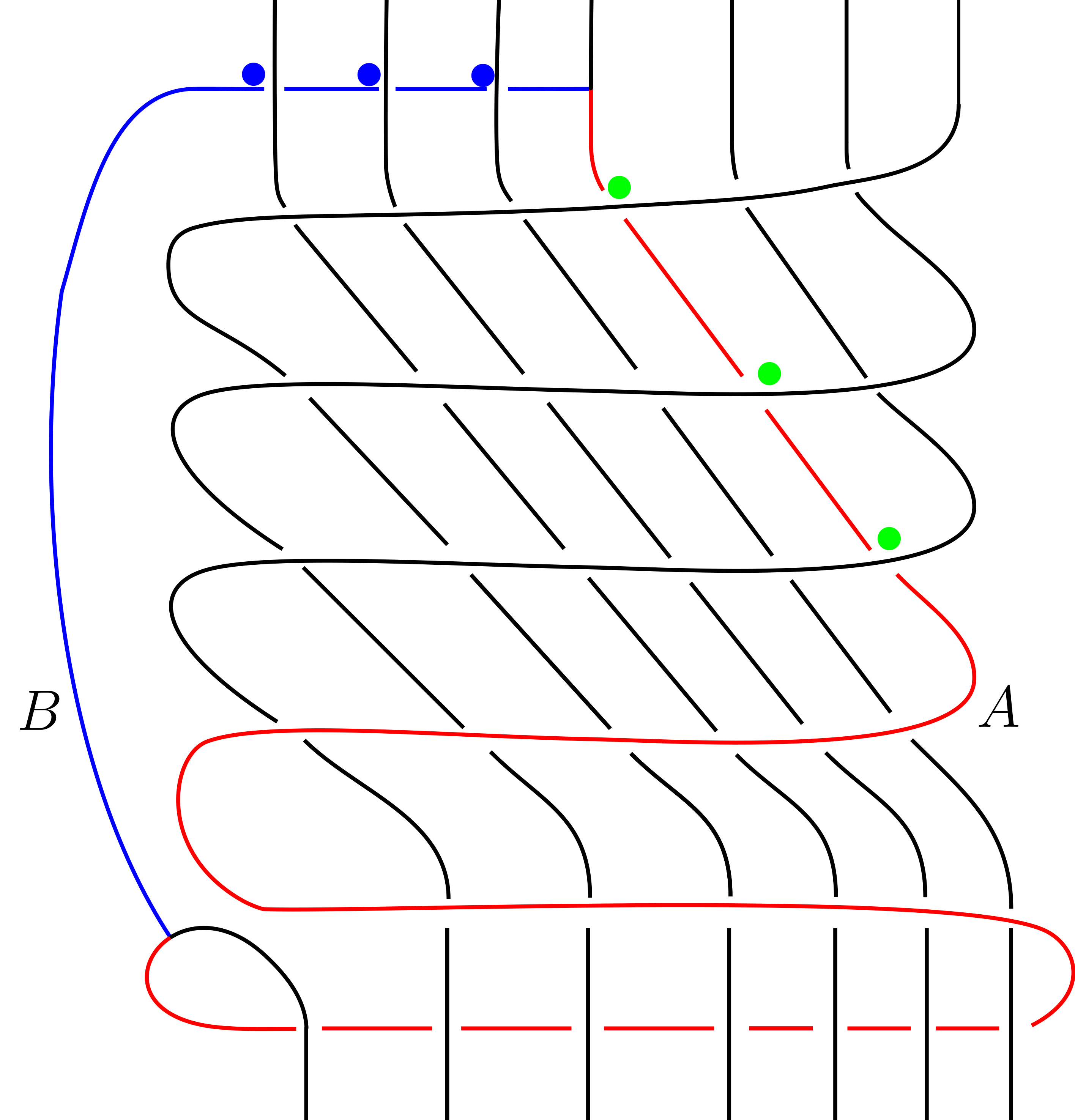}
    \caption{}
    \label{fig: 9}
\end{figure}

In Figure \ref{fig: 9}, the union of the black arcs and the red arc $A$ is the component $C'$ after banding $C_1$ to $C$ as before. Now when we isotope the red arc $A$ to the blue arc $B$, instead of losing the $k(p-1)$ crossings corresponding to those colored green in Figure \ref{fig: 4 old}, we only lose $k(q-1)$ crossings, corresponding to those marked by green dots in Figure \ref{fig: 9}. Furthermore, we gain an additional $k(p-q)$ crossings, at places marked by blue dots in Figure \ref{fig: 9}. Thus based on our computation in case (1) we get a braid diagram of $L'$ with braid index $kp$ and 
$$2(2kp-1) - k(p-1) + k(q-1) - k(p-q) =  2(2kp-1) - 2k(p-q)$$ fewer crossings than $D$. Equation (\ref{eqn: 2}) then gives
$$h(L) - h_4(L) \ge 2(2kp-1) - 2k(p-q) - 2kp = 2(kq-1)$$

Since $q > 1$, this is positive.
\end{proof}

\begin{proof}[Proof of Proposition \ref{prop: IV +}]

\begin{figure}[ht]
    \centering
   \includegraphics[scale=0.4]{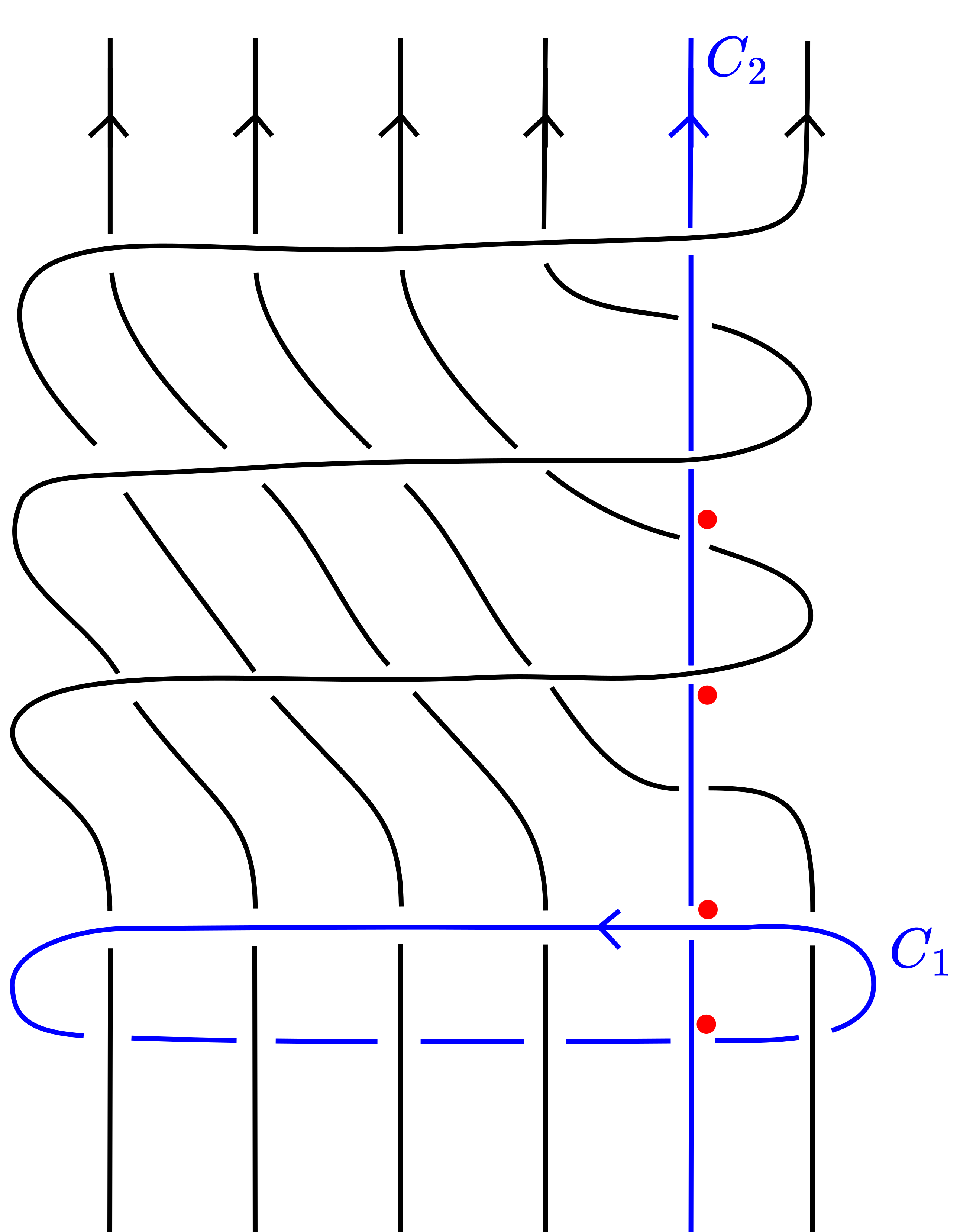}
    \caption{}
    \label{fig: 10}
\end{figure}

 Figure \ref{fig: 10} shows a diagram of $L(p,q;1,1;-,+)$, with $p = 5, q = 3$; a diagram $D$ of $L$ is obtained from this by taking $k$ parallel copies of the $C_{p,q}$ component and inserting $q$ full twists. The crossings of $D$ correspond to the crossings of the diagram $D$ in the proof of Proposition \ref{prop: III -}, together with the crossings of $C_2$ with the $C_{p,q}$ components and the crossings of $C_2$ with $C_1$. Thus
  $$c(D) = (k^{2}pq + 2kp -kq) + 2kq + 2$$
  $$= k^{2}pq + 2kp + kq +2$$

  Carrying out the same procedure as in the proof of Proposition \ref{prop: III -}, we get a link $L'$ with $|L'| = k + 1$, having a braid diagram of braid index $kp + 1$ with $l$, say, fewer crossings than $D$. By Lemma \ref{lemma: deg delta L} $\deg \Delta_L = 1 + k(kpq - p + q)$, and hence
  $$\deg \Delta_L -2 + (kp+1) - c(D) = -2(kp+1)$$
  Then (\ref{eqn: 1}) gives

\begin{equation}
\label{eqn: 3} 
h(L) - h_4(L) \ge l - 2(kp+1)
\end{equation} 
 Note that $l$ is the number of crossings lost in the proof of Proposition \ref{prop: III -} plus the 4 crossings with $C_2$ marked by the red dots in Figure \ref{fig: 10}. Thus, if $p < q$ then $l = 2(2kp-1) + 4 = 2(2kp+1)$, and if $p > q$ then $l = 2(2kp-1) -2k(p-q) + 4$. This gives $h(L) - h_4(L) \ge 2kp$ or $2kq$, respectively.  
\end{proof} 

\begin{proof}[Proof  of Proposition \ref {prop: IV -}]

  Here we first band $C_2$ to a $C_{p,q}$ component $C$, and then band $C_1$ to the resulting component $C'$. This gives a link $L'$ with $|L'| = |L| - 2 = k$, which we show has a braid diagram $D'$ of braid index $b = kp-1$ with $c = c(D')$ crossings. Let $L'' = L \setminus (C_1 \cup C_2 \cup C)$. We then get a cobordism in $S^3 \times I$ from $L \times \{0\}$ to $L' \times \{1\}$ consisting of $L'' \times I$ and a 4-punctured sphere $P$ with $P \cap (S^3 \times \{0\}) = (C_1 \cup C_2 \cup C) \times \{0\}$. Arguing as in the previous case we get
  $$h_4(L) \le 3 - b + c$$
  If $c = c(D) - l$, then we have
  $$h(L) - h_4(L) \ge \deg \Delta_L - 3 + (kp-1) -c(D) + l$$

  By Lemma \ref{lemma: deg delta L} $\deg \Delta_L = 1 + k(kpq - p - q)$. Since as an unoriented diagram $D$ is the same as the diagram $D$ in the proof of the previous proposition, $c(D) = k^{2}pq + 2kp + kq +2$ and we then easily compute that
$$\deg \Delta_L - 3 + (kp-1) - c(D) = -2k(p+q) - 5$$
  giving
\begin{equation}
\label{eqn: 4} 
h(L) - h_4(L) \ge l - 2k(p+q) - 5
\end{equation}

We now give a detailed description of the bandings and the diagram manipulations that result in the braid diagram $D'$, enabling us to compute $l$.

Recall that for the family of links $L(p,q;k,k;-,-)$, because of the symmetry in $p$ and $q$ we often assume $p < q$ (see Remark \ref{rem: notation conventions}(2)). However, because of our conventions on Seifert link diagrams it is convenient in the present proof to interchange $p$ and $q$ and assume $p > q$.

\begin{figure}[ht]
\centering
  \begin{minipage}{0.3\textwidth}
    \centering
   \includegraphics[scale=0.3]{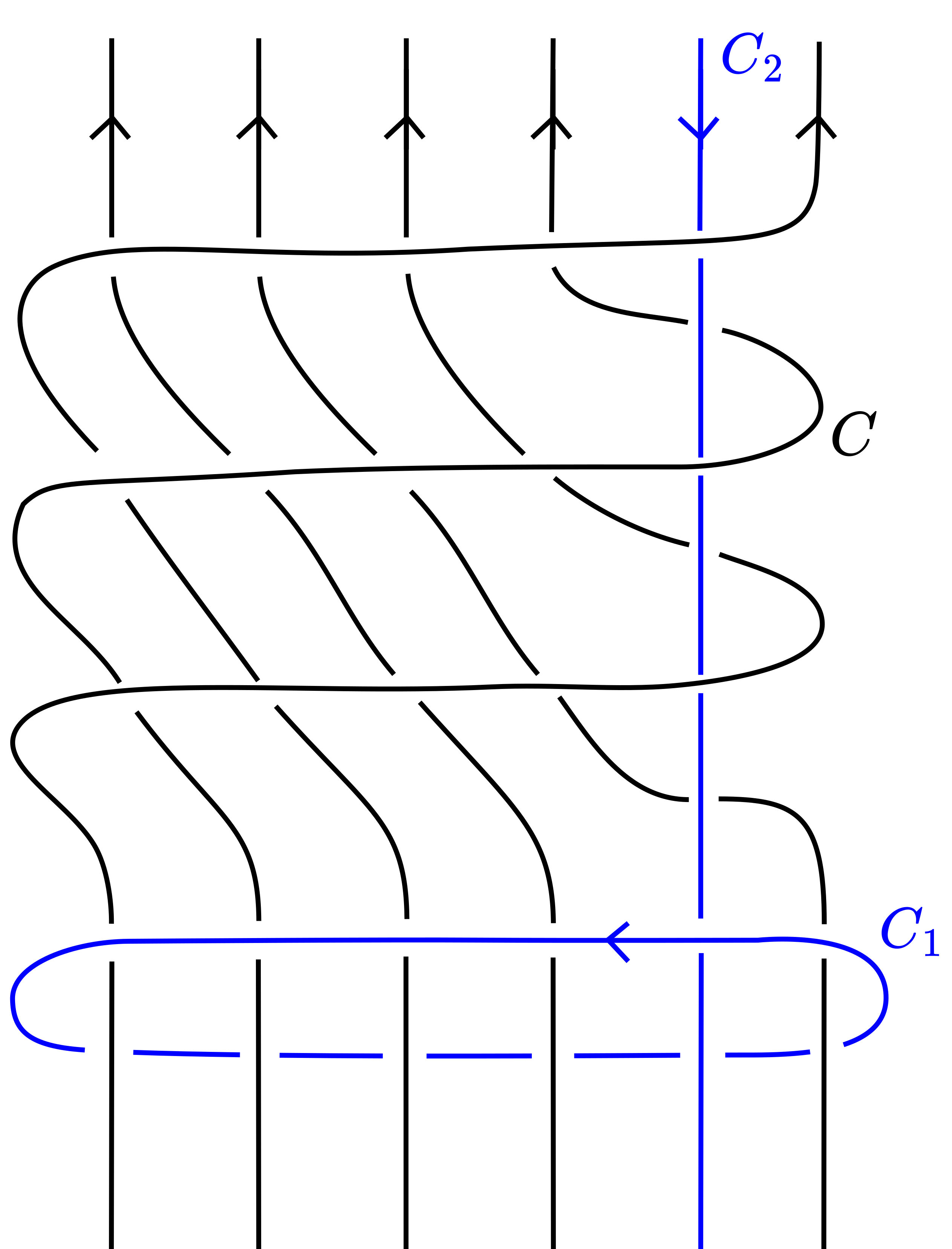}
    \caption{}
    \label{fig: 191}
  \end{minipage}
  \begin{minipage}{0.3\textwidth}
    \centering
   \includegraphics[scale=0.3]{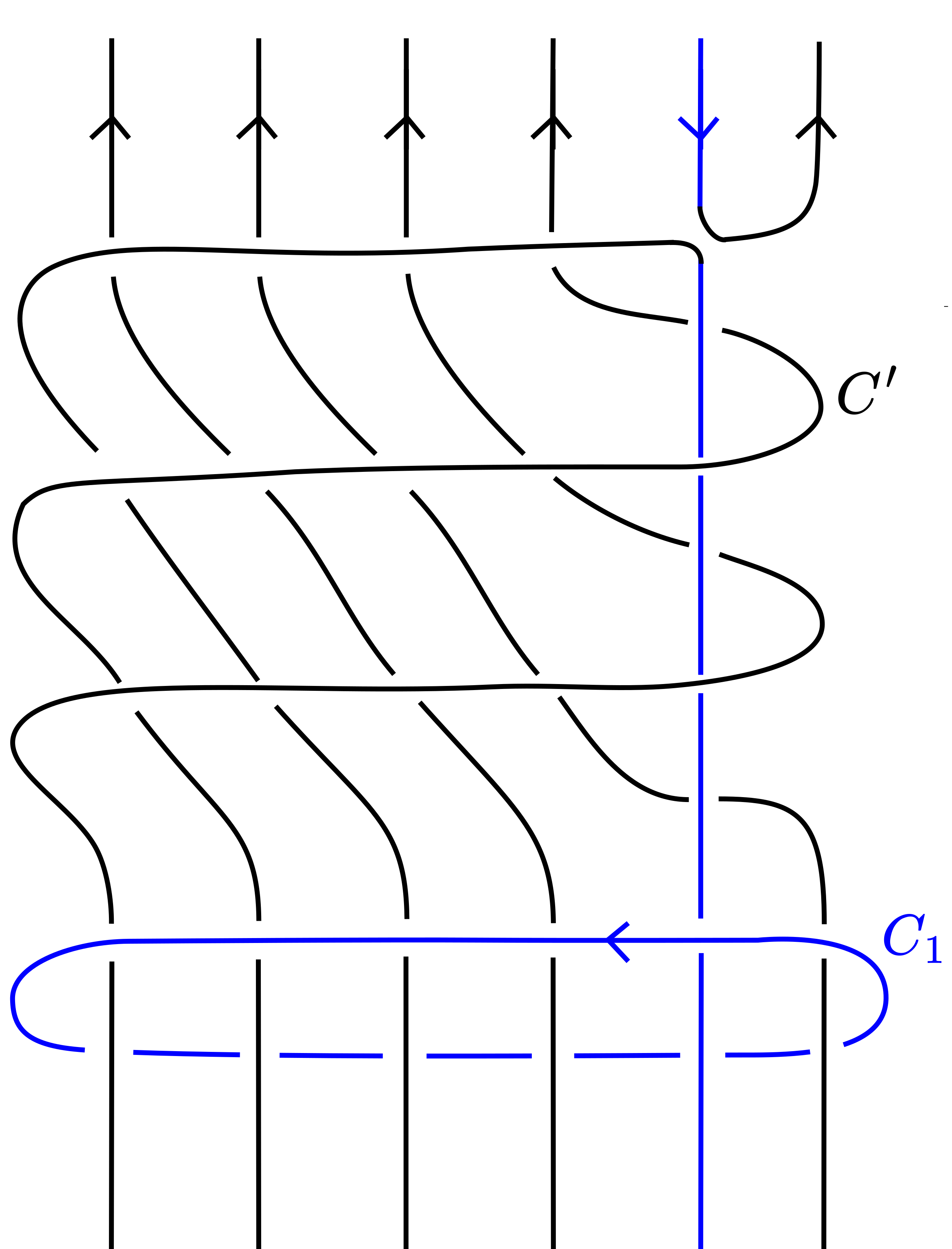}
    \caption{}
    \label{fig: 192}
  \end{minipage}
   \begin{minipage}{0.3\textwidth}
    \centering
   \includegraphics[scale=0.3]{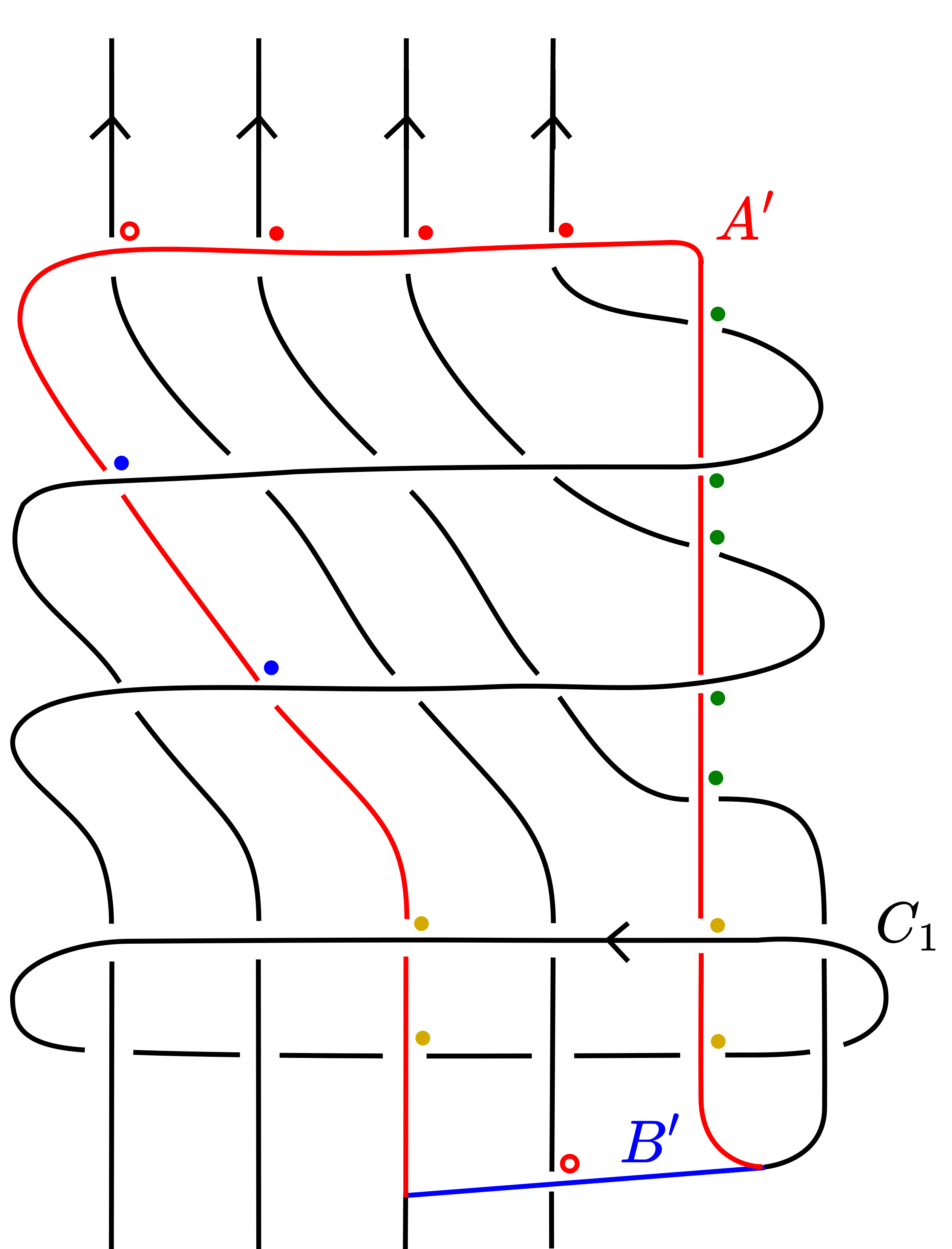}
    \caption{}
    \label{fig: 193}
  \end{minipage}
 \end{figure}
 
We start with a diagram $D$ of $L=L(p,q;k,k;-,-)$, which is obtained from the diagram of $L(p, q;1,1; -, -)$ (cf. Figure \ref{fig: 191}) by taking $k$ parallel copies of the $C_{p,q}$ component and inserting $q$ full twists as before. This diagram is the same as the diagram $D$ of $L$ in the proof of Proposition \ref{prop: IV +} except that the orientation of $C_2$ is reversed.  Let $C$ be the rightmost $C_{p,q}$ component in $D$; Figure \ref{fig: 191} shows the component $C$ together with $C_1$ and $C_2$. Band $C_2$ to $C$ as shown in Figure \ref{fig: 192} and call the resulting component $C'$; Figure \ref{fig: 192} shows the component $C'$ and $C_1$ after the band move.

Next, we isotope $C'$ by first moving the arc in the upper right corner of Figure \ref{fig: 192} to the bottom right. We then  replace the subarc $A'$ with $B'$ as shown in Figure \ref{fig: 193}. Let $D_1$ be the resulting diagram as shown in Figure \ref{fig: 194}. Then $c(D_1) = c(D) - l_1$ for some $l_1$.

To compute $l_1$, note that in going from $D$ to $D_1$ we lose the following numbers of crossings:

\begin{itemize}
\setlength\itemsep{0.3em}  
\item $kq$, corresponding to crossings marked by the solid red dots in Figure \ref{fig: 193}. (The loss of crossings  of $A'$ at the locations marked by the open red circle is cancelled with the gain of the crossings between the arc $B'$ and $C'$,  also marked by an open red circle in Figure \ref{fig: 193}.) 
\item $k(q-1)$, corresponding to crossings marked by the blue dots,
\item  $2k(q-1) + k$, corresponding to crossings marked by the green dots, 
\item  the 4 crossings  marked by yellow dots,
\item  the crossing that was eliminated in the band move.
\end{itemize}
Adding these gives $l_1 = 2k(2q-1) + 5$.

Now we band $C_1$ to $C'$  as shown in Figure \ref{fig: 195}. Then isotope the arc $A$ rel endpoints to $B$ in Figure \ref{fig: 196}. This gives a diagram $D'$ with $c(D') = c(D_1) - l_2$.
  
  \begin{figure}[ht]
\centering
  \begin{minipage}{0.3\textwidth}
    \centering
   \includegraphics[scale=0.3]{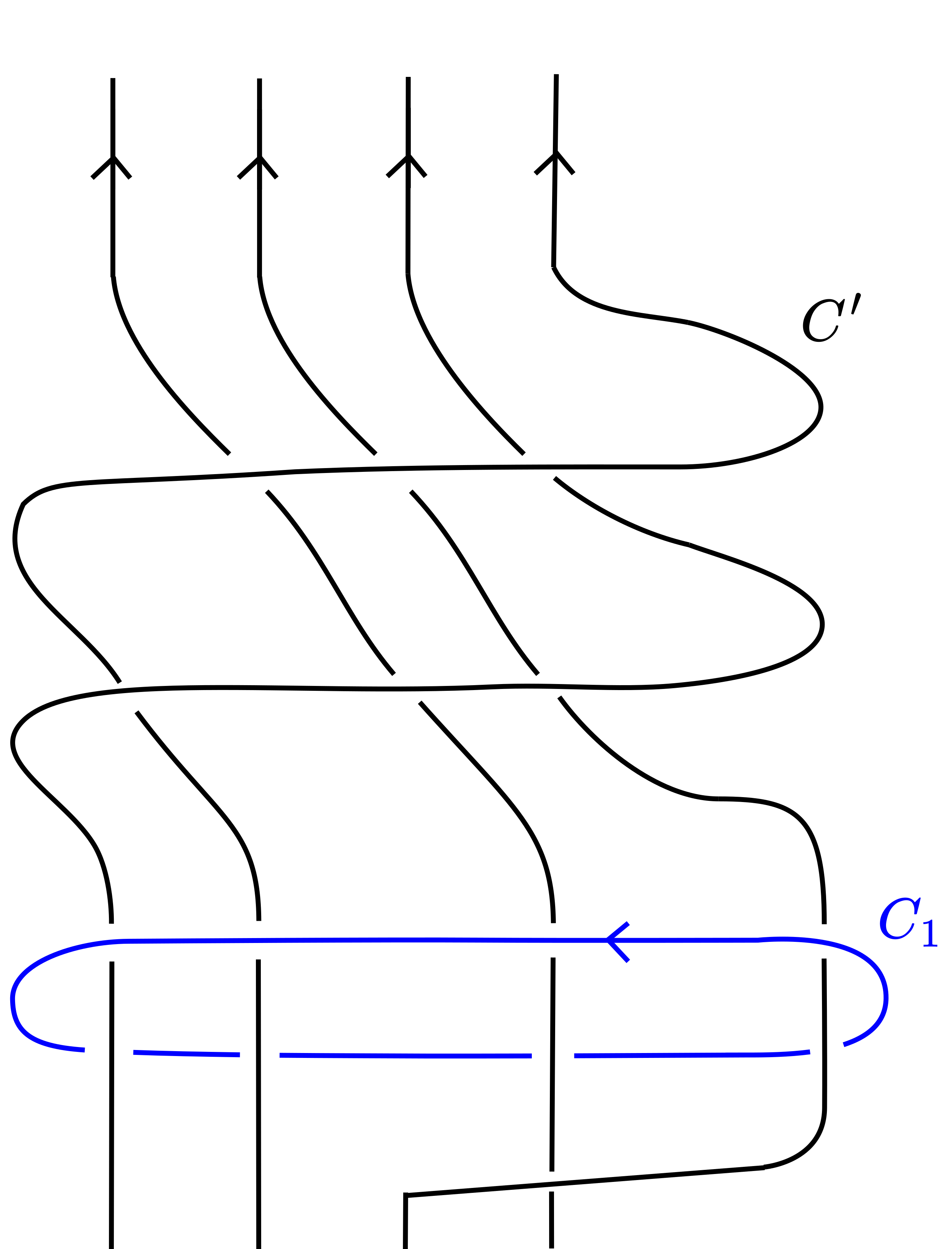}
    \caption{}
    \label{fig: 194}
  \end{minipage}
  \begin{minipage}{0.3\textwidth}
    \centering
   \includegraphics[scale=0.3]{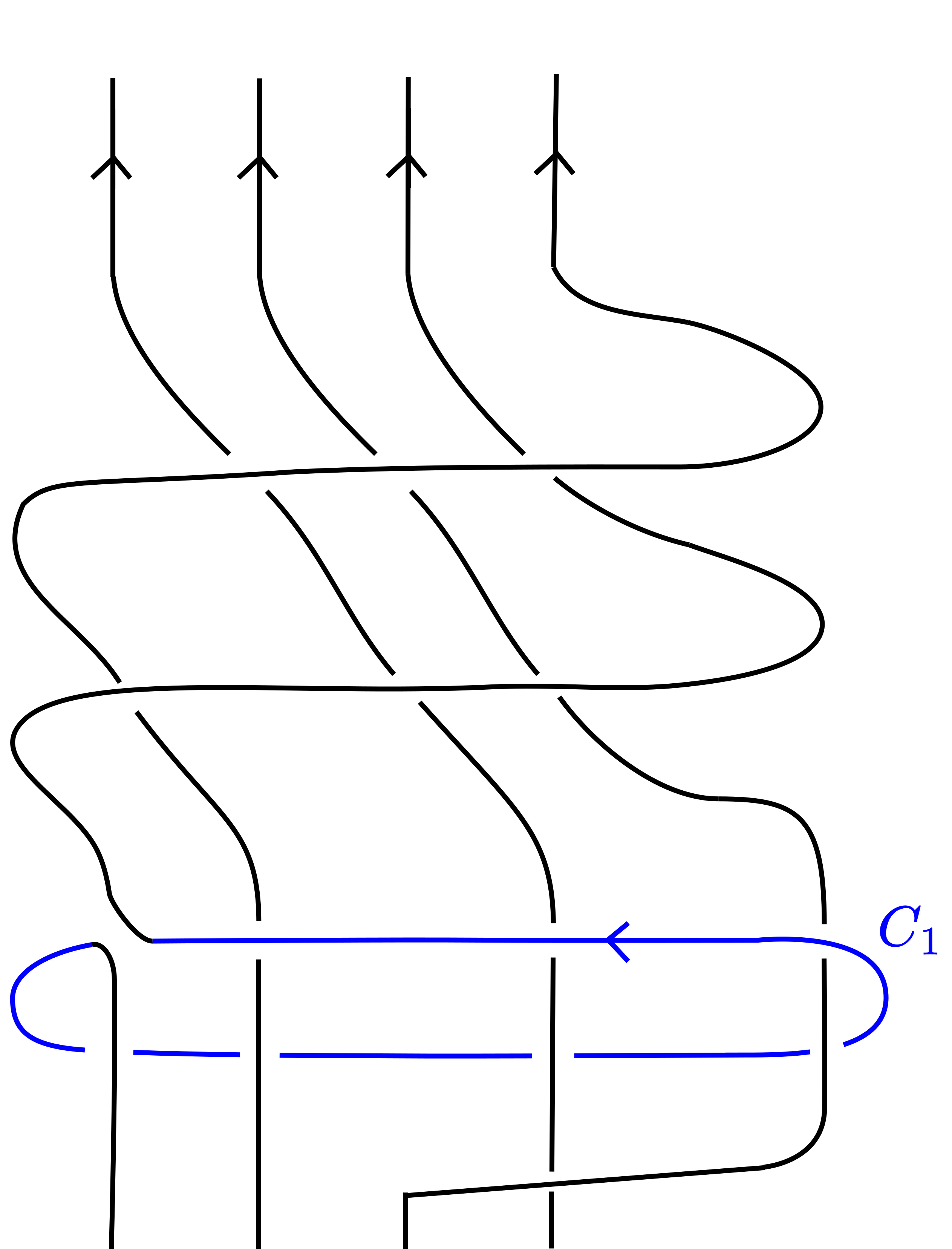}
    \caption{}
    \label{fig: 195}
  \end{minipage}
   \begin{minipage}{0.3\textwidth}
    \centering
   \includegraphics[scale=0.3]{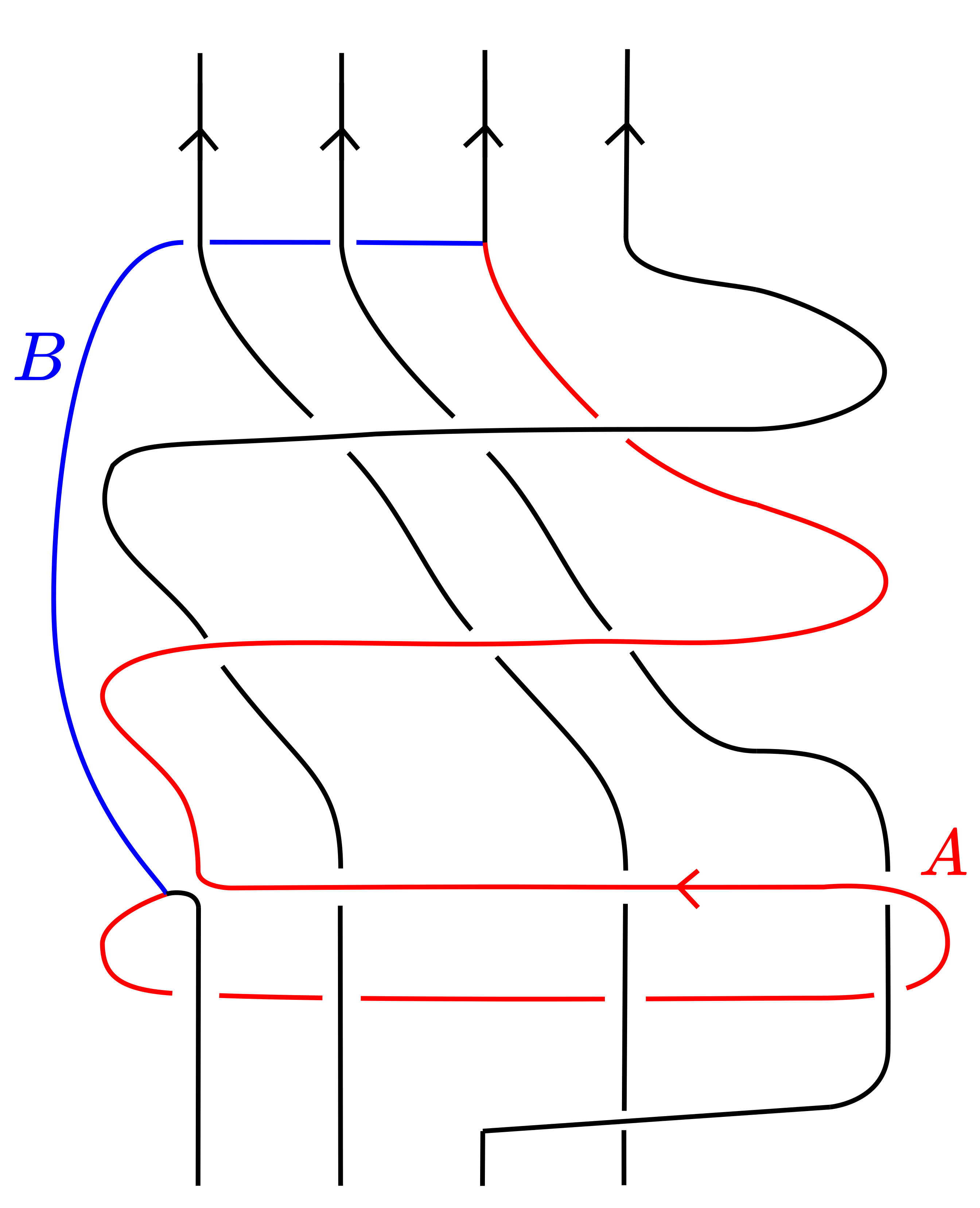}
    \caption{}
    \label{fig: 196}
  \end{minipage}
 \end{figure}

  The computation of $l_2$ is similar to that in the proof of the previous proposition in the case $p > q$. We get
  $$l_2 = 3(k(p-1) - 1) + k + k(q-1) + k + 2(k-1) -k(p-q) = 2k(p+q) - 6$$
  Then (\ref{eqn: 4}) gives
  $$h(L) - h_4(L) \ge 2k(2q-1) - 6$$

  which is positive unless $k = 1$ and $q = 2$.

  When $k = 1$ and $q = 2$ it is easy to see that $L'$ is the unknot. Therefore $L$ bounds a pair of pants in $B^4$, giving $h_4(L) = 2$, i.e. $g_4(L) = 0$. Also, $h(L) = \deg \Delta_L = p - 1$ by Lemma \ref{lemma: deg delta L}, and hence $h(L) - h_4(L) = p - 3$. So $g(L) > g_4(L)$ unless $k=1$, $q=2$, and $p=3$, in which case $g(L) = g_4(L) = 0$. (Note that this is consistent with Theorem \ref{thm: g = 0}.)
\end{proof}

 \begin{proof}[Proof  of Theorem \ref{thm: g_4 = g results}]
 The ``if'' direction is clear.

 For the converse, suppose $L$ is a Seifert link with $g_4(L) = g(L) \ne 0$. Then $L$ is not a connected sum of Hopf links, so $L$ is $0$-,  $1$-, or  $2$-core.

 By Propositions \ref{prop: w = 0} and \ref{prop: 0 < w} we have $w = k$. Then, by Propositions \ref{prop: III -}, \ref{prop: IV +} and \ref{prop: IV -}, either $L \in \mathcal P$ or $L = L(2,3;1,1;-,-)$. But if the latter holds then $g(L) = 0$ by Theorem \ref{thm: g = 0}.
 \end{proof}

\subsection{Strongly quasipositive Seifert links}
\label{subsec: sqp sl}

Recall that we abuse terminology by saying that a link $L$ is strongly quasipositive if either $L$ or $-L$, the mirror of $L$, is strongly quasipositive in the usual sense. In this section it will be necessary to consider orientations more carefully, so if a link $L$ is strongly quasipositive in the usual sense we will say it is {\em strictly} strongly quasipositive.

\begin{lemma}
\label{lemma: 0 framed}
  Let $L$ be a fibred link with fibre $F$. Suppose $F$ contains a simple loop $C$ such that the framing of $C$ induced by $F$ is the $0$-framing and $C$ is unknotted in $S^3$. Then $L$ is not strongly quasipositive.
\end{lemma}

\begin{proof}
  An annular neighborhood $A$ of $K$ in $F$ is unknotted and untwisted in $S^3$, and therefore $A$ is not quasipositive. Since $A$ is a full subsurface of $F$, $F$ is not quasipositive \cite{Ru2}, and hence $L$ is not strictly strongly quasipositive. Since this argument also applies to $-L$, $L$ is not strongly quasipositive.
\end{proof}

Let $TB(K)$ be the maximal Thurston-Bennequin number of $K$. If $K$ is a component of link $L$ with at least two components, $l(K)$ will denote $\mbox{lk} (K, L \setminus K)$.

\begin{lemma}
\label{lemma: tb}
Let $K$ be a component of a non-split strictly strongly quasipositive link $L$. Then $-l(K) \le TB(K)$. In particular, if $K$ is unknotted then $l(K) \ge 1$.
\end{lemma}

\begin{proof}
Let $F$ be a compact, oriented surface in $S^3$ with $\partial F = L$, and let $A$ be a collar neighborhood of $K$ in $F$. Then $A$ is an annulus with core isotopic to $K$ and $-l(K)$ positive twists. The result then follows from \cite{Ru2}.
\end{proof}

\begin{prop}
\label{prop: sqp g = 0} 
  If $L$ is a Seifert link with $g(L) = 0$ then $L$ is strongly quasipositive if and only if $L$ is either $\#_k H_+$, $k \ge 1$, or of the form $L(p,q;k,0)$.
\end{prop}

\begin{proof} 
  
  First we note that a connected sum $L$ of Hopf links is strongly quasipositive if and only if the summands all have the same sign. For if $L$ has both an $H_+$ and an $H_-$ summand, then it is easy to see that the planar fibre  of $L$ contains a curve $C$ as in Lemma \ref{lemma: 0 framed}.

  To complete the proof we need to show that of the links listed in Theorem \ref{thm: g = 0}, none of those in the second family in part (2), or those in parts (3) or (4), are strongly quasipositive.

\begin{itemize}[leftmargin=*] 
\setlength\itemsep{0.4em}  
\item $L = L(1,1;k,1)$.  Let $C$ be a positively oriented component of $L$. Then $C$ is unknotted and $l(C) = 0$. Therefore $L$ is not strongly quasipositive by Lemma \ref{lemma: 0 framed}.
\item $L = L(1,q;k,0;\varepsilon)$. $C_1$ is unknotted and $l(C_1) = 0$, and hence $L$ is not strongly quasipositive by Lemma \ref{lemma: 0 framed}.
\item $L = L(p,p+1;k,0;-,+)$. First note that since $l(C_1) = -1$, if $L$ is strongly quasipositive then by Lemma \ref{lemma: tb}, $-L$ is strictly strongly quasipositive. Let $C$ be a positively oriented $C_{p,p+1}$ component of $L$. Then $l(C) = -p(p+1) + 1$. In $-L$, $C$ is $T(-p,p+1)$, and $TB(T(-p,p+1)) = -p(p+1)$  by \cite[Theorem 4.1]{EH}. But $-p(p+1) + 1 \not\le -p(p+1)$, so $-L$ is not strictly strongly quasipositive by Lemma \ref{lemma: tb}.
\item   $L = L(2,3;k,1;-,-)$. Since $l(C_1) = -1$, if $L$ is strongly quasipositive then $-L$ is strictly strongly quasipositive. Let $C$ be a positively oriented $C_{2,3}$ component of $L$. Then $l(C) = -5$. But, as in the previous case, $TB(T(-2,3)) = -6$ and $-5 \not\le -6$, so $-L$ is not strictly strongly quasispositive.

\begin{figure}[ht]
\centering
\begin{minipage}{0.45\textwidth}
\centering
\includegraphics[scale=1.1]{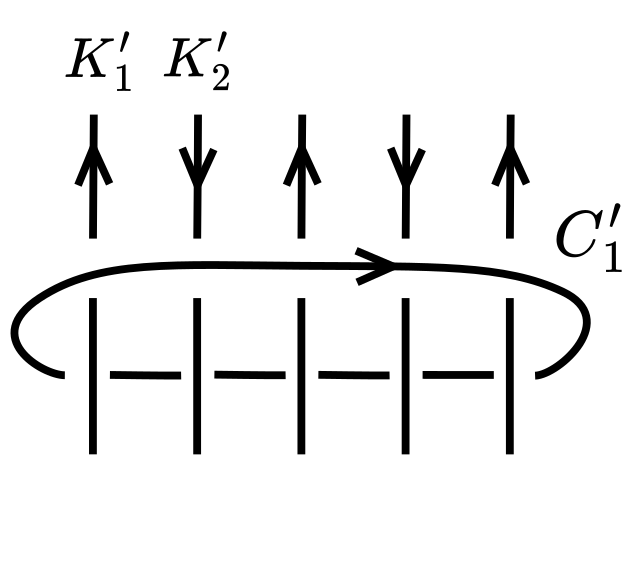}
\caption{}
\label{fig: A}
\end{minipage}
\begin{minipage}{0.45\textwidth}
\centering
\includegraphics[scale=1.1]{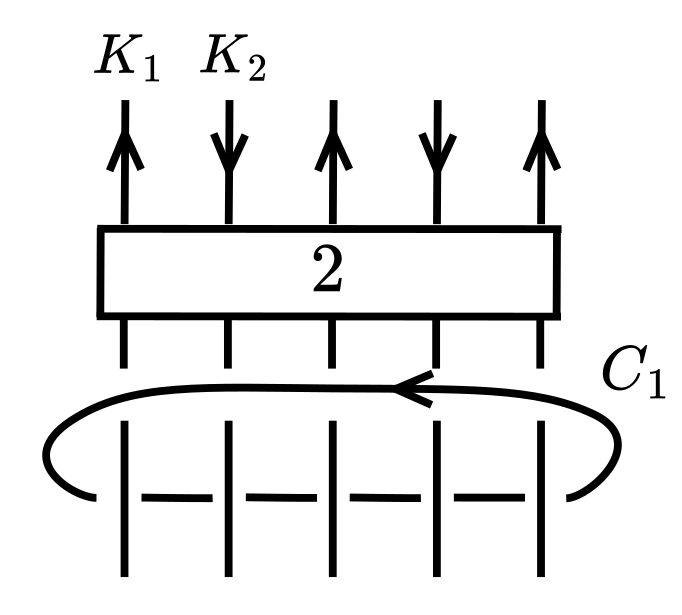}
\caption{}
\label{fig: B}
\end{minipage}
\end{figure}

\begin{figure}[ht]
\centering
\includegraphics[scale=1.1]{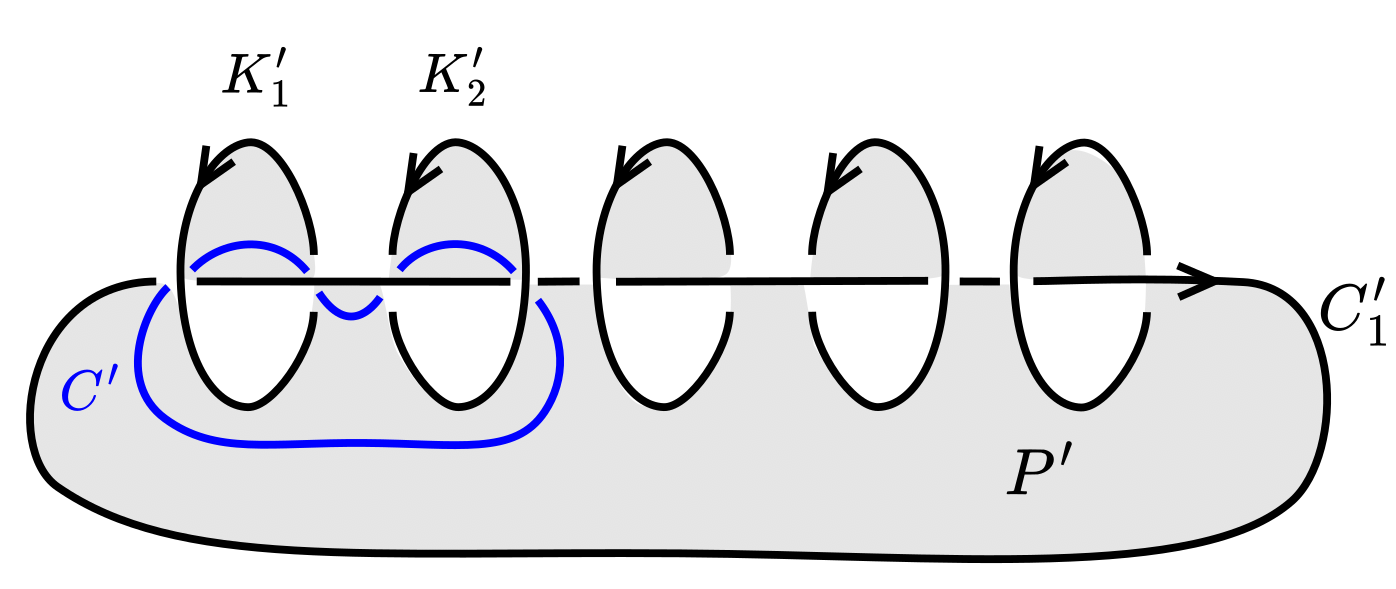}
\caption{}
\label{fig: C}
\end{figure}

\item   $L = L(1,2;k,1;-)$. Note that if $k=1$, then $L(1,2;1,1;-)$ is the negative Hopf link. So we assume that $k \ge 3$. Let $L'$ be the connected sum of $(k+1)/2$ copies of $H_{+}$ and $(k-1)/2$ copies of $H_-$ shown in Figure \ref{fig: A}, which illustrates the case $k = 5$. Then $L$ is obtained from $L'$ by putting two full positive twists in the $k$ braid strands and reversing the orientation of $C'_1$; see Figure \ref{fig: B}.  Let $V$ (resp. $V'$) be the solid torus that is the exterior of $C_1$ (resp. $C'_1$), and let $\mu, \lambda$ (resp. $\mu', \lambda'$) be the meridian and longitude of $C_1$ (resp. $C'_1$).  

\medskip 

\noindent $L'$ is fibred with fibre a planar surface $P'$, and the framing on $C'_1$ induced by $P'$ is $-l(C'_1) = -1$. Thus $[P' \cap \partial N(C'_1)] \in H_1(\partial N(C'_1))$ is $- \mu' + \lambda'$. There is a homeomorphism $h : V' \to V$ such that $h_{*}(\lambda') = - \lambda$, $h_{*}(\mu') = - \mu - 2 \lambda$, and $h(L' \setminus C'_1) = L \setminus C_1$. Hence if $P = h(P')$ then $[P \cap \partial N(C_1)] \in H_1(\partial N(C_1))$ is $h_{*}(- \mu' + \lambda') = \mu + \lambda$. Since $- l(C_1) = 1$, it follows that $L$ is fibred with fibre $P$.

\medskip 

\noindent Let $C' \subset P'$ be a simple loop as shown in  Figure \ref{fig: C}.  Note that $C'$ bounds a disk in $V'$. Hence $C = h(C')$ bounds a disk in $V$ and is therefore unknotted in $S^3$. In $P'$, $C'$ cobounds a pair of pants with $K'_1 \cup K'_2$, so in $P$, $C$ cobounds a pair of pants with $K_1 \cup K_2$. It follows that the framing on $C$ induced by $P$ is $- \mbox{lk} (C, K_1 \cup K_2) = - \mbox{lk} ( C', K'_1 \cup K'_2) = 0 $. Therefore $L$ is not strongly quasipositive by Lemma \ref{lemma: 0 framed}.

\end{itemize}
\end{proof}

\begin{proof}[Proof of Theorem \ref{thm: sqp seifert results intro}]

It suffices to prove (2).

If $L \in \mathcal P$ then $L$ is strongly quasipositive. If $L = L(p,q;k,0)$ then $L$ bounds an oriented surface $F$ that is a disjoint union of $k/2$ annuli. In \cite[Figures 19 and 20]{Ish}, this surface $F$ is explicitly shown to be quasipositive.

Conversely, suppose $L$ is a strongly quasipositive Seifert link, not of the form $L(p,q;k,0)$. Then $L$ is fibred by Proposition \ref{prop: basic props}(3). Therefore $g_4(L) = g(L)$ by Lemma \ref{lemma: conn qp}, and hence either $L \in \mathcal P$ or $g(L) = 0$ by Theorem \ref{thm: g_4 = g results}. But if $g(L) = 0$ then $L = \#_k H_+$, and hence $L \in \mathcal P$, by Proposition \ref{prop: sqp g = 0}.
 \end{proof}

\subsection{Definite Seifert links} 
\label{subsec: def sl}

In this section we prove Theorem \ref{thm: def results}.

We first determine the definite Seifert links with genus zero.

\begin{prop}
\label{prop: def g = 0}  

If $L$ is a Seifert link with $g(L) = 0$ then $L$ is definite if and only if it is either
\begin{enumerate}[leftmargin=*] 
\setlength\itemsep{0.3em}  
\item[{\rm (1)}] $\#_k H_+$
\item[{\rm (2)}] $L(p,q;2,0)$
\item[{\rm (3)}] $L(2,3;1,1;-,-)$. 
\end{enumerate}
\end{prop}

The proof will use the following 4-dimensional generalization of the inequality $|\sigma(L)| \leq h(L)$, due to Murasugi \cite{Mu} and Tristram \cite{Tri}. 
 
\begin{thm} {\rm (Murasugi-Tristram Inequality)}
\label{thm: mt ineq background} 
Let $L$ be a link and $F$ a compact, oriented surface with $\partial F = L$ and no closed components properly and smoothly embedded in $B^4$. Then
$$|\sigma(L)| + |\eta(L) - (\beta_0(F) - 1)| \leq \beta_1(F)$$
\end{thm}
Taking $F$ to be a connected surface realising $h_4(L)$ we deduce that 
$$|\sigma(L)| \leq h_4(L)$$ 

\begin{lemma}
\label{lemma: sqpdef background}
If $L$ is definite, then $g_{4}(L) = g(L)$.
\end{lemma}

\begin{proof}
If $L$ is definite then $|\sigma(L)| \leq h_4(L) \leq h(L) = |\sigma(L)|$ and therefore $0 = h(L) - h_{4}(L) = 2(g(L) - g_{4}(L))$.
\end{proof}

Here is an immediate consequence of the Murasugi-Tristram inequality (Theorem \ref{thm: mt ineq background}). 

\begin{lemma}
\label{lemma: sig is zero} 
If $\eta(L) = 0$ and $L$ bounds a smooth surface in $B^4$ consisting of annuli and a disk then $\sigma(L) = 0$.
\end{lemma}

\begin{prop}
\label{prop: same sigs} 
Let $L_0$ be a link and $L_1$ a non-empty sublink. Let $F$ be a smooth surface properly embedded in $S^3 \times I$ such that $F \cap (S^3 \times \{0\}) = L_0 \times \{0\}$, and $F$ consists of $L_1 \times I$ together with annuli having their boundaries in $S^3 \times \{0\}$. If $\eta(L_0) = \eta(L_1) = 0$ then $\sigma(L_0) = \sigma(L_1)$. 
\end{prop}

\begin{proof}
  
Consider the connected sum $L_0 \# -L_1$, formed along a component of $-L_1$ and the corresponding component of $L_1 \subset L_0$. From the surface $F$ we can construct a surface $F'$ in $S^3 \times I$ such that $F' \cap (S^3 \times \{0\}) = (L_0 \# -L_1) \times \{0\}$, and $F'$ is the disjoint union of $(L_1 \# -L_1) \times I$ and annuli with their boundaries in $S^3 \times \{0\}$. It is easy to see that $L_1 \# - L_1$ bounds a surface $F''$ in $B^4$ consisting of a disk and $|L_1| - 1$ annuli. Capping off $(S^3 \times \{1\}, (L_1 \# - L_1) \times \{1\})$ with $(B^4, F'')$, $F' \cup F''$ is a surface $F^*$ in a $4$-ball with $\partial F^* = L_0 \# - L_1$, consisting of annuli and a disk. Since $\eta(L_0 \# - L_1) = \eta(L_0) + \eta(L_1) = 0$, Lemma \ref{lemma: sig is zero}  gives 
$\sigma(L_0) - \sigma(L_1) = \sigma(L_0 \# - L_1) = 0$. 
\end{proof}

\begin{cor}
\label{cor: indef}
  Let $L_0$ and $L_1$ be links as in Proposition \ref{prop: same sigs} such that $\eta(L_1) = 0$ and $|\sigma(L_1)| < |L_0| - 1$. Then $L_0$ is indefinite.
\end{cor}

\begin{proof}

  If $\eta(L_0) \ne 0$ then $L_0$ is indefinite.

  If $\eta(L_0) = 0$ then $\sigma(L_0) = \sigma(L_1)$ by Proposition \ref{prop: same sigs}, giving $|\sigma(L_0)| = |\sigma(L_1)| < |L_0| - 1 \le h(L_0)$. So $L_0$ is again indefinite.
\end{proof}

\begin{proof}[Proof of Proposition \ref{prop: def g = 0}]

The Seifert links $L$ with $g(L) = 0$ are listed in Theorem \ref{thm: g = 0}. 

  If $L$ is a connected sum of Hopf links then clearly $L$ is definite if and only if the Hopf links all have the same sign. So suppose $L$ is not a connected sum of Hopf links.

  First consider $L(p,q;k,0)$. If $k = 2$ then $L$ bounds, in the Heegaard torus, an annulus whose core is $T(p,q)$, so it has Seifert matrix $[pq]$. Hence $L$ is definite. If $k \ge 4$ then $\Delta_L = 0$, so $\eta(L) \ne 0$, and therefore $L$ is indefinite.

  In the other cases with $w = 0$, we band together pairs of oppositely oriented $C_{p,q}$ components to get a cobordism in $S^3 \times I$ as in Proposition \ref{prop: same sigs} from $L$ to a link $L'$, where, if $L = L(1,q;k,0; \varepsilon)$, then $L' = C_1$ is the unknot, and if $L = L(p,p+1;k,0;-,+)$ then $L' = C_1 \cup C_2 = H_-$. In both cases $L$ is indefinite by Corollary \ref{cor: indef}. 

  In the cases where $w = 1$, in parts (2) and (3) of Theorem \ref{thm: g = 0} we have $k \ge 3$ (see Section 3.1). The above construction gives a similar cobordism from $L$ to $L_1$, where $L_1$ is the corresponding link with $k = 1$. If $L = L(1,1;k,1)$ then $L_1$ is the unknot, and if $L = L(1,2;k,1;-)$ then $L_1 = H_-$. In both cases $L$ is indefinite by Corollary \ref{cor: indef}. 

  Finally, suppose $L = L(2,3;k,1;-,-)$. If $k = 1$ then the framings on $C_1$ and $C_2$ are 1 and 2 respectively. So taking $[C_1], \; [C_2]$ as a basis for $H_1(P)$ we see that $L$ has Seifert matrix $\begin{bmatrix}
    1 & 1\\
    1 & 2\\
  \end{bmatrix}$. 

  Hence $\sigma(L) = 2$, so $L$ is definite. If $k \ge 3$ then the usual banding construction gives a cobordism as in Proposition \ref{prop: same sigs} from $L$ to $L_1 = L(2,3;1,1;-,-)$. Since $\eta(L_1) = 0$ and $\sigma(L_1) = 2 < k = |L| - 1$, $L$ is indefinite by Corollary \ref{cor: indef}.
\end{proof} 
  
\begin{proof}[Proof of Theorem \ref{thm: def results}]
If $g(L) = 0$ this is Proposition \ref{prop: def g = 0}. 

Also, the $ADE$  links are definite. On the other hand, if $L$ is definite and $g(L) \ne 0$ then $L \in \mathcal P$ by Lemma \ref{lemma: sqpdef background} and Theorem \ref{thm: g_4 = g results}. Hence $L$ is an $ADE$  link by \cite{Baa}.
\end{proof}

\section{The \texorpdfstring{$ADE$}{ADE} link conjecture} 
\label{sec: $ADE$ lk conj}
We first recall the $ADE$ Link Conjecture below.
\begin{customconj}{\ref{conj: $ADE$ conj}} 
If $L$ is a prime, fibred, strongly quasipositive link that is not an $ADE$ link then all cyclic branched covers $\Sigma_{\psi}(L)$ of $L$ are $NLS$, $LO$, and $CTF$.
\end{customconj}
 
The $ADE$ links are Seifert, and are not toroidal since they are prime and their 2-fold branched covers do not contain incompressible tori \cite{GLith}. Hence, after Theorem \ref{thm: $ADE$ for seifert links intro}, to complete the proof of the $ADE$ link conjecture it remains to show that if $L$ is a prime, fibred, strongly quasipositive link that is either toroidal or hyperbolic, then $\Sigma_n(L)$ is $NLS$, $LO$, and $CTF$ for all $n \ge 2$. In the toroidal case one expects that a much more general statement holds (see Conjecture \ref{conj: GL}). On the other hand, the conditions of being fibred and strongly quasipositive are both necessary for the Seifert and hyperbolic cases (Section \ref{subsec: necessity of hyp}). 

Let $L$ be a toroidal link and let $T$ be an essential torus in $X(L)$. Then $S^3 = X_1 \cup_T X_2,$ where  $T$ is the boundary of both $X_i$.  Let $L_i = L \cap X_i$, $i = 1,2$. Without loss of generality, we assume that $X_1$ is a solid torus; then $L_1$ is non-empty. If $X_2$ is a solid torus then $T$ is a Heegaard torus in $S^3$ and we say that $T$ is {\it unknotted}. Otherwise, $X_2$ is the exterior of a non-trivial knot $K$ and we say that $T$ is {\it knotted}. If $L_2$ is empty then $T$ is necessarily knotted and $L$ is a {\it satellite} link, with {\it companion} $K$. So a satellite link is a toroidal link, though the converse is not necessarily true. However, a toroidal knot must be a satellite knot.  

\begin{conj}
\label{conj: GL}
If $L$ is a prime toroidal link then all cyclic branched covers $\Sigma_\psi(L)$ of $L$ are $NLS$, $LO$, and $CTF$. 
\end{conj}

Conjecture \ref{conj: GL} had its origins in \cite[Conjecture 1.7]{GLid}, where the authors conjectured that the conclusion of Conjecture \ref{conj: GL} holds for all prime satellite knots when $n$ is sufficiently large. It was successively upgraded to all $n \ge 2$, then to all canonical cyclic branched covers $\Sigma_n(L)$ for toroidal links $L$, and finally to the present version.

\subsection{Results on the \texorpdfstring{$ADE$}{ADE} Link Conjecture and Conjecture \ref{conj: GL}}
\label{subsec: results $ADE$ conjecture}

First of all, in the direction of Conjecture \ref{conj: GL}, the following is proved in \cite{BGH1}. Although the result there is only stated for the canonical cyclic branched covers $\Sigma_n(L)$, it is easy to see that the proof applies to arbitrary cyclic branched covers $\Sigma_\psi(L)$. 

\begin{thm}[\cite{BGH1}]
\label{thm: prime sat}
If $L$ is an  prime satellite link then any cyclic branched cover $\Sigma_\psi(L)$ of $L$ is $NLS$ and $LO$. If $L$ has a fibred companion then $\Sigma_\psi(L)$ is $CTF$.
\end{thm}

\begin{cor}[The $ADE$ Link Conjecture for Satellite Links, \cite{BGH1}]
  \label{cor: $ADE$ sat}
The $ADE$ link conjecture is true for satellite links.
\end{cor}

Next, we consider $NLS$, $LO$, and $CTF$ versions of the $ADE$ conjecture and Conjecture \ref{conj: GL} separately.

\subsubsection{The \texorpdfstring{$LO$}{LO} version}
The proof of the $LO$ version of the $ADE$ link conjecture is almost complete: the only case remaining is stated in the following problem. 
\begin{problem}
Show that the $\Sigma_n(L)$ is $LO$ when $L$ is prime, fibred, strongly quasipositive, and toroidal, and every essential torus in the exterior of $L$ is unknotted.
\end{problem}

Here we provide more details. It turns out that the $LO$ version of Theorem \ref{thm: prime sat} holds more generally. See \cite[Theorem 2.12]{BGH1}.
\begin{thm}[\cite{BGH1}]
  If $L$ is an  prime link whose exterior contains a knotted essential torus then all cyclic branched covers $\Sigma_\psi(L)$ of $L$ are $LO$.
\end{thm}

Turning to the hyperbolic case, let $L = K_1 \cup \ldots \cup K_m$ be an  hyperbolic link. Gabai and Mosher have shown, independently, that there is a pseudo-Anosov flow $\Phi$ on $S^3 \setminus L$. (No proof has been published; but see \cite{Mo} and \cite{LT}.) For each $i,\, 1 \le i \le m$, $\Phi$ determines a slope with multiplicity on $\partial N(K_i)$, the {\it degeneracy locus} $\delta_i(\Phi)$ of $\Phi$ on $\partial N(K_i)$. We say that $\delta_i(\Phi)$ is {\it meridional} if the corresponding slope is a meridian of $K_i$. The following is proved in \cite{BGH2}.

\begin{thm}[\cite{BGH2}]
\label{thm: pA flow}
  Let $L$ be an  hyperbolic link whose complement admits a pseudo-Anosov flow none of whose degeneracy loci are meridional. Then all cyclic branched covers $\Sigma_\psi(L)$ of $L$ are $LO$.
\end{thm}

Fibred hyperbolic links with non-zero fractional Dehn twist coefficient on each boundary component of the fibre  satisfy the hypothesis of Theorem \ref{thm: pA flow}. This gives the following.

\begin{cor}[\cite{BGH2}]
\label{cor: $LO$ hyp FSQP}
  Let $L$ be a hyperbolic, fibred, strongly quasipositive link. Then all cyclic branched covers $\Sigma_\psi(L)$ of $L$ are $LO$.
\end{cor}

In particular we get

\begin{cor}[The LO version for Hyperbolic Links \cite{BGH2}]
\label{cor: $LO$ $ADE$ hyp}
  The $LO$ version of the $ADE$ link conjecture is true for hyperbolic links.
\end{cor}

Theorem \ref{thm: $ADE$ for seifert links intro}, Corollary \ref{cor: $ADE$ sat}, and Corollary \ref{cor: $LO$ $ADE$ hyp} combine to give the $LO$ version of the $ADE$ link conjecture for knots.

 \begin{thm}[LO version of the $ADE$ Knot Conjecture]
 \label{thm: $LO$ $ADE$ knots}
Let $K$ be a prime, fibred, strongly quasipositive knot. Then $\Sigma_n(K)$ is $LO$ for all $n \ge 2$ unless $K$ is $T(3,4)$, $T(3,5)$, or $T(2,q)$ for some odd $q \ge 3$.
\end{thm}

The cyclic branched covers of prime toroidal links contain essential tori, so their fundamental groups are infinite. Since finite groups are not left-orderable, our $LO$ result for hyperbolic links (Corollary \ref{cor: $LO$ $ADE$ hyp}) leads to the following interesting characterization of $ADE$ links.

\begin{cor}
The $ADE$ links can be characterised as the only fibred, strongly quasipositive links $L$ that have a cyclic branched cover $\Sigma_\psi(L)$ with finite fundamental group. 
\end{cor}

\subsubsection{The \texorpdfstring{$NLS$}{NLS} version}
Here our results only apply to the canonical cyclic branched covers $\Sigma_n(L)$.

First we have the following \cite[Theorem 1.1]{BBG1}. 

\begin{thm}[\cite{BBG1}]
  \label{thm: FSQP > 5}
  If $L$ is a fibred, strongly quasipositive link of positive genus then $\Sigma_n(L)$ is $NLS$ for all $n \ge 6$.
\end{thm}

Note that this is best possible since $\Sigma_n(T(2,3))$ has finite fundamental group for $2 \le n \le 5$. 
When $L$ is a toroidal link, in addition to Theorem \ref{thm: prime sat} we have
\begin{thm}[{\rm \cite[Theorem 2.13]{BGH1}}]
\label{thm: tor $NLS$ 2 3}
If $L$ is a prime toroidal link then $\Sigma_{2^a}(L)$ is $NLS$ for $a \geq 1$ and $\Sigma_{2^a3}(L)$ is $NLS$ for $a \geq 0$.
\end{thm}

Theorems \ref{thm: FSQP > 5} and \ref{thm: tor $NLS$ 2 3} give

\begin{thm}
  If $L$ is a prime, toroidal, fibred, strongly quasipositive link of positive genus then $\Sigma_n(L)$ is $NLS$ for $n \ne 5$.
\end{thm}

Thus for the $NLS$ version of the $ADE$ link conjecture for the canonical cyclic branched covers of toroidal links $L$ the only cases remaining are when $n = 5$ and when $g(L) = 0$. 

A link that is the closure of a positive braid is fibred and strongly quasipositive. Combining \cite{BBG1} and \cite{Baa} shows that the $NLS$ version of the $ADE$ link conjecture for canonical branched covers holds for such links.

\begin{thm}
\label{thm: $NLS$ PB}
  If $L$ is a prime positive braid link that is not an $ADE$ link then $\Sigma_n(L)$ is $NLS$ for all $n \ge 2$.
\end{thm}

\subsubsection{The \texorpdfstring{$CTF$}{CTF} version} When $L$ is a toroidal link, the best known result for the $CTF$ version is what has been stated in Theorem \ref{thm: prime sat}. In the hyperbolic case, for knots we have the following, which is contained in \cite[Corollary 1.3]{BH}.

\begin{thm}[\cite{BH}]
  \label{thm: $CTF$ 2(2g-1)}
  Let $K$ be a hyperbolic, fibred, strongly quasipositive knot whose monodromy $h$ has fractional Dehn twist coefficient $c(h)$. Then $\Sigma_n(K)$ is $CTF$ for $n \ge 1/|c(h)|$. In particular, this holds for $n \ge 2(2g(K) - 1)$.
\end{thm}

\subsubsection{Moore's question}
Finally, we recall the following question of Allison Moore, which was one of the motivations for \cite{BBG1}.

\begin{question} (Moore)
\label{question: moore}
  If $K$ is a hyperbolic $L$-space knot, is it true that $\Sigma_2(K)$ is not an $L$-space?
\end{question}

Since $L$-space knots are prime \cite{Krc}, fibred \cite{Ni}, and strongly quasipositive \cite{Hed}, Theorem \ref{thm: $LO$ $ADE$ knots} implies:

\begin{cor}[LO version of Moore's question]
\label{cor: $LO$ Moore}
If $K$ is an $L$-space knot then $\Sigma_n(K)$ is $LO$ for all $n \ge 2$ unless $K$ is $T(3,4)$, $T(3,5)$, or $T(2,q)$ for some odd $q \ge 3$.  
\end{cor}

In particular this shows that, modulo the $L$-space conjecture, the answer to Moore's question is affirmative. We also note that combining \cite[Corollary 1.4]{BBG1} with the main result of \cite{FRW} gives an affirmative answer to the analog of Moore's question for all higher order cyclic branched covers:

\begin{thm}
\label{thm: Moore > 2}
If $K$ is an $L$-space knot then $\Sigma_n(K)$ is $NLS$ for all $n \ge 3$ unless $K$ is $T(2,3)$ or $T(2,5)$.
\end{thm}

Despite Corollary \ref{cor: $LO$ Moore} and Theorem \ref{thm: Moore > 2}, Moore's original question remains open.

\subsection{Necessity of the hypotheses}
\label{subsec: necessity of hyp}
Although we expect that the cyclic branched covers of all prime toroidal links are $NLS$, $LO$, and $CTF$ (Conjecture \ref{conj: GL}), in the Seifert and hyperbolic cases of the $ADE$ link conjecture both the fibred and strongly quasipositive hypotheses are necessary.

For Seifert links, the links that are $ADE$ links up to orientation but are not $ADE$ links show this. First, consider the case $L = L(1,q;2,0)$, $q > 1$. Note that $L$ is obtained from the torus link $T(2,2q)$ by reversing the orientation of one of the components. $L$ is strongly quasipositive by Theorem \ref{thm: sqp seifert results intro} and $\Sigma_n(L)$ is an $L$-space for all $n$ \cite{Pe}. But $L$ is not fibred by Proposition \ref{prop: basic props}(3). It is also easy to see this directly since $L$  bounds an annulus with corresponding Seifert matrix $[q]$ where $q > 1$.

The other Seifert links that are $ADE$ up to orientation but not $ADE$ are $L(2,q;1,1;-)$, $q \ge 3$ odd, $L(1,q;2,2,-)$, $L(1,q;2,0;+)$ and $L(3,2;1,1;-)$. These are fibred, but not strongly quasipositive since they are not in $\mathcal P$.

In the hyperbolic case, to see that the fibred hypothesis is necessary we note that there are many hyperbolic links $L$ that are strongly quasipositive but not fibred and have $\Sigma_2(L)$ an $L$-space. The most obvious examples are alternating; then $\Sigma_2(L)$ is an $L$-space \cite{OSz2}, hence not $CTF$ \cite{OSz1}, and not $LO$ by \cite{BGW}. Also, an alternating link is strongly quasipositive if and only if it is special alternating if and only if it is definite \cite[Section 7]{BBG1}, and there are many such that are not fibred, even among the 2-bridge knots (for example $5_2$). For non-alternating examples, Dunfield has produced a list of the non-alternating knots with at most 13 crossings whose 2-fold branched covers are $L$-spaces (private communication), and of the 11- and 12-crossing examples in that list, KnotInfo gives the following as being strongly quasipositive but not fibred: $11n_k$ for $k$ = 93, 126, 136, 169, 171, 180, 181, and $12n_k$ for $k$ = 406, 453, 585, 806, and 881.

There are also hyperbolic examples showing that the hypothesis of strong quasipositivity is necessary, and in fact cannot be weakened to quasipositivity. Note that here we are restricted to non-alternating links, as there are no known examples of alternating links that are quasipositive but not strongly quasipositive; see \cite[Question 7.5]{BBG1}. Here, of the 11- and 12-crossing knots in Dunfield's list of non-alternating knots $K$ with $\Sigma_2(K)$ an $L$-space, KnotInfo identifies the following as being fibred, quasipositive, but not strongly quasipositive: $11n_k$ for $k$ = 95, 108, 109, 113, 118, 127, 144, 146, 159, 172, 176, and $12n_k$ for $k$ = 316, 375, 407, 409, 441, 451, 510, 514, 543, 621, 666, 683, 684, 707, 708, 717, 719, 730, 747, 748, 767, 768, 820, 822, 829, 831, 838, 871, 882, and 887.

\subsection{Behavior of the \texorpdfstring{$n$}{n}-fold cyclic branched covers of a link as \texorpdfstring{$n$}{n} varies}
Known results, including Theorem \ref{thm: canonical br covers}, suggest the possibility that for a given prime link $L$, its canonical cyclic branched covers $\Sigma_n(L)$ satisfy one of the following, where $\ast$ denotes any of the properties $LO$, $NLS$, or $CTF$.

\begin{enumerate}[leftmargin=*] 
\setlength\itemsep{0.4em}  
\item[{\rm (A)}] $\Sigma_n(L)$ is $\ast$ for all $n \ge 2$;
\item[{\rm (B)}] $\Sigma_n(L)$ is not $\ast$ for all $n \ge 2$;
\item[{\rm (C)}] for some integer $N$ such that $2 \le N \le 5$, $\Sigma_n(L)$ is $\left\{ \begin{array}{cc} 
    \mbox{not } \ast & \mbox{ for } 2 \le n \le N, \\ 
    \ast & \mbox{ for } n > N.  
    \end{array} \right.$
\end{enumerate}

A consequence of this would be

\begin{enumerate}[leftmargin=*] 
\setlength\itemsep{0.3em}  
\item[{\rm (D)}] if $\Sigma_2(L)$ is $\ast$ then $\Sigma_n(L)$ is $\ast$ for all $n \ge 2$.
\end{enumerate}

Since it is relatively rare for a 3-manifold to be an $L$-space, one expects that for most prime links $L$, $\Sigma_2(L)$ is $NLS$, $LO$, and $CTF$, and therefore that most links should satisfy (A). Indeed, we expect (A) to hold for all prime toroidal links (Conjecture \ref{conj: GL}). On the other hand, for Seifert and hyperbolic links all three cases occur. For Seifert links, see Theorem \ref{thm: canonical br covers}.

In the hyperbolic case all three cases (A), (B), and (C) occur even for knots. There are many examples illustrating this; we mention just three. The figure eight knot satisfies (B), for all three values of $\ast$. The $L$-space knot $P(-2,3,7)$ satisfies (A) for $\ast = $ $LO$ by Corollary \ref{cor: $LO$ Moore}, and for $\ast = $ $CTF$ by [Z, Proposition 1.13]. The latter implies that it satisfies (A) for $\ast$ = $NLS$; this was previously known by \cite[Corollary 1.2]{BBG1} and the fact that the Alexander polynomial of $P(-2,3,7)$ is not a product of cyclotomic polynomials \cite[Subsection 3.2]{Nie}. Finally, $9_3$, which is the 2-bridge knot with rational fraction 19/6, satisfies (C) for $\ast = $ $NLS$ or $LO$, with $N = 2$: $\Sigma_n(9_3)$ is an $L$-space and not $LO$ for $n = 2$, while for $n \ge 3$ it is $NLS$ by \cite[Corollary 10.11]{BBG1} and $LO$ by \cite[Theorem 1.3]{Tu}. It is the knot $K(2,3)$ in the notation of \cite[Subsection 10.4]{BBG1} and $J(3,-6)$ in the notation of \cite{Tu}.  

Considering arbitrary cyclic branched covers $\Sigma_\psi(L)$ of  links $L$, our results suggest that the following question might have a positive answer.

\begin{question}
\label{question: 2-fold}
  
Is it true that if $L$ is a prime  link such that $\Sigma_2(L)$ is $NLS$, $LO$, or $CTF$, then the same holds for all cyclic branched covers $\Sigma_\psi(L)$ of $L$?
\end{question}

Theorem \ref{thm: pA flow} shows that both the hypothesis and the conclusion hold for many hyperbolic links, and we expect that the same is true for all toroidal links (Conjecture \ref{conj: GL}). However, if $\Sigma_2(L)$ is not $\ast$, for some $\ast =$ $NLS$, $LO$, or $CTF$, the situation is more complicated, for both Seifert and hyperbolic links. This can be seen even with the canonical cyclic branched covers corresponding to different orientations of the same underlying unoriented link. 

For example, in the Seifert case Theorem \ref{thm: canonical br covers} shows that the $ADE$ link $L(1,q;2,2)$, $q > 1$, (which is the torus link $T(2,2q)$) satisfies (C) while $L(1,q;2,0)$ satisfies (B). For hyperbolic examples, let $L$ be the 2-bridge link with continued fraction $[2k,2l_1,2k,2l_2,...,2k,2l_r,2k]$, where $k \ge 3$, $l_i > 0$, and $r > 0$, with the canonical orientation, and let $L'$ be obtained from $L$ by reversing the orientation of one of the components. Then $\Sigma_2(L) = \Sigma_2(L')$ is a lens space. On the other hand, by \cite[Theorem 6.13 and Theorem 7.4]{BGH2}, $\Sigma_n(L)$ is an $L$-space (hence not $CTF$) and not $LO$ for all $n \ge 2$, while $\Sigma_n(L')$ is $NLS$ and $LO$ for all $n \ge 3$.  

In light of Question \ref{question: 2-fold}, it is of interest to determine which links have 2-fold branched covers that are $L$-spaces. Examples are the (successively strictly larger) families of alternating, quasi-alternating, and $\mathbb{Z}/2$-Khovanov thin links. The torus knots $T(3,4)$ and $T(3,5)$ have 2-fold branched covers that are not $\ast$ for all three values of $\ast$ but are not $\mathbb{Z}/2$-Khovanov thin. Also, Dunfield has shown (private communication) that the hyperbolic knot $K = 11n_{126}$ has $\Sigma_2(K)$ an $L$-space, while Jablan has shown that it is not $\mathbb{Z}/2$-Khovanov thin \cite{Jbn}; see \cite[Examples 11.1]{BBG1}. 

\begin{question}
  For which hyperbolic links $L$ is $\Sigma_2(L)$ an $L$-space?
\end{question}

\end{document}